# DIFFUSION APPROXIMATION FOR A PROCESSOR SHARING QUEUE IN HEAVY TRAFFIC[1]

By H. Christian Gromoll

*EURANDOM*

Consider a single server queue with renewal arrivals and i.i.d. service times in which the server operates under a processor sharing service discipline. To describe the evolution of this system, we use a measure valued process that keeps track of the residual service times of all jobs in the system at any given time. From this measure valued process, one can recover the traditional performance processes, including queue length and workload. We show that under mild assumptions, including standard heavy traffic assumptions, the (suitably rescaled) measure valued processes corresponding to a sequence of processor sharing queues converge in distribution to a measure valued diffusion process. The limiting process is characterized as the image under an appropriate lifting map, of a one-dimensional reflected Brownian motion. As an immediate consequence, one obtains a diffusion approximation for the queue length process of a processor sharing queue.

**1. Introduction.** Consider a queueing system which consists of a single server with an infinite capacity buffer, to which jobs arrive according to a delayed renewal process. The $i$th such arrival requires an amount of processing time that is the $i$th member of a sequence of independent and identically distributed strictly positive random variables. The server, rather than providing service to just one job at a time, operates under a processor sharing discipline; that is, it works simultaneously on all jobs currently in the buffer, providing an equal fraction of its attention to each. Thus, at any given time that the buffer is nonempty, each job in the buffer is being processed at a rate that is the reciprocal of the number of jobs in the buffer. When the

Received September 2002; revised June 2003.
[1]Supported in part by NSF Grants DMS-97-03891 and DMS-00-71408 and a gift from the David and Holly Mendel fund.

server has fulfilled a given job's service time requirement, the job exits the buffer. This system is known as a processor sharing queue.

The processor sharing service discipline can be viewed as an idealization of a round-robin or time-sharing protocol used in computer and communication systems. There is a considerable literature on processor sharing queues (see [18] for a survey up to 1987), much of which assumes either Poisson arrivals or exponential service times. For a discussion of more recent work, including a handful of results for the case of generally distributed interarrival and service times (the $GI/GI/1$ processor sharing queue), see [6].

In this paper, we present a heavy traffic diffusion approximation for a measure valued process that keeps track of the "state" of a $GI/GI/1$ processor sharing queue. A direct consequence of this is the existence of a diffusion approximation for the queue length process. Note that since the workload process in a $GI/GI/1$ queue is the same for all nonidling service disciplines, the heavy traffic approximation for the workload process under a processor sharing service discipline is the same as the well-known approximation under a FIFO (first-in-first-out) service discipline [9]. However, this simple relationship does not hold for the queue length process.

The measure valued process that we study keeps track of the *residual service times* of jobs in the buffer. The residual service time at time $t \geq 0$ of a job which has entered the buffer by time $t$, is given by the amount of processing time originally requested by the job minus the total amount of processing time it has received by time $t$. Jobs with residual service times at time $t$ equal to zero have received enough processing time to fulfill their requirement and have departed the buffer. Let $\mathcal{M}_{\mathrm{F}}$ denote the space of finite, nonnegative Borel measures on $\mathbb{R}_+ = [0, \infty)$. The measure valued process $\{\mu(t) : t \geq 0\}$ is such that for each $t \geq 0$, $\mu(t)$ is the random element of $\mathcal{M}_{\mathrm{F}}$ that has a unit of mass at the residual service time of each job currently in the buffer at time $t$. From this process, one can recover information about the performance of the system. For example, let $Z(t)$ denote the number of jobs in the buffer, or queue length, at time $t \geq 0$. Then $Z(t)$ can be recovered as the integral of $\mu(t)$ against the function that is identically one: $Z(t) = \langle 1, \mu(t) \rangle$, where $\langle 1, \mu(t) \rangle = \int_{\mathbb{R}_+} 1 \mu(t)(dx)$. Similarly, let $W(t)$ denote the sum of the residual service times of all jobs in the buffer, or workload, at time $t \geq 0$. Then $W(t) = \langle \chi, \mu(t) \rangle$, where $\langle \chi, \mu(t) \rangle = \int_{\mathbb{R}_+} \chi(x) \mu(t)(dx)$ and $\chi(x) = x$. Although the process $\mu(\cdot)$ includes information about the queue length and workload processes, it also provides a more detailed description of the state of the system than is available from these one-dimensional performance processes alone. It is this level of detail which facilitates our analysis. The process $\mu(\cdot)$ is called the *state descriptor* for the processor sharing queue. Note that this terminology is not intended to imply that $\mu(\cdot)$ is necessarily a Markovian state descriptor (it may not be, since we do not include the residual interarrival time, the time remaining until the next job arrives to the buffer, in the state description). The process $\mu(\cdot)$ has previously been used by Grishechkin [5], along with other measure valued descriptors, in his heavy traffic analysis of the steady state distribution of a processor sharing queue. It has also been used recently by Gromoll, Puha and Williams [6], Puha and Williams [16] and Puha, Stolyar and Williams [15] for obtaining fluid limit results. Similar measure valued descriptors have also been used recently to describe other queueing systems. Doytchinov, Lehoczky and Shreve [3] and Kruk, Lehoczky, Shreve and Yeung [11] have used measure valued processes in the context of a queueing system with deadlines. Limic [12, 13] has used measure valued processes in studying the heavy traffic behavior of a LIFO (last-in-first-out) preemptive resume queue.

In [6], a fluid (or law of large numbers) approximation for the state descriptor $\mu(\cdot)$ of a heavily loaded processor sharing queue was studied. It was conjectured there that an understanding of the steady state behavior of this fluid approximation would provide a key ingredient for establishing a heavy traffic diffusion approximation for the process $\mu(\cdot)$. This is verified in the present paper and the role played by fluid approximations is informally described in the following paragraphs. The reader is referred to Sections 3 and 4 for the complete treatment.

For each $r>0$ in some sequence which tends to infinity, the *fluid scaled* and *diffusion scaled* versions of $\mu(\cdot)$, denoted $\bar{\mu}^r(\cdot)$ and $\hat{\mu}^r(\cdot)$, are defined by $\bar{\mu}^r(t)=r^{-1}\mu(rt)$ and $\hat{\mu}^r(t)=r^{-1}\mu(r^2t)$ for $t\geq 0$. Since $\hat{\mu}^r(t)=\bar{\mu}^r(rt)$ for each $t\geq 0$, one can think of the diffusion scaled process $\hat{\mu}^r(\cdot)$ over a finite time interval $[0,T]$, where $T>1$, as corresponding to the fluid scaled process $\bar{\mu}^r(\cdot)$ over the time interval $[0,rT]$. The latter interval grows without bound as $r\to\infty$. To study $\bar{\mu}^r(\cdot)$ over this "long" time interval, one considers it in sections by covering $[0,rT]$ with the overlapping finite time intervals $[m,m+L]$, where $m=0,\ldots,\lfloor rT\rfloor$ and $L>1$ is fixed (here $\lfloor\cdot\rfloor$ denotes the integer part). The *shifted* fluid scaled processes $\bar{\mu}^{r,m}(\cdot)$ are then defined for each $m\leq\lfloor rT\rfloor$ and $t\in[0,L]$ by $\bar{\mu}^{r,m}(t)=\bar{\mu}^r(m+t)$. Thus, we can study the diffusion scaled process $\hat{\mu}^r(\cdot)$ over $[0,T]$ by studying the family of shifted fluid scaled processes $\{\bar{\mu}^{r,m}(\cdot): m\leq\lfloor rT\rfloor\}$ over $[0,L]$. By building on the techniques developed in [6], it is shown that as $r\to\infty$, "good" sample paths of the processes $\{\bar{\mu}^{r,m}(\cdot), m\leq\lfloor rT\rfloor\}$ can be uniformly approximated on $[0,L]$ by measure valued functions $\bar{\zeta}(\cdot):[0,\infty)\longrightarrow\mathcal{M}_{\mathrm{F}}$ known as *fluid model solutions* (see Definition 4.2 in Section 4). By "good" sample paths, we mean sample paths in a set whose probability approaches 1 as $r\to\infty$.

Recently, Puha and Williams [16] have shown that under mild conditions, such fluid model solutions $\bar{\zeta}(\cdot)$ converge to a steady state; that is, for any such $\bar{\zeta}(\cdot)$, there is a $\bar{\zeta}_\infty\in\mathcal{M}_{\mathrm{F}}$ such that as $t\to\infty$, $\bar{\zeta}(t)\longrightarrow\bar{\zeta}_\infty$ in the topology of weak convergence of measures. Moreover, $\bar{\zeta}_\infty$ depends only on the limiting service time distribution $\nu$, and on the first moment of $\bar{\zeta}(t)$, which remains constant in $t$. Indeed,

$$\bar{\zeta}_\infty = \Delta_\nu \langle \chi, \bar{\zeta}(t) \rangle \qquad \text{for all } t \geq 0,$$

where $\Delta_\nu : \mathbb{R}_+ \longrightarrow \mathcal{M}_{\mathrm{F}}$ is the lifting map defined by $\Delta_\nu w = w\langle \chi, \nu_e \rangle^{-1} \nu_e$ and $\nu_e$ denotes the residual excess lifetime distribution associated with $\nu$; that is, $\nu_e$ is the probability measure on $\mathbb{R}_+$ satisfying $\nu_e([0,x]) = \langle \chi, \nu \rangle^{-1} \int_0^x \nu(y,\infty)\, dy$ for all $x \in \mathbb{R}_+$. In [16], conditions are also identified under which convergence to steady state occurs *uniformly* for certain collections of fluid model solutions.

In the present paper, the steady state results in [16] are combined with the fluid approximations of $\{\bar{\mu}^{r,m}(\cdot) : m \leq \lfloor rT \rfloor\}$ described above to show that if $L > 1$ is sufficiently large and $t \in [L-1, L]$, then for "good" sample paths,

$$\bar{\mu}^{r,m}(t) \approx \Delta_\nu \langle \chi, \bar{\mu}^{r,m}(t) \rangle. \tag{1.1}$$

This fact is applied to the diffusion scaled state descriptors in the following way. For any $t \in [(L-1)r^{-1}, T]$, there is an $m \in \{0, \ldots, \lfloor rT \rfloor\}$ and an $s \in [L-1, L]$ such that

$$\hat{\mu}^r(t) = \bar{\mu}^r(rt) = \bar{\mu}^r(m+s) = \bar{\mu}^{r,m}(s). \tag{1.2}$$

Combining (1.2) with the approximation (1.1) yields that for "good" sample paths, and for $t \in [(L-1)r^{-1}, T]$,

$$\hat{\mu}^r(t) \approx \Delta_\nu \langle \chi, \hat{\mu}^r(t) \rangle. \tag{1.3}$$

It is shown that under mild conditions, including standard heavy traffic assumptions, the above approximation is valid on the entire interval $[0,T]$, with probability approaching 1 as $r \to \infty$. Recall that $\langle \chi, \hat{\mu}^r(t) \rangle = r^{-1} W(r^2 t)$ is the (diffusion scaled) workload at time $t$, which we denote by $\hat{W}^r(t)$. Thus, the essence of (1.3) is that as the diffusion scaled process $\hat{\mu}^r(\cdot)$ approaches the heavy traffic limit, it can be recovered from the diffusion scaled workload process $\hat{W}^r(\cdot)$ on $[0,T]$ by an appropriate lifting map. This phenomenon is known as *state space collapse*. Since it is known that under the assumptions we will impose, the process $\hat{W}^r(\cdot)$ converges in distribution to a reflected Brownian motion $W^*(\cdot)$ on $\mathbb{R}_+$, a diffusion approximation for $\hat{\mu}^r(\cdot)$ follows quickly from state space collapse. (See Sections 3 and 4 for more details.) This yields a limiting measure valued diffusion process which is confined to the one-dimensional subspace $\{c\nu_e : c \geq 0\}$ of $\mathcal{M}_{\mathrm{F}}$. It is given by

$$\mu^*(\cdot) = \Delta_\nu W^*(\cdot).$$

A consequence of the above diffusion approximation for the diffusion scaled state descriptor $\hat{\mu}^r(\cdot)$ is the existence of a diffusion approximation for the diffusion scaled queue length process $\hat{Z}^r(\cdot) = \langle 1, \hat{\mu}^r(\cdot)\rangle$. The limiting process is the one-dimensional reflected Brownian motion

$$Z^*(\cdot) = \langle 1, \mu^*(\cdot)\rangle.$$

The precise form of $Z^*(\cdot)$ is given in Corollary 2.4, which verifies a conjecture of Harrison and Williams [8] in the case of a single processor sharing queue with a single class of jobs. We note that in cases where the process $Z^*(\cdot)$ has negative drift, the steady state distribution of $Z^*(t)$ as $t \to \infty$ is consistent with that obtained by Grishechkin [5] using other means (see Section 2.4).

The method outlined above for obtaining a heavy traffic diffusion approximation using fluid approximations as an intermediate step is analogous to the method developed by Bramson [2] and Williams [17] for proving diffusion approximations for open multiclass networks with HL (head-of-the-line) service disciplines. Indeed, the framework developed in [2, 17] serves as the primary motivation for the approach taken here for the processor sharing discipline. Note that since processor sharing is not an HL service discipline, a large part of the machinery needs to be developed from first principles.

The paper is organized as follows. In Section 2 we give a precise description of our model and assumptions, and we state our main results. Section 3 is devoted to fluid scale analysis, which prepares us for the proof of state space collapse. The latter is given in Section 4.

1.1. *Notation.* The following notation will be used throughout the paper. Let $\mathbb{N} = \{1, 2, \dots\}$ and let $\mathbb{R}$ denote the set of real numbers. For $a, b \in \mathbb{R}$, we write $a \vee b$ for the maximum of $a$ and $b$, $a \wedge b$ for the minimum of $a$ and $b$, $a^+$ for the positive part of $a$, $\lfloor a \rfloor$ for the largest integer less than or equal to $a$, and $\lceil a \rceil$ for the smallest integer greater than or equal to $a$. The nonnegative real numbers $[0, \infty)$ will be denoted by $\mathbb{R}_+$. For a function $g : \mathbb{R}_+ \longrightarrow \mathbb{R}$, let $\|g\|_\infty = \sup_{x \in \mathbb{R}_+} |g(x)|$ and $\|g\|_K = \sup_{x \in [0,K]} |g(x)|$ for each $K \geq 0$.

For a set $B \subset \mathbb{R}_+$, we denote the indicator of the set $B$ by $\mathbb{1}_B$. We also define the following real valued functions on $\mathbb{R}_+ : \chi(x) = x$, for $x \in \mathbb{R}_+$, and $\varphi(x) = 1/x$, for $x \in (0, \infty)$ with $\varphi(0) = 0$. For a topological space $A$, denote by $\mathbf{C}_b(A)$ the set of continuous, bounded, real-valued functions defined on $A$. In addition, for an interval $I \subset \mathbb{R}$, $\mathbf{C}_b^1(I)$ is the set of once continuously differentiable, real-valued functions defined on $I$ that together with their first derivatives are bounded on $I$.

Recall that $\mathcal{M}_{\mathrm{F}}$ is the set of finite, nonnegative Borel measures on $\mathbb{R}_+$. Consider $\zeta \in \mathcal{M}_{\mathrm{F}}$ and a Borel measurable function $g : \mathbb{R}_+ \longrightarrow \mathbb{R}$ which is integrable with respect to $\zeta$. We define $\langle g, \zeta\rangle = \int_{\mathbb{R}_+} g(x)\zeta(dx)$. Our equations will involve expressions of the form $\int_{[a,\infty)} g(x-a)\zeta(dx)$ for $a > 0$. To ease notation throughout, we write this as $\langle g(\cdot - a), \zeta\rangle$, making the convention that such a $g$ is always extended to be identically zero on $(-\infty, 0)$. The space $\mathcal{M}_{\mathrm{F}}$ is endowed with the topology of weak convergence of measures; that is, for $\zeta_n, \zeta \in \mathcal{M}_{\mathrm{F}}$, $n \in \mathbb{N}$, we have $\zeta_n \xrightarrow{\mathrm{w}} \zeta$ if and only if $\langle g, \zeta_n \rangle \longrightarrow \langle g, \zeta \rangle$ as $n \to \infty$ for all $g : \mathbb{R}_+ \longrightarrow \mathbb{R}$ that are bounded and continuous. With this topology, $\mathcal{M}_{\mathrm{F}}$ is a Polish space [14]. It will be convenient to deal with weak convergence on $\mathcal{M}_{\mathrm{F}}$ by means of a suitable metric which we now define. First, choose a countable set of nonnegative functions $\{g_k^{\mathcal{G}}\}_{k=1}^{\infty} \subset \mathbf{C}_b^1(\mathbb{R}_+)$ that are each bounded by 1, and such that $\{g_k^{\mathcal{G}}\}_{k=1}^{\infty}$ is convergence determining for $\mathcal{M}_{\mathrm{F}}$, that is, such that $\zeta_n \xrightarrow{\mathrm{w}} \zeta$ in $\mathcal{M}_{\mathrm{F}}$, as $n \to \infty$, if and only if for each $k \in \mathbb{N}$, $\langle g_k^{\mathcal{G}}, \zeta_n \rangle \longrightarrow \langle g_k^{\mathcal{G}}, \zeta \rangle$, as $n \to \infty$. We refer the reader to [4], Chapter 3, Section 4, for a possible construction. Note that the construction in the proof of Proposition 4.2 in [4], Chapter 3, Section 4, can be modified to require that $g_k^{\mathcal{G}} \in \mathbf{C}_b^1(\mathbb{R}_+)$ and $\|g_k^{\mathcal{G}}\|_\infty \le 1$ for $k \in \mathbb{N}$. Next, for each $k \in \mathbb{N}$, let $h_k^{\mathcal{G}} \in \mathbf{C}_b^1(\mathbb{R}_+)$ be a nonnegative function that is bounded by 1 and satisfies $h_k^{\mathcal{G}} \equiv 0$ on $[0, k-1]$, $h_k^{\mathcal{G}} \equiv 1$ on $[k, \infty)$ and $\|(h_k^{\mathcal{G}})'\|_\infty \le 2$. Letting $\mathcal{G} = \{g_k^{\mathcal{G}}\}_{k=1}^{\infty} \cup \{h_k^{\mathcal{G}}\}_{k=1}^{\infty}$, we define a metric on $\mathcal{M}_{\mathrm{F}}$ in terms of $\mathcal{G}$ as follows. For $\zeta_1, \zeta_2 \in \mathcal{M}_{\mathrm{F}}$, define

$$
\begin{aligned}
\mathbf{d}[\zeta_1, \zeta_2] = {} & \sum_{k=1}^{\infty} 2^{-k} (|\langle g_k^{\mathcal{G}}, \zeta_1 \rangle - \langle g_k^{\mathcal{G}}, \zeta_2 \rangle| \wedge 1) \\
& + \sup_{k \in \mathbb{N}} |\langle h_k^{\mathcal{G}}, \zeta_1 \rangle - \langle h_k^{\mathcal{G}}, \zeta_2 \rangle|.
\end{aligned} \tag{1.4}
$$

It is straightforward to see that $\mathbf{d}[\cdot,\cdot]$ is a complete metric on $\mathcal{M}_{\mathrm{F}}$ that induces the topology of weak convergence of measures. Note that although $\mathbf{d}[\cdot,\cdot]$ depends on the set $\mathcal{G}$, we do not indicate this explicitly in the notation. We denote the zero measure in $\mathcal{M}_{\mathrm{F}}$ by $\mathbf{0}$ and the measure in $\mathcal{M}_{\mathrm{F}}$ that puts one unit of mass at the point $x \in \mathbb{R}_+$ (or Dirac measure at $x$) by $\delta_x$. We will use $\delta_x^+$ to denote the following truncated Dirac measure:

$$
\delta_x^+ = \begin{cases} \delta_x, & x \in (0, \infty), \\ \mathbf{0}, & x = 0. \end{cases}
$$

We use "$\Longrightarrow$" to denote convergence in distribution of random elements of a metric space. For random elements $X, Y$ of a metric space, $X \sim Y$ will denote equivalence in distribution. Following Billingsley [1], we will use $\mathbf{P}$ and $\mathbf{E}$, respectively, to denote the probability measure and expectation operator associated with whatever space the relevant random element is defined on. All stochastic processes used in this paper will be assumed to have paths that are right continuous with finite left limits (r.c.l.l.). For $L \in (1, \infty)$ and a Polish space $\mathcal{S}$, we denote by $\mathbf{D}_L(\mathcal{S}) = \mathbf{D}([0, L], \mathcal{S})$ [resp. $\mathbf{D}_\infty(\mathcal{S}) = \mathbf{D}([0, \infty), \mathcal{S})$] the space of r.c.l.l. functions from $[0, L]$ (resp. $[0, \infty)$) into $\mathcal{S}$, and we endow this space with the usual Skorohod $J_1$-topology [4]. We use "$\xrightarrow{J_1}$" to denote convergence in this topology.

**2. The processor sharing queue.** In this section we give a precise description of the processor sharing queue, specify our assumptions and state our main result. A formal definition of the processor sharing queue, as considered in this paper, was previously given in [6]. The reader is referred there for a detailed discussion. However, since the diffusion approximation presented here requires stronger assumptions on the model than were made in [6], it is necessary to briefly review the definition. This is the subject of Section 2.1, where we introduce a sequence of processor sharing queueing models and the associated notation. Section 2.2 describes the scaling and time shifts we will apply to this sequence of models, and Section 2.3 specifies our asymptotic assumptions on the sequence. The statement of our main result appears in Section 2.4.

2.1. *A sequence of processor sharing queues.* We now specify a sequence of processor sharing queueing models indexed by $r \in \mathcal{R}$, where $\mathcal{R} \subset (0, \infty)$ is a sequence which increases to infinity. Each model in the sequence may be defined on a separate probability space and we use $\mathbf{P}^r$ and $\mathbf{E}^r$, respectively, for the probability and expectation operator on each of these spaces. The $r$th model in the sequence consists of the following: A server, which processes jobs from an infinite capacity buffer according to the processor sharing discipline, a pair of stochastic primitive processes $E^r(\cdot), V^r_\cdot$, which describe the arrival of work to the buffer, and a random initial condition, which specifies the state of the system at time 0. A measure valued process $\mu^r(\cdot)$ together with a set of descriptive equations describe the time evolution of the state of the system.

The *exogenous arrival process* $E^r(\cdot)$ is a rate $\alpha^r$ delayed renewal process associated with a sequence $\{u_i^r\}_{i=1}^\infty$ of finite nonnegative *interarrival times*. For $t \geq 0$, $E^r(t)$ represents the total number of jobs which have arrived to the buffer during the time interval $(0, t]$. We refer to the $i$th job to arrive to the buffer after time 0 simply as the $i$th job. In contrast, jobs which are already in the buffer at time 0 will always be referred to as *initial* jobs. The quantity $u_1^r$ is the arrival time of the first job, and $u_i^r$, $i \geq 2$, is the elapsed time between the arrival of the $(i-1)$st and $i$th jobs. Thus, for $i \geq 1$, $U_i^r = \sum_{j=1}^i u_j^r$ is the arrival time of the $i$th job. We define $U_0^r = 0$. So for $t \geq 0$,

$$E^r(t) = \sup\{i \geq 0 : U_i^r \leq t\}. \tag{2.1}$$

It is assumed that $\{u_i^r\}_{i=1}^\infty$ is a sequence of independent random variables and that $\{u_i^r\}_{i=2}^\infty$ is i.i.d. with mean $(\alpha^r)^{-1} \in (0, \infty)$ and standard deviation $a^r < \infty$. The first element $u_1^r$ of the sequence is assumed to be strictly positive with finite mean.

The *service process* $\{V_i^r, i = 1, 2, \dots\}$ records the cumulative amount of processing time required from the server by the first $i$ jobs. It is defined from a sequence $\{v_i^r\}_{i=1}^\infty$ of strictly positive *service times* by

$$V_i^r = \sum_{j=1}^{i} v_j^r \qquad \text{for } i = 1, 2, \dots. \tag{2.2}$$

The quantity $v_i^r$ represents the amount of processing time that the $i$th job requires from the server. It is assumed that $\{v_i^r\}_{i=1}^\infty$ is a sequence of strictly positive i.i.d. random variables with common distribution given by a Borel probability measure $\nu^r$ on $\mathbb{R}_+$. We assume $\nu^r$ has mean $(\beta^r)^{-1} \in (0, \infty)$ and standard deviation $b^r < \infty$. Finally, we assume that the service process $V_\cdot^r$ is independent of the exogenous arrival process $E^r(\cdot)$.

The *initial condition* specifies $Z^r(0)$, the number of initial jobs present in the buffer at time zero, as well as the service time requirement for each of these initial jobs. We assume $Z^r(0)$ is a nonnegative, integer valued random variable. The service times for initial jobs are taken to be the first $Z^r(0)$ elements of a sequence $\{\tilde{v}_j^r\}_{j=1}^\infty$ of strictly positive random variables. The random variables $Z^r(0)$ and $\{\tilde{v}_j^r\}_{j=1}^\infty$ are assumed to be independent of $\{u_i^r\}_{i=2}^\infty$ and $\{v_i^r\}_{i=1}^\infty$, but are not assumed to be independent of one another. A convenient way to express the initial condition is to define an initial random measure $\mu^r(0)$, which is the nonnegative, finite random Borel measure on $\mathbb{R}_+$ given by

$$\mu^r(0) = \sum_{j=1}^{Z^r(0)} \delta_{\tilde{v}_j^r}. \tag{2.3}$$

Recall that $\delta_x$ denotes the Dirac measure at $x \in \mathbb{R}_+$. Henceforth, $\mu^r(0)$ will be used as the initial condition. We assume that $\mu^r(0)$ satisfies

$$\mathbf{E}^r[\langle 1, \mu^r(0)\rangle] < \infty, \tag{2.4}$$

$$\mathbf{E}^r[\langle \chi, \mu^r(0)\rangle] < \infty. \tag{2.5}$$

Note that since $\langle 1, \mu^r(0)\rangle = Z^r(0)$ and $\langle \chi, \mu^r(0)\rangle = \sum_{j=1}^{Z^r(0)} \tilde{v}_j^r$, assumptions (2.4) and (2.5) mean, respectively, that the expected initial queue length and expected initial workload are finite.

The *state descriptor* is a performance process which describes the time evolution of the state of the processor sharing queue. It is a measure valued process $\mu^r(\cdot) : [0, \infty) \longrightarrow \mathcal{M}_{\mathrm{F}}$ such that for each $t \geq 0$, the random Borel measure $\mu^r(t)$ has one unit of mass located at the *residual service time* of each job that is in the buffer at time $t$. The residual service time at time $t \geq 0$ of job $i \leq E^r(t)$ [resp., of initial job $j \leq Z^r(0)$], denoted $R_i^r(t)$ [resp., $\tilde{R}_j^r(t)$], is the remaining amount of processing time required to fulfill the service time requirement of the job. If this residual service time is zero, the job has completed service and has departed the buffer. Thus at $t \geq 0$, the state descriptor is given by

$$\mu^r(t) = \sum_{j=1}^{Z^r(0)} \delta^+_{\tilde{R}^r_j(t)} + \sum_{i=1}^{E^r(t)} \delta^+_{R^r_i(t)}. \tag{2.6}$$

Recall that $\delta^+_x$ is the Dirac measure at $x$ for $x \in (0,\infty)$ with $\delta^+_0 = \mathbf{0}$, which ensures in the above that $\mu^r(t)$ has one unit of mass only for each job with nonzero residual service time at time $t \geq 0$. Let $Z^r(t)$ denote the number of jobs in the buffer, or queue length, at time $t \geq 0$. Then clearly, $Z^r(t) = \langle 1, \mu^r(t)\rangle$ for all $t \geq 0$. Thus, under the processor sharing discipline, any job that is present in the buffer at time $t$ receives service at the instantaneous rate $\langle 1, \mu^r(t)\rangle^{-1}$. Note that if a job is present in the buffer at time $t$, then $\langle 1, \mu^r(t)\rangle \neq \mathbf{0}$. The *cumulative service per job* provided by the server up to time $t \geq 0$ is defined by

$$S^r(t) = \int_0^t \varphi(\langle 1, \mu^r(s)\rangle)\, ds, \tag{2.7}$$

where $\varphi(x) = 1/x$ for $x \in (0,\infty)$ with $\varphi(0) = 0$. Then, since job $i \leq E^r(t)$ arrives at time $U^r_i \leq t$, the cumulative amount of processing time that job $i$ receives by time $t$ is equal to $v^r_i \wedge (S^r(t) - S^r(U^r_i))$. Similarly, the cumulative amount of processing time that an initial job $j \leq Z^r(0)$ receives by time $t$ is equal to $\tilde{v}^r_j \wedge S^r(t)$. Therefore, for $t \geq 0$, the residual service times are given by the equations

$$R^r_i(t) = (v^r_i - S^r(t) + S^r(U^r_i))^+ \qquad \text{for } i = 1, \dots, E^r(t), \tag{2.8}$$

$$\tilde{R}^r_j(t) = (\tilde{v}^r_j - S^r(t))^+ \qquad \text{for } j = 1, \dots, Z^r(0). \tag{2.9}$$

Given the primitive processes $E^r(\cdot), V^r_\cdot$, and the initial condition $\mu^r(0)$, equations (2.6)–(2.9) determine the state descriptor $\mu^r(\cdot)$, the process $S^r(\cdot)$ and the residual service times. This fact is not difficult, although somewhat tedious, to show. We note that the definition of the state descriptor given here is equivalent to the formulation in [6].

Let $W^r(t)$ denote the (immediate) workload in the buffer at time $t \geq 0$. This is defined as the amount of time the server would have to work in order to complete the remaining service time requirements of all jobs in the buffer at time $t$, assuming no new arrivals take place. Since this is given by the sum of the residual service times of all jobs which are in the buffer at time $t$, we have

$$W^r(t) = \langle \chi, \mu^r(t)\rangle \qquad \text{for } t \geq 0. \tag{2.10}$$

The cumulative service per job process $S^r(\cdot)$ will play a particularly important role in our analysis. We will find it convenient to have notation for the increments of this process. For $t,h \geq 0$, define the *cumulative service per job in* $[t,t+h]$ by

$$S^r_{t,t+h} = S^r(t+h) - S^r(t) = \int_t^{t+h} \varphi(\langle 1, \mu^r(s)\rangle)\, ds. \tag{2.11}$$

Then for $U_i^r \leq t$, the amount of service received by the $i$th job by time $t$ can be written as $v_i^r \wedge S^r_{U_i^r,t}$, and the residual service time at time $t$ can be written

$$R_i^r(t) = (v_i^r - S^r_{U_i^r,t})^+.$$

2.2. *Scaling.* Our result concerns the asymptotic behavior of processor sharing queues on *diffusion* scale. In particular, the focus of our attention is the diffusion scaled state descriptor, which is defined for $t \geq 0$ by

$$\hat{\mu}^r(t) = \frac{1}{r}\mu^r(r^2 t). \tag{2.12}$$

We will also be interested in diffusion scaled versions of the workload and queue length processes, defined for $t \geq 0$ by $\hat{W}^r(t) = r^{-1}W^r(r^2t) = \langle \chi, \hat{\mu}^r(t)\rangle$ and $\hat{Z}^r(t) = r^{-1}Z^r(r^2t) = \langle 1, \hat{\mu}^r(t)\rangle$, respectively. Much of our understanding of the diffusion scaled state descriptor $\hat{\mu}^r(\cdot)$ will be derived from results about the *fluid* scaled state descriptor, which is defined for $t \geq 0$ by

$$\bar{\mu}^r(t) = \frac{1}{r}\mu^r(rt). \tag{2.13}$$

The relationship $\hat{\mu}^r(t) = \bar{\mu}^r(rt)$ is essential for bootstrapping results from fluid scale up to diffusion scale. Recall that our approach for studying the diffusion scaled process $\hat{\mu}^r(\cdot)$ over a fixed finite time interval $[0,T]$, for $T > 1$, is to look at the fluid scaled process $\bar{\mu}^r(\cdot)$ over the time interval $[0, rT]$. Specifically, we study overlapping sections of $\bar{\mu}^r(\cdot)$, each defined on a finite time interval of fixed length $L > 1$. For each $r \in \mathcal{R}$, $t \geq 0$ and $m \in \{0,1,\dots\}$, define

$$\bar{\mu}^{r,m}(t) = \bar{\mu}^r(m+t). \tag{2.14}$$

Then since for each $r \in \mathcal{R}$, the time interval $[0,rT]$ is covered by the $\lfloor rT \rfloor + 1$ overlapping time intervals $[m, m+L]$, for $0 \leq m \leq \lfloor rT \rfloor$, we have that for any $t \in [0, rT]$, there is (at least one) $m \leq \lfloor rT \rfloor$ and $s \in [0,L]$ such that

$$\bar{\mu}^r(t) = \bar{\mu}^{r,m}(s). \tag{2.15}$$

Much of Section 3 is devoted to extending results in [6] on the asymptotics (as $r \to \infty$) of $\bar{\mu}^r(\cdot)$ over $[0,L]$, so that they hold for $\bar{\mu}^{r,m}(\cdot)$ over $[0,L]$ for all $m \leq \lfloor rT \rfloor$. For this discussion, we define the following fluid scaled and shifted fluid scaled processes. For all $r \in \mathcal{R}$, $t \geq 0$, $h \geq 0$ and $m \in \{0, 1, \dots\}$, let

$$\bar{E}^r(t) = \frac{1}{r} E^r(rt), \qquad \bar{E}^{r,m}(t) = \bar{E}^r(m+t), \tag{2.16}$$

$$\bar{Z}^r(t) = \frac{1}{r} Z^r(rt), \qquad \bar{Z}^{r,m}(t) = \bar{Z}^r(m+t), \tag{2.17}$$

$$\bar{W}^r(t) = \frac{1}{r} W^r(rt), \qquad \bar{W}^{r,m}(t) = \bar{W}^r(m+t), \tag{2.18}$$

$$\bar{S}^r_{t,t+h} = S^r_{rt,r(t+h)}, \qquad \bar{S}^{r,m}_{t,t+h} = \bar{S}^r_{m+t,m+t+h}. \tag{2.19}$$

Note that by (2.19),

$$\begin{aligned} \bar{S}^{r,m}_{t,t+h} = S^r_{r(m+t),r(m+t+h)} &= \int_{r(m+t)}^{r(m+t+h)} \varphi(\langle 1, \mu^r(s)\rangle)\, ds \\ &= \int_t^{t+h} \varphi(\langle 1, \bar{\mu}^{r,m}(s)\rangle)\, ds. \end{aligned} \tag{2.20}$$

2.3. *Heavy traffic assumptions.* In this section, we specify the assumptions under which the diffusion approximation will be proved. Let $\alpha > 0$, $a > 0$ and $\theta > 0$ be fixed constants and let $\nu$ be a probability measure on $\mathbb{R}_+$ that does not charge the origin, and satisfies

$$\langle \chi^{4+\theta}, \nu\rangle < \infty. \tag{2.21}$$

Then $\beta = \langle \chi, \nu\rangle^{-1}$ is positive and finite and $b = (\langle \chi^2, \nu\rangle - \langle \chi, \nu\rangle^2)^{1/2}$ is finite. In order to obtain convergence in distribution of the diffusion scaled state descriptors $\hat{\mu}^r(\cdot)$ to a diffusion process, we impose the following asymptotic assumptions on the sequence of processor sharing queues defined in Section 2.1. For the sequence of arrival processes, we assume that as $r \to \infty$,

$$(\alpha^r, a^r) \longrightarrow (\alpha, a), \tag{2.22}$$

$$\mathbf{E}^r[u_1^r]/r \to 0, \tag{2.23}$$

$$\limsup_{r\to\infty} \mathbf{E}^r[(u_2^r)^{2+\theta}] < \infty. \tag{2.24}$$

For the sequence of service processes, we assume that as $r \to \infty$,

$$\nu^r \xrightarrow{\text{w}} \nu, \tag{2.25}$$

$$(\beta^r, b^r) \longrightarrow (\beta, b), \tag{2.26}$$

$$\limsup_{r\to\infty} \langle \chi^{4+\theta}, \nu^r\rangle < \infty. \tag{2.27}$$

Define the *traffic intensity parameter* for the $r$th system by $\rho^r = \alpha^r/\beta^r$, and let $\rho = \alpha/\beta$. In addition to the above, our assumption that the sequence of processor sharing queues approaches *heavy traffic* requires that

$$\rho = 1. \tag{2.28}$$

We also require the sequence to approach heavy traffic at an appropriate rate and assume that as $r \to \infty$,

$$r(1-\rho^r) \longrightarrow \lambda \qquad \text{for some } \lambda \in \mathbb{R}. \tag{2.29}$$

Assumption (2.23) implies that the initial residual interarrival time vanishes on diffusion scale. Assumptions (2.24) and (2.27) imply Lindeberg-type conditions, which along with (2.22), (2.23) and (2.26), are needed to imply functional central limit theorems for the triangular arrays $\{u_i^r; i = 1, 2, \dots\}_{r\in\mathcal{R}}$ and $\{v_i^r; i = 1, 2, \dots\}_{r\in\mathcal{R}}$ and ultimately, convergence in distribution of the diffusion scaled workload processes to a reflected Brownian motion (see Proposition 3.1). Assumption (2.27) is two moments stronger than what would normally be used for such functional central limit theorems. The additional restriction is used in a separate part of our analysis to estimate moments of the shifted fluid scaled state descriptors $\{\bar{\mu}^{r,m}(\cdot),\ m \le \lfloor rT \rfloor\}$ (see Lemma 3.6).

We also make assumptions on the asymptotic behavior of the diffusion scaled initial measure $\hat{\mu}^r(0)$. Since $\hat{\mu}^r(0) = \bar{\mu}^r(0)$, we make our assumptions in terms of $\bar{\mu}^r(0)$ since they will be used in that form. The following definition is central to these assumptions.

Definition 2.1 (Invariant manifold). Let $\nu_e$ denote the excess lifetime distribution associated with $\nu$, that is, $\nu_e$ is the probability measure on $\mathbb{R}_+$ satisfying $\langle 1\!\!1_{[0,x]}, \nu_e\rangle = \beta\int_0^x \langle 1\!\!1_{(y,\infty)}, \nu\rangle\, dy$ for all $x \in \mathbb{R}_+$. Let

$$\mathbf{M}_\nu = \{c\nu_e \in \mathcal{M}_{\mathrm{F}} : c \in \mathbb{R}_+\}.$$

Following usage in [2], we refer to the one parameter family of measures $\mathbf{M}_\nu \subset \mathcal{M}_{\mathrm{F}}$ as the invariant manifold associated with $\nu$.

Let $\Theta$ be a random measure taking values in $\mathcal{M}_{\mathrm{F}}$, such that

$$\Theta \in \mathbf{M}_\nu \qquad \text{a.s.,} \tag{2.30}$$

$$\mathbf{E}[\langle 1, \Theta\rangle] < \infty. \tag{2.31}$$

Note that assumption (2.21) implies that $\langle \chi, \nu_e\rangle < \infty$ and $\langle \chi^{k+\theta}, \nu_e\rangle < \infty$, for $k = 1, 2, 3$. Combining this with (2.30) and (2.31) implies that $\Theta$ also satisfies

$$\mathbf{E}[\langle \chi, \Theta\rangle] < \infty, \tag{2.32}$$

$$\mathbf{E}[\langle \chi^{k+\theta}, \Theta\rangle] < \infty \qquad \text{for } k = 1, 2, 3. \tag{2.33}$$

An important property of $\Theta$ that we will need is that it has no atoms a.s. This is a trivial consequence of (2.30) since $\nu_e$ has a density function. This "no atoms" property of $\Theta$ will be used in the following form:

$$\lim_{\kappa\downarrow 0}\mathbf{P}\left(\sup_{x\in\mathbb{R}_+}\langle \mathbb{1}_{[x,x+\kappa]},\Theta\rangle<\varepsilon\right)=1 \qquad \text{for all } \varepsilon>0. \tag{2.34}$$

The equivalence of (2.34) to the fact that $\Theta$ has no atoms a.s. is proved in [6], Lemma A.1.

For the sequence of (scaled) initial measures $\hat{\mu}^r(0)=\bar{\mu}^r(0)$, we assume that as $r\to\infty$,

$$(\bar{\mu}^r(0),\langle\chi,\bar{\mu}^r(0)\rangle,\langle\chi^{1+\theta},\bar{\mu}^r(0)\rangle) \quad\Longrightarrow\quad (\Theta,\langle\chi,\Theta\rangle,\langle\chi^{1+\theta},\Theta\rangle). \tag{2.35}$$

Note that for any $g\in\mathbf{C}_b(\mathbb{R}_+)$, $\Psi_g:\mathcal{M}_{\mathrm{F}}\longrightarrow\mathbb{R}$ defined by $\Psi_g(\zeta)=\langle g,\zeta\rangle$ is a continuous function. So the first component of (2.35) implies by the continuous mapping theorem ([1], Theorem 5.1) that for any such $g$, $\langle g,\bar{\mu}^r(0)\rangle\Longrightarrow\langle g,\Theta\rangle$ as $r\to\infty$. The second component of (2.35) implies that the fluid/diffusion scaled initial workload converges in distribution, that is, $\hat{W}^r(0)=\bar{W}^r(0)\Longrightarrow\langle\chi,\Theta\rangle$ as $r\to\infty$.

To simplify the statements of results for the remainder of the paper, we now summarize our assumptions in the following.

(Q.1) There is a sequence of processor sharing queues defined in Section 2.1 such that for some constants $\theta>0$, $\alpha>0$, $a>0$, some probability measure $\nu$ on $\mathbb{R}_+$ that does not charge the origin, and some random measure $\Theta$ taking values in $\mathcal{M}_{\mathrm{F}}$, (2.21)–(2.31) and (2.35) hold.

2.4. *Main result.* In order to state our main result, we will need the following definition.

DEFINITION 2.2. Assume (Q.1). Let $\Delta_\nu:\mathbb{R}_+\longrightarrow\mathcal{M}_{\mathrm{F}}$ be the lifting map associated with $\nu$ given by

$$\Delta_\nu w=\frac{w}{\langle\chi,\nu_e\rangle}\nu_e \qquad \text{for } w\in\mathbb{R}_+. \tag{2.36}$$

THEOREM 2.3. *Assume* (Q.1). *Then as $r\to\infty$, the sequence of diffusion scaled state descriptors $\{\hat{\mu}^r(\cdot)\}$ converges in distribution to the measure valued process $\mu^*(\cdot)=\Delta_\nu W^*(\cdot)$, where $W^*(\cdot)$ is a reflected Brownian motion on $\mathbb{R}_+$ with drift $-\lambda$, variance $\alpha a^2+\beta b^2$ and initial value $W^*(0)$ that is equal in distribution to $\langle\chi,\Theta\rangle$.*

COROLLARY 2.4. *Assume* (Q.1). *Then as $r\to\infty$, the sequence of diffusion scaled queue length processes $\{\hat{Z}^r(\cdot)\}$ converges in distribution to the process $Z^*(\cdot)=C_\nu W^*(\cdot)$, where $W^*(\cdot)$ is the reflected Brownian motion specified above, and $C_\nu=2\beta(1+\beta^2b^2)^{-1}$.*

PROOF. Since $\hat{Z}^r(\cdot) = \langle 1, \hat{\mu}^r(\cdot)\rangle$ and since $\Psi_1 : \mathcal{M}_{\mathrm{F}} \longrightarrow \mathbb{R}$ defined by $\Psi_1(\xi) = \langle 1, \xi\rangle$ is continuous, the result follows from Theorem 2.3 by the continuous mapping theorem. The form of the constant $C_\nu$ follows from (2.36) by rewriting the mean of $\nu_e$ in terms of the mean $\beta^{-1}$, and standard deviation $b$, of $\nu$. □

Harrison and Williams [8] have conjectured that for a queue operating under the processor sharing discipline (possibly with several different job classes which may arrive to the buffer), a heavy traffic limit theorem should hold for the queue length process associated with each job class. They predict that for each job class, the limiting queue length process should be a constant multiple of the reflected Brownian motion obtained as the limiting workload process for the queue. Moreover, they specify the form of the constants for which the above claim should hold (see A.58, A.60 and A.61 in [8]). Corollary 2.4 verifies this conjecture in the case of a single job class. In particular, the constant $C_\nu$ validates (after some rewriting) the prediction made in [8], A.60 and A.61.

In cases where $\lambda$ is positive, the process $Z^*(\cdot)$ has a steady state distribution which can be computed as in [7] to be

$$\begin{aligned}\lim_{t\to\infty} \mathbf{P}(Z^*(t) \le x) \\ = 1 - \exp(-\lambda(\beta^{-2} + b^2)(a^2 + b^2)^{-1} x).\end{aligned} \tag{2.37}$$

It is interesting to note that in the case $\beta = 1$, this is consistent with the steady state distribution obtained by Grishechkin [5]. Indeed, since the approach taken in [5] was to first take the limit as $t \to \infty$ to get the steady state distribution of $\hat{Z}^r(t)$ and then take the limit of the resulting steady state distributions as $r \to \infty$, the consistency with (2.37) suggests that an interchange of limits is possible in this situation.

The strategy for proving Theorem 2.3 can be divided into two main tasks. The first task is to establish a tightness property for the processes $\{\bar{\mu}^{r,m}(\cdot), r \in \mathcal{R}, m \le \lfloor rT \rfloor\}$. This is accomplished in Section 3. In Section 4, we undertake the second task, which is to identify limit points obtained in Section 3 as *fluid model solutions* (see Definition 4.2) and to show that, consequently, the process $\hat{\mu}^r(\cdot)$ can be approximated, as $r \to \infty$, by "overlapping" fluid model solutions. This will be done in such a way that as $r \to \infty$, the random measure $\hat{\mu}^r(t)$ is close to the steady state of some fluid model solution for all $t \in [0, T]$, with arbitrarily high probability. Combined with a result in [16] which provides rates of convergence to steady state for such fluid model solutions, this leads to the proof of state space collapse. Finally, the proof of Theorem 2.3 collecting these results appears at the end of Section 4.

**3. Compactness.** The aim of this section is to prove a tightness property for the shifted and fluid scaled state descriptors $\{\bar{\mu}^{r,m}(\cdot), r \in \mathcal{R}, 0 \le m \le \lfloor rT \rfloor\}$. Specifically, our goal is to show in Corollary 3.16 at the end of this section that sequences indexed by $r \in \mathcal{R}$ of "good" sample paths of the processes $\{\bar{\mu}^{r,m}(\cdot),$
$r \in \mathcal{R}, m \le \lfloor rT \rfloor\}$ are relatively compact in $\mathbf{D}_L(\mathcal{M}_{\mathrm{F}})$, for $L > 1$. Here, sample paths are "good" with arbitrarily high probability as $r \to \infty$. The road to Corollary 3.16 involves a detailed analysis of the behavior of the processor sharing queue on fluid scale and we rely heavily on the ideas developed in [6], Section 5. Nevertheless, many of the techniques developed there need to be refined and several new ideas are needed for the present analysis. This is to account for the fact that for each $r \in \mathcal{R}$, we are dealing with order $r$ measure valued processes (one for each $m \in \{0, \ldots, \lfloor rT \rfloor\}$) instead of just one.

The section is organized similarly to [6], Section 5. In Sections 3.1 and 3.2 we describe a dynamic equation satisfied by the fluid scaled state descriptors $\{\bar{\mu}^{r,m}(\cdot)\}$ and a well-known heavy traffic limit theorem for the workload processes $\{\langle \chi, \bar{\mu}^{r,m}(\cdot)\rangle\}$, respectively. Section 3.3 contains a generalization of Lemma 5.3 in [6]. It states that under certain conditions, the fluid scaled queue length process $\langle 1, \bar{\mu}^r(\cdot)\rangle$ can be bounded above over long time intervals. This will yield an upper bound for the shifted fluid scaled queue length process $\langle 1, \bar{\mu}^{r,m}(\cdot)\rangle$ on $[0, L]$ for each $m \le \lfloor rT \rfloor$. Section 3.4 contains a functional weak law of large numbers estimate, which is a refinement of Lemma A.2
in [6]. In Section 3.5, we obtain an upper bound for certain moments of the processes $\{\bar{\mu}^{r,m}(\cdot)\}$. In Section 3.6, we state and prove Lemma 3.8. This lemma combines the results of Sections 3.1–3.5 for subsequent easy reference. Lastly, Section 3.7 contains further analysis of the fluid scaled processes $\{\bar{\mu}^{r,m}(\cdot)\}$. It consists of four secondary lemmas which are consequences of the results summarized in Lemma 3.8. These four lemmas lead up to an oscillation bound (Theorem 3.14) for the processes $\{\bar{\mu}^{r,m}(\cdot)\}$ as well as Corollary 3.16.

3.1. *Dynamic equation.* We begin by specifying a dynamic equation satisfied by the processes $\{\bar{\mu}^{r,m}(\cdot)\}$ for $m \le \lfloor rT \rfloor$. Starting with (2.6) and substituting in the definition of the residual service times (2.8) and (2.9), one obtains, after some simplification, that for each $r$, a.s. for each Borel measurable function $g : \mathbb{R}_+ \longrightarrow \mathbb{R}$ and all $t, h \ge 0$,

$$\begin{aligned}\langle g, \mu^r(t+h)\rangle &= \langle (\mathbb{1}_{(0,\infty)} g)(\cdot - S^r_{t,t+h}), \mu^r(t)\rangle \\ &\quad + \sum_{i=E^r(t)+1}^{E^r(t+h)} (\mathbb{1}_{(0,\infty)} g)(v_i^r - S^r_{U_i^r, t+h}).\end{aligned} \tag{3.1}$$

Recall that we always assume $g$ is extended to be identically zero on $(-\infty,0)$ so that functions of the form $g(\cdot-a)$ are well defined on $\mathbb{R}_+$ for any $a>0$.

Equation (3.1) takes the following form for the shifted fluid scaled processes $\bar{\mu}^{r,m}(\cdot)$. For each $r\in\mathcal{R}$ and each $m\le\lfloor rT\rfloor$, a.s. for each Borel measurable function $g:\mathbb{R}_+\longrightarrow\mathbb{R}$, and all $t,h\ge 0$,

$$\begin{aligned}\langle g,\bar{\mu}^{r,m}(t+h)\rangle &= \langle(\mathbb{1}_{(0,\infty)}g)(\cdot-\bar{S}^{r,m}_{t,t+h}),\bar{\mu}^{r,m}(t)\rangle\\ &\quad+\frac{1}{r}\sum_{i=r\bar{E}^{r,m}(t)+1}^{r\bar{E}^{r,m}(t+h)}(\mathbb{1}_{(0,\infty)}g)(v_i^r-\bar{S}^{r,m}_{U_i^r r^{-1}-m,t+h}).\end{aligned}\tag{3.2}$$

We refer to (3.2) as the *dynamic equation* for $\bar{\mu}^{r,m}(\cdot)$. Frequently we will set $g\equiv 1$ in this equation, in which case it will look like

$$\begin{aligned}\langle 1,\bar{\mu}^{r,m}(t+h)\rangle &= \langle\mathbb{1}_{(\bar{S}^{r,m}_{t,t+h},\infty)},\bar{\mu}^{r,m}(t)\rangle\\ &\quad+\frac{1}{r}\sum_{i=r\bar{E}^{r,m}(t)+1}^{r\bar{E}^{r,m}(t+h)}\mathbb{1}_{(0,\infty)}(v_i^r-\bar{S}^{r,m}_{U_i^r r^{-1}-m,t+h}).\end{aligned}\tag{3.3}$$

An important bound which we will use often is obtained from (3.3) by ignoring any processing of jobs during the fluid scaled time interval $[t,t+h]$:

$$\langle 1,\bar{\mu}^{r,m}(t+h)\rangle\le\langle 1,\bar{\mu}^{r,m}(t)\rangle+\bar{E}^{r,m}(t+h)-\bar{E}^{r,m}(t).\tag{3.4}$$

3.2. *Diffusion limit for the workload process.* We will take advantage of the following well-known result.

PROPOSITION 3.1. *Assume* (Q.1). *Then as $r\to\infty$, the sequence of diffusion scaled workload processes $\{\hat{W}^r(\cdot)\}=\{\langle\chi,\hat{\mu}^r(\cdot)\rangle\}$ converges in distribution to a process $W^*(\cdot)$, which is a reflected Brownian motion on $\mathbb{R}_+$ with drift $-\lambda$, variance $\alpha a^2+\beta b^2$ and initial value $W^*(0)$ that is equal in distribution to $\langle\chi,\Theta\rangle$.*

PROOF. It is well known that assumptions (Q.1) are sufficient to imply the above result for the workload process of any single server, single class queue operating under a work conserving service discipline (including processor sharing). The use of functional central limit theorems and continuous mappings to prove such results goes back to [9]. For a detailed account, see, for example, [17], which covers a much more general setting than that considered here. □

COROLLARY 3.2. *Assume* (Q.1) *and let $T>1$ and $0<\eta'<1$ be given. Then:*

(i) *There exists* $M > 1$ *such that for any* $L > 1$,

$$\liminf_{r\to\infty} \mathbf{P}^r\left(\sup_{m\le\lfloor rT\rfloor} \|\langle\chi,\bar{\mu}^{r,m}(\cdot)\rangle\|_L \le M\right) \ge 1-\eta'. \tag{3.5}$$

(ii) *For any* $L > 1$ *and any* $\gamma > 0$,

$$\liminf_{r\to\infty} \mathbf{P}^r\left(\sup_{m\le\lfloor rT\rfloor,\ t\in[0,L]} |\langle\chi,\bar{\mu}^{r,m}(t)\rangle - \langle\chi,\bar{\mu}^{r,m}(0)\rangle| \le \frac{\gamma}{4}\right) \ge 1-\eta'. \tag{3.6}$$

PROOF. By expanding the definition of $\bar{\mu}^{r,m}(\cdot)$ and rewriting it in terms of diffusion scaling, we see that it suffices to show that:

(i) There exists $M > 1$ such that for any $L > 1$,

$$\liminf_{r\to\infty} \mathbf{P}^r\left(\sup_{t\in[0,T+(L/r)]} |\hat{W}^r(t)| \le M\right) \ge 1-\eta'. \tag{3.7}$$

(ii) For any $L > 1$ and any $\gamma > 0$,

$$\liminf_{r\to\infty} \mathbf{P}^r\left(\sup_{t\in[0,T],\ h\in[0,L/r]} |\hat{W}^r(t+h) - \hat{W}^r(t)| \le \frac{\gamma}{4}\right) \ge 1-\eta', \tag{3.8}$$

both of which follow easily from Proposition 3.1, using the tightness of $\{\hat{W}^r(\cdot)\}_{r\in\mathcal{R}}$ and the fact that the limiting process $W^*(\cdot)$ is a.s. continuous. Note that $M$ can indeed be chosen independent of $L$. □

3.3. *Stability of the queue length process.* The following lemma is a simple generalization of Lemma 5.3 in [6]. It will be used in two separate places to bound the total mass of the fluid scaled state descriptor over certain time intervals.

LEMMA 3.3. *Assume* (Q.1) *and let* $T > 1$ *be given. Suppose that for a given* $r \in \mathcal{R}$, *there is an event* $A^r$ *and there are constants* $t_0 \in [0, rT]$ *and* $t_1, c, l > 0$ *such that on* $A^r$,

(i) $\sup_{t\in[t_0,t_0+t_1]} \bar{E}^r(t+l) - \bar{E}^r(t) \le \frac{c}{4}$,
(ii) $\sup_{t\in[t_0,t_0+t_1]} \langle\chi,\bar{\mu}^r(t)\rangle \le \frac{l}{4}$,
(iii) $\langle 1,\bar{\mu}^r(t_0)\rangle \le \frac{c}{2}$.

*Then on* $A^r$,

$$\sup_{t\in[t_0,t_0+t_1]} \langle 1,\bar{\mu}^r(t)\rangle \le c. \tag{3.9}$$

Proof. Cover $[t_0, t_0+t_1]$ with the time intervals $I_k = [t_0+kl, t_0+(k+1)l]$, where $k \in \{0, \dots, \lfloor t_1/l \rfloor\}$. We prove by induction on $k$ that on $A^r$,

$$\sup_{t \in I_k} \langle 1, \bar{\mu}^r(t) \rangle \le c. \tag{3.10}$$

We first verify the case $k = 0$. We have on $A^r$, for $t \in I_0$,

$$\begin{aligned} \langle 1, \bar{\mu}^r(t) \rangle &\le \langle 1, \bar{\mu}^r(t_0) \rangle + \bar{E}^r(t) - \bar{E}^r(t_0) \\ &\le \langle 1, \bar{\mu}^r(t_0) \rangle + \bar{E}^r(t_0+l) - \bar{E}^r(t_0) \\ &\le \frac{c}{2} + \frac{c}{4} < c, \end{aligned} \tag{3.11}$$

by (3.4) (with $m = 0$), as well as (i) and (iii) above. We now proceed by induction and assume that (3.10) holds on $A^r$ for some $0 \le k < \lfloor t_1/l \rfloor$. To show that the statement holds with $(k+1)$ in place of $k$, we argue analogously to (3.11). We must first show that $\langle 1, \bar{\mu}^r(t_0+(k+1)l) \rangle \le c/2$. To this end, we use an argument inspired in part by an idea of Grishechkin ([5], page 542). The idea is to consider two cases separately: the case where the queue length becomes zero during the interval $I_k$, and the case where it does not. If $\langle 1, \bar{\mu}^r(s) \rangle$ is never zero for $s \in I_k$, we can write

$$\begin{aligned} \bar{S}^r_{t_0+kl, t_0+(k+1)l} &= \int_{t_0+kl}^{t_0+(k+1)l} \varphi(\langle 1, \bar{\mu}^r(s) \rangle)\, ds \\ &= \int_{t_0+kl}^{t_0+(k+1)l} \langle 1, \bar{\mu}^r(s) \rangle^{-1}\, ds \ge l/c, \end{aligned}$$

since we have assumed that (3.10) holds for this $k$. By Markov's inequality and (ii) above, this implies that

$$\begin{aligned} \Big\langle 1\!\!1_{(\bar{S}^r_{t_0+kl, t_0+(k+1)l}, \infty)}, \bar{\mu}^r(t_0+kl) \Big\rangle &\le \langle 1\!\!1_{(l/c, \infty)}, \bar{\mu}^r(t_0+kl) \rangle \\ &\le \frac{c}{l} \langle \chi, \bar{\mu}^r(t_0+kl) \rangle \le \frac{c}{4}. \end{aligned} \tag{3.12}$$

If, on the other hand, $\langle 1, \bar{\mu}^r(s) \rangle = 0$ for some $s \in I_k$, then all mass present in the system at time $t_0 + kl$ is gone by time $s$. More precisely, we have in this case, using $m = 0$, $t = t_0 + kl$ and $h = s - (t_0 + kl)$ in (3.3), that

$$\Big\langle 1\!\!1_{(\bar{S}^r_{t_0+kl, s}, \infty)}, \bar{\mu}^r(t_0+kl) \Big\rangle = 0.$$

Since $\bar{S}^r_{t_0+kl, s} \le \bar{S}^r_{t_0+kl, t_0+(k+1)l}$, we again obtain the estimate

$$\Big\langle 1\!\!1_{(\bar{S}^r_{t_0+kl, t_0+(k+1)l}, \infty)}, \bar{\mu}^r(t_0+kl) \Big\rangle = 0 < \frac{c}{4}. \tag{3.13}$$

Thus, we have on $A^r$,

$$\begin{aligned}\langle 1,\bar{\mu}^r(t_0+(k+1)l)\rangle &\le \left\langle \mathbb{1}_{(\bar{S}^r_{t_0+kl,t_0+(k+1)l},\infty)},\bar{\mu}^r(t_0+kl)\right\rangle \\ &\quad + \bar{E}^r(t_0+(k+1)l)-\bar{E}^r(t_0+kl) \\ &\le \frac{c}{4}+\frac{c}{4},\end{aligned} \tag{3.14}$$

where in the first inequality, we have used (3.3) and ignored any processing of jobs in the second term there. The second inequality follows by (3.12), (3.13) and (i) above.

Now we can complete the proof in a similar manner to that for (3.11). By (3.4), (3.14) and (i) above, we have on $A^r$ for any $t\in I_{k+1}$,

$$\begin{aligned}\langle 1,\bar{\mu}^r(t)\rangle &\le \langle 1,\bar{\mu}^r(t_0+(k+1)l)\rangle + \bar{E}^r(t)-\bar{E}^r(t_0+(k+1)l) \\ &\le \langle 1,\bar{\mu}^r(t_0+(k+1)l)\rangle + \bar{E}^r(t_0+(k+2)l)-\bar{E}^r(t_0+(k+1)l) \\ &\le \frac{c}{2}+\frac{c}{4}<c.\end{aligned}$$

□

3.4. *Weak law estimate.* The following lemma provides a generalized functional weak law of large numbers estimate, which will be used in several ways throughout the paper.

LEMMA 3.4. *Assume* (Q.1). *Let* $T,L>1$, $l>0$ *and* $0<\varepsilon',\ \eta'<1$. *Suppose that* $g:\mathbb{R}_+\longrightarrow\mathbb{R}_+$ *is a Borel measurable function such that* $\langle g,\nu\rangle<\infty$,

$$\langle g,\nu^r\rangle \longrightarrow \langle g,\nu\rangle \qquad \text{as } r\to\infty, \tag{3.15}$$

*and such that for some* $p'>0$, $\limsup_{r\to\infty}\langle g^{2+p'},\nu^r\rangle<\infty$. *Then*

$$\limsup_{r\to\infty}\mathbf{P}^r\left(\sup_{m\le\lfloor rT\rfloor}\left\|\frac{1}{r}\sum_{i=r\bar{E}^{r,m}(\cdot)+1}^{r\bar{E}^{r,m}(\cdot+l)} g(v_i^r)-\alpha l\langle g,\nu\rangle\right\|_L>\varepsilon'\right)\le\eta'. \tag{3.16}$$

PROOF. Define

$$c_g=\begin{cases}\langle g,\nu\rangle, & \langle g,\nu\rangle>0,\\ 1, & \langle g,\nu\rangle=0.\end{cases} \tag{3.17}$$

Let $\varepsilon_g=\varepsilon'(2\alpha c_g)^{-1}\wedge(l/2)$ and let $l_1=l-\varepsilon_g$ and $l_2=l+\varepsilon_g$. Note that for any $t\in[0,L]$, we always have $[t,t+l]\supset[k_1\varepsilon_g,k_1\varepsilon_g+l_1]$ for some $k_1\in\{0,\dots,\lceil L/\varepsilon_g\rceil\}$, and $[t,t+l]\subset[k_2\varepsilon_g,k_2\varepsilon_g+l_2]$ for some $k_2\in\{0,\dots,\lceil L/\varepsilon_g\rceil\}$. For any $r\in\mathcal{R}$, $m\le\lfloor rT\rfloor$, and $t\in[0,L]$, consider the event

$$\left\{\left|\frac{1}{r}\sum_{i=r\bar{E}^{r,m}(t)+1}^{r\bar{E}^{r,m}(t+l)} g(v_i^r)-\alpha l\langle g,\nu\rangle\right|>\varepsilon'\right\}. \tag{3.18}$$

This event must be contained in the event

$$\left\{\frac{1}{r}\sum_{i=r\bar{E}^{r,m}(t)+1}^{r\bar{E}^{r,m}(t+l)} g(v_i^r) > \alpha l\langle g,\nu\rangle + \varepsilon'\right\} \cup \left\{\frac{1}{r}\sum_{i=r\bar{E}^{r,m}(t)+1}^{r\bar{E}^{r,m}(t+l)} g(v_i^r) < \alpha l\langle g,\nu\rangle - \varepsilon'\right\}. \tag{3.19}$$

Thus, for some $k_1, k_2 \in \{0, \ldots, \lceil L/\varepsilon_g\rceil\}$, the event in (3.18) is contained in

$$\left\{\frac{1}{r}\sum_{i=r\bar{E}^{r,m}(k_2\varepsilon_g)+1}^{r\bar{E}^{r,m}(k_2\varepsilon_g+l_2)} g(v_i^r) > \alpha l_2\langle g,\nu\rangle + \frac{\varepsilon'}{2}\right\} \cup \left\{\frac{1}{r}\sum_{i=r\bar{E}^{r,m}(k_1\varepsilon_g)+1}^{r\bar{E}^{r,m}(k_1\varepsilon_g+l_1)} g(v_i^r) < \alpha l_1\langle g,\nu\rangle - \frac{\varepsilon'}{2}\right\}, \tag{3.20}$$

where we have used the definitions of $\varepsilon_g$, $l_1$, $l_2$ and $c_g$. Therefore, to prove (3.16), it suffices to show that for $j \in \{1,2\}$,

$$\limsup_{r\to\infty} \mathbf{P}^r\left(\sup_{m\le\lfloor rT\rfloor,\ k\le\lceil L/\varepsilon_g\rceil}\left|\frac{1}{r}\sum_{i=r\bar{E}^{r,m}(k\varepsilon_g)+1}^{r\bar{E}^{r,m}(k\varepsilon_g+l_j)} g(v_i^r) - \alpha l_j\langle g,\nu\rangle\right| > \frac{\varepsilon'}{2}\right) \le \frac{\eta'}{2}. \tag{3.21}$$

Furthermore, by summing over $m \le \lfloor rT\rfloor$ and $k \le \lceil L/\varepsilon_g\rceil$, it suffices to show for $j \in \{1,2\}$ that for all sufficiently large $r$, all $m \le \lfloor rT\rfloor$ and all $k \le \lceil L/\varepsilon_g\rceil$,

$$\mathbf{P}^r\left(\left|\frac{1}{r}\sum_{i=r\bar{E}^{r,m}(k\varepsilon_g)+1}^{r\bar{E}^{r,m}(k\varepsilon_g+l_j)} g(v_i^r) - \alpha l_j\langle g,\nu\rangle\right| > \frac{\varepsilon'}{2}\right) \le \frac{\eta'}{2(\lceil L/\varepsilon_g\rceil+1)(\lfloor rT\rfloor+1)}. \tag{3.22}$$

To this end, we now fix $j \in \{1,2\}$ and $k \in \{0, \ldots, \lceil L/\varepsilon_g\rceil\}$ for the remainder of the proof. For better readability, we will suppress notation indicating dependence on $k$ and $j$ for certain objects, when it is clear from the definition.

For each $r \in \mathcal{R}$ and $m \le \lfloor rT\rfloor$, let

$$A^{r,m} = \left\{\left|\frac{1}{r}\sum_{i=r\bar{E}^{r,m}(k\varepsilon_g)+1}^{r\bar{E}^{r,m}(k\varepsilon_g+l_j)} g(v_i^r) - \alpha l_j\langle g,\nu\rangle\right| > \frac{\varepsilon'}{2}\right\} \tag{3.23}$$

be the event appearing in (3.22). For notational convenience, we define $E^{r,m} = r\bar{E}^{r,m}(k\varepsilon_g)$ and $N^{r,m} = r\bar{E}^{r,m}(k\varepsilon_g + l_j) - r\bar{E}^{r,m}(k\varepsilon_g)$. Define $c_+ = 1 + \varepsilon'(8\alpha l_j c_g)^{-1}$ and $c_- = 1 - \varepsilon'(16(\alpha l_j c_g \vee 1))^{-1}$, and let

$$A_E^{r,m} = \{\lfloor c_- r\alpha l_j \rfloor \le N^{r,m} \le \lfloor c_+ r\alpha l_j \rfloor\}.$$

Let $\breve{A}_E^{r,m}$ denotes the complement of $A_E^{r,m}$. For each $r \in \mathcal{R}$ and $m \le \lfloor rT \rfloor$,

$$\mathbf{P}^r(A^{r,m}) \le \mathbf{P}^r(\breve{A}_E^{r,m}) + \mathbf{P}^r(A^{r,m} \cap A_E^{r,m}). \tag{3.24}$$

We will show that each of the two terms in (3.24) is bounded above by $\eta'[4(\lceil L/\varepsilon_g \rceil + 1)(\lfloor rT \rfloor + 1)]^{-1}$ for all sufficiently large $r$.

For the first term in (3.24), note that on $\breve{A}_E^{r,m}$, either $E^r(r(m + k\varepsilon_g + l_j)) - E^r(r(m + k\varepsilon_g)) < \lfloor c_- r\alpha l_j \rfloor$ or $E^r(r(m + k\varepsilon_g + l_j)) - E^r(r(m + k\varepsilon_g)) > \lfloor c_+ r\alpha l_j \rfloor$ by definition. Inversion of the renewal process $E^r(\cdot)$ implies that on $\breve{A}_E^{r,m}$, either

$$\sum_{i=E^{r,m}+1}^{E^{r,m}+\lfloor c_- r\alpha l_j \rfloor} u_i^r > r(m + k\varepsilon_g + l_j) - r(m + k\varepsilon_g) = rl_j \tag{3.25}$$

or

$$\sum_{i=E^{r,m}+2}^{E^{r,m}+\lfloor c_+ r\alpha l_j \rfloor} u_i^r \le r(m + k\varepsilon_g + l_j) - r(m + k\varepsilon_g) = rl_j. \tag{3.26}$$

Since $\alpha^r \longrightarrow \alpha$, as $r \to \infty$ [see (2.22)], and since by definition, $0 < c_- < 1 < c_+$, (3.25) and (3.26) imply that on $\breve{A}_E^{r,m}$, for some $c_2 > 0$ and all sufficiently large $r$, either

$$\sum_{i=E^{r,m}+1}^{E^{r,m}+\lfloor c_- r\alpha l_j \rfloor} \left(u_i^r - \frac{1}{\alpha^r}\right) > rl_j - \frac{\lfloor c_- r\alpha l_j \rfloor}{\alpha^r} > rc_2 \tag{3.27}$$

or

$$\sum_{i=E^{r,m}+2}^{E^{r,m}+\lfloor c_+ r\alpha l_j \rfloor} \left(u_i^r - \frac{1}{\alpha^r}\right) \le rl_j - \frac{\lfloor c_+ r\alpha l_j \rfloor - 1}{\alpha^r} \le -rc_2. \tag{3.28}$$

Define the random variables

$$x_i^r = u_i^r \mathbb{1}_{\{u_i^r \le r^{1/8}\}} \quad \text{and} \quad y_i^r = u_i^r \mathbb{1}_{\{u_i^r > r^{1/8}\}}, \tag{3.29}$$

and let

$$X_-^{r,m} = \sum_{i=E^{r,m}+1}^{E^{r,m}+\lfloor c_- r\alpha l_j \rfloor} (x_i^r - \mathbf{E}^r[x_i^r]), \tag{3.30}$$

$$Y_-^{r,m} = \sum_{i=E^{r,m}+1}^{E^{r,m}+\lfloor c_- r\alpha l_j \rfloor} (y_i^r - \mathbf{E}^r[y_i^r]), \tag{3.31}$$

$$X_{+}^{r,m}=\sum_{i=E^{r,m}+2}^{E^{r,m}+\lfloor c_{+}r\alpha l_j\rfloor}(x_i^r-\mathbf{E}^r[x_i^r]), \tag{3.32}$$

$$Y_{+}^{r,m}=\sum_{i=E^{r,m}+2}^{E^{r,m}+\lfloor c_{+}r\alpha l_j\rfloor}(y_i^r-\mathbf{E}^r[y_i^r]). \tag{3.33}$$

Then using (3.27)–(3.33), we have

$$\begin{aligned}
\mathbf{P}^r(\breve{A}_E^{r,m}) &\le \mathbf{P}^r(|X_{-}^{r,m}+Y_{-}^{r,m}|\ge rc_2)+\mathbf{P}^r(|X_{+}^{r,m}+Y_{+}^{r,m}|\ge rc_2)\\
&\le \mathbf{P}^r(|X_{-}^{r,m}|\ge rc_2/2)+\mathbf{P}^r(|X_{+}^{r,m}|\ge rc_2/2)\\
&\quad+\mathbf{P}^r(|Y_{-}^{r,m}|\ge rc_2/2)+\mathbf{P}^r(|Y_{+}^{r,m}|\ge rc_2/2)\\
&\le (rc_2/2)^{-4}(\mathbf{E}^r[(X_{-}^{r,m})^4]+\mathbf{E}^r[(X_{+}^{r,m})^4])\\
&\quad+(rc_2/2)^{-2}(\mathbf{E}^r[(Y_{-}^{r,m})^2]+\mathbf{E}^r[(Y_{+}^{r,m})^2])\\
&\le (rc_2/2)^{-4}(\lfloor c_{-}r\alpha l_j\rfloor^2+(\lfloor c_{+}r\alpha l_j\rfloor-1)^2)(r^{1/8})^4\\
&\quad+(rc_2/2)^{-2}(\lfloor c_{-}r\alpha l_j\rfloor+(\lfloor c_{+}r\alpha l_j\rfloor-1))\mathbf{E}^r[(y_2^r)^2] \qquad (3.34)\\
&\le r^{-3/2}(c_2/2)^{-4}2(c_{+}\alpha l_j)^2\\
&\quad+r^{-1}(c_2/2)^{-2}2c_{+}\alpha l_j\mathbf{E}^r[(u_2^r)^2;u_2^r>r^{1/8}]\\
&\le \frac{1}{\lfloor rT\rfloor+1}\biggl(\frac{\lfloor rT\rfloor+1}{r}\biggl(r^{-1/2}\frac{2(c_{+}\alpha l_j)^2}{(c_2/2)^4}\\
&\qquad\qquad+\frac{2c_{+}\alpha l_j}{(c_2/2)^2}r^{-\theta/8}\mathbf{E}^r[(u_2^r)^{2+\theta}]\biggr)\biggr).
\end{aligned}$$

The fourth inequality above follows from independence of the members of the sequence $\{x_i^r\}$ (resp. $\{y_i^r\}$). The fifth inequality follows since $c_{-}\le c_{+}$, and the last uses Markov's inequality in the final term. Letting $r\to\infty$ in (3.34), we see that the quantity in the outer parentheses on the last line tends to zero uniformly in $m$ [see (2.24)], and is thus smaller than $\eta'[4(\lceil L/\varepsilon_g\rceil+1)]^{-1}$ for all sufficiently large $r$. So $\mathbf{P}^r(\breve{A}_E^{r,m})$ is bounded as desired for all sufficiently large $r$ and all $m\le\lfloor rT\rfloor$.

We handle the second term in (3.24) in a similar fashion. Write $A^{r,m}=A_{-}^{r,m}\cup A_{+}^{r,m}$, where

$$A_{-}^{r,m}=\left\{\sum_{i=E^{r,m}+1}^{E^{r,m}+N^{r,m}}g(v_i^r)-r\alpha l_j\langle g,\nu\rangle<-\frac{r\varepsilon'}{2}\right\} \tag{3.35}$$

and

$$A_{+}^{r,m}=\left\{\sum_{i=E^{r,m}+1}^{E^{r,m}+N^{r,m}}g(v_i^r)-r\alpha l_j\langle g,\nu\rangle>\frac{r\varepsilon'}{2}\right\}. \tag{3.36}$$

Define the random variables

$$g_i^r = g(v_i^r)\mathbb{1}_{\{g(v_i^r)\le r^{1/8}\}} \quad \text{and} \quad h_i^r = g(v_i^r)\mathbb{1}_{\{g(v_i^r)>r^{1/8}\}} \tag{3.37}$$

and let

$$G_-^{r,m} = \sum_{i=E^{r,m}+1}^{E^{r,m}+\lfloor c_- r\alpha l_j\rfloor} (g_i^r - \mathbf{E}^r[g_i^r]), \tag{3.38}$$

$$H_-^{r,m} = \sum_{i=E^{r,m}+1}^{E^{r,m}+\lfloor c_- r\alpha l_j\rfloor} (h_i^r - \mathbf{E}^r[h_i^r]), \tag{3.39}$$

$$G_+^{r,m} = \sum_{i=E^{r,m}+1}^{E^{r,m}+\lfloor c_+ r\alpha l_j\rfloor} (g_i^r - \mathbf{E}^r[g_i^r]), \tag{3.40}$$

$$H_+^{r,m} = \sum_{i=E^{r,m}+1}^{E^{r,m}+\lfloor c_+ r\alpha l_j\rfloor} (h_i^r - \mathbf{E}^r[h_i^r]). \tag{3.41}$$

Let $c_3 = \varepsilon'(8c_+\alpha l_j)^{-1}$ and note that by assumption,

$$|\langle g,\nu^r\rangle - \langle g,\nu\rangle| \le c_3, \tag{3.42}$$

for all sufficiently large $r$. So for all sufficiently large $r$, and all $m \le \lfloor rT\rfloor$,

$$\begin{aligned}
&\mathbf{P}^r(A^{r,m}\cap A_E^{r,m})\\
&\quad= \mathbf{P}^r(A_-^{r,m}\cap A_E^{r,m}) + \mathbf{P}^r(A_+^{r,m}\cap A_E^{r,m})\\
&\quad\le \mathbf{P}^r(G_-^{r,m} + H_-^{r,m} + \lfloor c_- r\alpha l_j\rfloor\langle g,\nu^r\rangle - r\alpha l_j\langle g,\nu\rangle < -r\varepsilon'/2)\\
&\qquad + \mathbf{P}^r(G_+^{r,m} + H_+^{r,m} + \lfloor c_+ r\alpha l_j\rfloor\langle g,\nu^r\rangle - r\alpha l_j\langle g,\nu\rangle > r\varepsilon'/2)\\
&\quad\le \mathbf{P}^r(|G_-^{r,m} + H_-^{r,m}| + |\lfloor c_- r\alpha l_j\rfloor(\langle g,\nu^r\rangle - \langle g,\nu\rangle)\\
&\qquad\qquad + \langle g,\nu\rangle(\lfloor c_- r\alpha l_j\rfloor - r\alpha l_j)| > r\varepsilon'/2)\\
&\qquad + \mathbf{P}^r(|G_+^{r,m} + H_+^{r,m}| + |\lfloor c_+ r\alpha l_j\rfloor(\langle g,\nu^r\rangle - \langle g,\nu\rangle)\\
&\qquad\qquad + \langle g,\nu\rangle(\lfloor c_+ r\alpha l_j\rfloor - r\alpha l_j)| > r\varepsilon'/2)\\
&\quad\le \mathbf{P}^r(|G_-^{r,m} + H_-^{r,m}| + c_+ r\alpha l_j c_3 + \langle g,\nu\rangle(|c_- - 1|r\alpha l_j + 1) > r\varepsilon'/2)\\
&\qquad + \mathbf{P}^r(|G_+^{r,m} + H_+^{r,m}| + c_+ r\alpha l_j c_3 + \langle g,\nu\rangle|c_+ - 1|r\alpha l_j > r\varepsilon'/2) \qquad (3.43)\\
&\quad\le \mathbf{P}^r(|G_-^{r,m} + H_-^{r,m}| > r\varepsilon'/4) + \mathbf{P}^r(|G_+^{r,m} + H_+^{r,m}| > r\varepsilon'/4)\\
&\quad\le \mathbf{P}^r(|G_-^{r,m}| > r\varepsilon'/8) + \mathbf{P}^r(|G_+^{r,m}| > r\varepsilon'/8)\\
&\qquad + \mathbf{P}^r(|H_-^{r,m}| > r\varepsilon'/8) + \mathbf{P}^r(|H_+^{r,m}| > r\varepsilon'/8)\\
&\quad\le (r\varepsilon'/8)^{-4}(\mathbf{E}^r[(G_-^{r,m})^4] + \mathbf{E}^r[(G_+^{r,m})^4])
\end{aligned}$$

$$
\begin{aligned}
&\quad + (r\varepsilon'/8)^{-2}(\mathbf{E}^r[(H_-^{r,m})^2] + \mathbf{E}^r[(H_+^{r,m})^2]) \\
&\le (r\varepsilon'/8)^{-4}(\lfloor c_- r\alpha l_j\rfloor^2 + \lfloor c_+ r\alpha l_j\rfloor^2)(r^{1/8})^4 \\
&\quad + (r\varepsilon'/8)^{-2}(\lfloor c_- r\alpha l_j\rfloor + \lfloor c_+ r\alpha l_j\rfloor)\mathbf{E}^r[(h_1^r)^2] \\
&\le \frac{1}{\lfloor rT\rfloor + 1}\left(\frac{\lfloor rT\rfloor + 1}{r}\left(r^{-1/2}\frac{2(c_+\alpha l_j)^2}{(\varepsilon'/8)^4} + \frac{2c_+\alpha l_j \mathbf{E}^r[(h_1^r)^2]}{(\varepsilon'/8)^2}\right)\right).
\end{aligned}
$$

In the first inequality above, we have used (3.38)–(3.41), as well as the fact that $\lfloor c_- r\alpha l_j\rfloor \le N^{r,m} \le \lfloor c_+ r\alpha l_j\rfloor$ on $A_E^{r,m}$. The third inequality uses (3.42) and the fact that $c_- < 1 < c_+$. The fourth inequality follows from the definitions of $c_3$, $c_-$, $c_+$ and $c_g$, where for the first term it is necessary that $r > 16c_g/\varepsilon'$. The seventh inequality follows from independence of the elements of the sequence $\{g_i^r\}$ (resp. $\{h_i^r\}$). Using the definition of $h_1^r$ and Markov's inequality,

$$
\begin{aligned}
\mathbf{E}^r[(h_1^r)^2] &= \mathbf{E}^r[g(v_1^r)^2; g(v_1^r) > r^{1/8}] \\
&\le r^{-p'/8}\langle g^{2+p'}, \nu^r\rangle,
\end{aligned}
\tag{3.44}
$$

which tends to zero as $r \to \infty$ by assumption. Thus, letting $r \to \infty$ in (3.43), we see that the quantity in the outer parentheses on the last line tends to zero uniformly in $m$ and is thus smaller than $\eta'[4(\lceil L/\varepsilon_g\rceil + 1)]^{-1}$ for all sufficiently large $r$. □

3.5. *Moment estimates.* We will need the following simple lemma.

LEMMA 3.5. *Suppose $\xi_k, \xi \in \mathcal{M}_F$ such that $\xi_k \xrightarrow{w} \xi$, as $k \to \infty$. Suppose further that for some $q' > 0$ and $1 < M < \infty$,*

$$
\limsup_{k\to\infty}(\langle 1, \xi_k\rangle \vee \langle \chi^{1+q'}, \xi_k\rangle) \le M. \tag{3.45}
$$

*Then for any $q \in [0, q')$, there exists an $M_q \in [M, \infty)$ such that $\langle \chi^{1+q}, \xi\rangle \le M_q$ and*

$$
\langle \chi^{1+q}, \xi_k\rangle \longrightarrow \langle \chi^{1+q}, \xi\rangle, \qquad k \to \infty. \tag{3.46}
$$

*Furthermore, $M_q$ depends only on $M$, $q'$ and $q$.*

PROOF. Since $\xi_k \xrightarrow{w} \xi$, as $k \to \infty$, (3.45) implies that for any $x \in \mathbb{R}_+$,

$$
\langle \mathbb{1}_{(x,\infty)}, \xi\rangle \le \langle 1, \xi\rangle = \lim_{k\to\infty}\langle 1, \xi_k\rangle \le M. \tag{3.47}
$$

For any $x \ge 1$, (3.45) and Markov's inequality yield

$$
\begin{aligned}
\limsup_{k\to\infty}\langle \mathbb{1}_{(x,\infty)}, \xi_k\rangle &\le \limsup_{k\to\infty} x^{-(1+q')}\langle \chi^{1+q'}, \xi_k\rangle \\
&\le x^{-(1+q')}M.
\end{aligned}
\tag{3.48}
$$

So by the Portmanteau theorem ([1], Theorem 2.1), for any $x \geq 1$,

$$\langle \mathbb{1}_{(x,\infty)}, \xi \rangle \leq \limsup_{k\to\infty} \langle \mathbb{1}_{(x,\infty)}, \xi_k \rangle \leq x^{-(1+q')} M. \tag{3.49}$$

Since $\chi^{1+q}(0) = 0$, we can use integration by parts together with (3.47) and (3.49) to obtain

$$\begin{aligned} \langle \chi^{1+q}, \xi \rangle &= \int_0^\infty (1+q) x^q \langle \mathbb{1}_{(x,\infty)}, \xi \rangle \, dx \\ &\leq (1+q) M \left( 1 + \int_1^\infty x^{-(1+q'-q)} \, dx \right) < \infty. \end{aligned} \tag{3.50}$$

We can take

$$\begin{aligned} M_q &= (1+q) M \left( 1 + \int_1^\infty x^{-(1+q'-q)} \, dx \right) \\ &< \infty, \end{aligned} \tag{3.51}$$

which depends only on $M$, $q'$ and $q$. To show (3.46), first choose $q_1 \in (q, q')$ and $M_{q_1} \in [M, \infty)$ such that

$$\langle \chi^{1+q_1}, \xi \rangle \leq M_{q_1}. \tag{3.52}$$

For any $y \geq 1$ that is not an atom of $\xi$,

$$\begin{aligned} &|\langle \chi^{1+q}, \xi_k \rangle - \langle \chi^{1+q}, \xi \rangle| \\ &\quad \leq |\langle \chi^{1+q} \mathbb{1}_{[0,y]}, \xi_k \rangle - \langle \chi^{1+q} \mathbb{1}_{[0,y]}, \xi \rangle| \\ &\qquad + |\langle \chi^{1+q} \mathbb{1}_{(y,\infty)}, \xi_k \rangle - \langle \chi^{1+q} \mathbb{1}_{(y,\infty)}, \xi \rangle| \\ &\quad \leq |\langle \chi^{1+q} \mathbb{1}_{[0,y]}, \xi_k \rangle - \langle \chi^{1+q} \mathbb{1}_{[0,y]}, \xi \rangle| \\ &\qquad + |y^{-(q_1-q)} \langle \chi^{1+q_1} \mathbb{1}_{[1,\infty)}, \xi_k \rangle| \\ &\qquad + |y^{-(q_1-q)} \langle \chi^{1+q_1} \mathbb{1}_{[1,\infty)}, \xi \rangle|. \end{aligned} \tag{3.53}$$

Since $y$ is not an atom of $\xi$, we see that the term on the fourth line of (3.53) tends to zero as $k \to \infty$. So since $q_1 < q'$ and $M \leq M_{q_1}$, (3.45), (3.52) and (3.53) imply

$$\limsup_{k\to\infty} |\langle \chi^{1+q}, \xi_k \rangle - \langle \chi^{1+q}, \xi \rangle| \leq 2 y^{-(q_1-q)} M_{q_1}. \tag{3.54}$$

Since $\xi$ can have at most countably many atoms, $y$ can be chosen arbitrarily large.

□

The next lemma provides an upper bound for three of the moments of the fluid scaled state descriptors $\{\bar{\mu}^{r,m}(\cdot)\}$. This estimate will be fundamental to many subsequent proofs.

LEMMA 3.6. *Assume* (Q.1) *and let* $T>1$ *and* $0<\eta<1$ *be given. Then there exists* $p>0$ *and* $M_T>1$ *such that for any* $L>1$,

$$\liminf_{r\to\infty} \mathbf{P}^r(\Omega_1^r) \geq 1-\frac{\eta}{6}, \tag{3.55}$$

*where*

$$\Omega_1^r = \left\{ \sup_{m\leq\lfloor rT\rfloor} \|\langle 1,\bar{\mu}^{r,m}(\cdot)\rangle \vee \langle \chi,\bar{\mu}^{r,m}(\cdot)\rangle \vee \langle \chi^{1+p},\bar{\mu}^{r,m}(\cdot)\rangle\|_L \leq M_T \right\}.$$

PROOF. Choose $p>0$ such that $6p+2p^2<\theta$. By (2.31)–(2.35) and Corollary 3.2, using $\eta'=\eta/36$ there, there exists $M_0>1$ such that for any $L>1$,

$$\liminf_{r\to\infty} \mathbf{P}^r(A_1^r) \geq 1-\frac{\eta}{12}, \tag{3.56}$$

where

$$A_1^r = \left\{ \langle 1,\bar{\mu}^r(0)\rangle \vee \sup_{m\leq\lfloor rT\rfloor} \|\langle \chi,\bar{\mu}^{r,m}(\cdot)\rangle\|_L \vee \langle \chi^{1+p},\bar{\mu}^r(0)\rangle \leq M_0 \right\}.$$

Define

$$M_1 = 32(\alpha\vee 1)M_0, \tag{3.57}$$

$$M_T = M_1 \vee (M_0 + M_1^{1+p}(1+p^{-1})2\alpha\langle\chi^{2+2p},\nu\rangle + 4\alpha\langle\chi^{1+p},\nu\rangle). \tag{3.58}$$

Let $L>1$ be fixed and let $A_2^r$ be the event on which

$$\sup_{m\leq\lfloor rT\rfloor} \|\bar{E}^{r,m}(\cdot+4M_0)-\bar{E}^{r,m}(\cdot)\|_L \leq 8\alpha M_0, \tag{3.59}$$

$$\sup_{m\leq\lfloor rT+L\rfloor} \frac{1}{r}\sum_{i=r\bar{E}^{r,m}(0)+1}^{r\bar{E}^{r,m}(1)} \chi^{1+p}(v_i^r) \leq 2\alpha\langle\chi^{1+p},\nu\rangle, \tag{3.60}$$

$$\sup_{m\leq\lfloor rT+L\rfloor} \frac{1}{r}\sum_{i=r\bar{E}^{r,m}(0)+1}^{r\bar{E}^{r,m}(1)} \chi^{2+2p}(v_i^r) \leq 2\alpha\langle\chi^{2+2p},\nu\rangle. \tag{3.61}$$

Define $A^r = A_1^r \cap A_2^r$. We will apply Lemma 3.4 three times in order to obtain a lower bound on the probability of $A_2^r$, and hence $A^r$. By choice of $p$, (2.27) and Lemma 3.5 [with $q'\in(p,\theta)$ or $q'\in(1+2p,3+\theta)$ there] imply that the assumptions of Lemma 3.4 are satisfied for each of the three choices $g\equiv 1$, $g\equiv\chi^{1+p}$ and $g\equiv\chi^{2+2p}$. Note that Lemma 3.4 still holds if, in (3.16), we replace the supremum over $m\leq\lfloor rT\rfloor$ by the supremum over $m\leq\lfloor rT+L\rfloor$ (see proof of Lemma 3.4). Thus, by applying Lemma 3.4 three times [resp., for $(l,\varepsilon',\eta',g)$ equal to $(4M_0, 4\alpha M_0\wedge(1/2),\eta/36,1)$, $(1,\alpha\langle\chi^{1+p},\nu\rangle\wedge$

$(1/2), \eta/36, \chi^{1+p})$ and $(1, \alpha\langle\chi^{2+2p}, \nu\rangle \wedge (1/2), \eta/36, \chi^{2+2p})]$ and combining with (3.56) one obtains

$$\liminf_{r\to\infty} \mathbf{P}^r(A^r) \geq 1 - \frac{\eta}{6}. \tag{3.62}$$

Thus, by (3.56) and (3.62) it suffices to show that on $A^r$, for all $m \leq \lfloor rT \rfloor$,

$$\|\langle 1, \bar{\mu}^{r,m}(\cdot)\rangle \vee \langle \chi^{1+p}, \bar{\mu}^{r,m}(\cdot)\rangle\|_L \leq M_T. \tag{3.63}$$

Note that by (3.57) and (3.59), $\|\bar{E}^r(\cdot + 4M_0) - \bar{E}^r(\cdot)\|_{\lfloor rT\rfloor + L} \leq M_1/4$ on $A^r$. So for $c = M_1$, $l = 4M_0$, $t_0 = 0$ and $t_1 = \lfloor rT \rfloor + L$, conditions (i)–(iii) of Lemma 3.3 are satisfied on $A^r$, which implies that on $A^r$,

$$\begin{aligned}\sup_{m\leq\lfloor rT\rfloor} \|\langle 1, \bar{\mu}^{r,m}(\cdot)\rangle\|_L &= \|\langle 1, \bar{\mu}^r(\cdot)\rangle\|_{\lfloor rT\rfloor + L} \\ &\leq M_1 \leq M_T.\end{aligned} \tag{3.64}$$

It remains to consider the $(1+p)$th moment. First, consider any job $i \geq 1$ such that $U_i^r r^{-1} \leq t_1 < t_2$ for some times $t_1, t_2 \in [0, \lfloor rT \rfloor + L]$. If $\bar{\mu}^r(t^*) = \mathbf{0}$ for some $t^* \in [t_1, t_2]$, then

$$\begin{aligned}&\chi^{1+p}\Big(v_i^r - \bar{S}^r_{U_i^r r^{-1}, t_2}\Big) \\ &\quad\leq \chi^{1+p}(v_i^r - \bar{S}^r_{U_i^r r^{-1}, t^*}) \\ &\quad= 0 \leq \chi^{1+p}\left(v_i^r - \frac{t_2 - t_1}{M_1}\right),\end{aligned} \tag{3.65}$$

where the equality above follows from (3.2), by setting $g \equiv \chi^{1+p}$, $m = 0$, $t = U_i^r r^{-1}$ and $h = t^* - U_i^r r^{-1}$ there. If $\bar{\mu}^r(s) \neq \mathbf{0}$ for all $s \in [t_1, t_2]$, then on $A^r$,

$$\begin{aligned}\bar{S}^r_{t_1,t_2} &\geq (t_2 - t_1) \inf_{s\in[t_1,t_2]} \varphi(\langle 1, \bar{\mu}^r(s)\rangle) \\ &\geq (t_2 - t_1)\Big(\sup_{s\in[t_1,t_2]} \langle 1, \bar{\mu}^r(s)\rangle\Big)^{-1} \\ &\geq (t_2 - t_1) M_1^{-1},\end{aligned} \tag{3.66}$$

where the last inequality is by (3.64). So in this case, on $A^r$, we also have

$$\begin{aligned}&\chi^{1+p}\Big(v_i^r - \bar{S}^r_{U_i^r r^{-1}, t_2}\Big) \\ &\quad\leq \chi^{1+p}(v_i^r - \bar{S}^r_{t_1,t_2}) \\ &\quad\leq \chi^{1+p}\left(v_i^r - \frac{t_2 - t_1}{M_1}\right).\end{aligned} \tag{3.67}$$

Combining (3.65) and (3.67) yields that on $A^r$,

$$\begin{aligned}
&\chi^{1+p}\Big(v_i^r - \bar{S}^r_{U_i^r r^{-1}, t_2}\Big) \\
&\quad\le \chi^{1+p}\left(v_i^r - \frac{t_2 - t_1}{M_1}\right) \\
&\quad\le \mathbb{1}_{[(t_2-t_1)M_1^{-1},\infty)}(v_i^r)\chi^{1+p}(v_i^r) \\
&\quad\le (\mathbb{1}_{[(t_2-t_1)M_1^{-1},\infty)}\chi^{-(1+p)})(v_i^r)\chi^{2+2p}(v_i^r) \\
&\quad\le \left(\frac{M_1}{t_2-t_1}\right)^{(1+p)}\chi^{2+2p}(v_i^r).
\end{aligned} \tag{3.68}$$

To complete the proof, we must show that $\langle \chi^{1+p}, \bar{\mu}^{r,m'}(s)\rangle \le M_T$ on $A^r$, for any $m' \le \lfloor rT \rfloor$ and $s \in [0, L]$. Note that since $L > 1$, $s > 1$ implies that $\bar{\mu}^{r,0}(s) = \bar{\mu}^{r,1}(s_0)$ for some $s_0 \in [0, L-1]$, and $\bar{\mu}^{r,1}(s) = \bar{\mu}^{r,2}(s_1)$ for some $s_1 \in [0, L - 1]$. Thus, it suffices to consider the case $m' \in \{0, 1\}$ and $s \le 1$ and the case $m' \in \{2, \dots, \lfloor rT \rfloor\}$ and $s \in [0, L]$. If $m' \in \{0, 1\}$ and $s \le 1$, we have on $A^r$,

$$\begin{aligned}
&\langle \chi^{1+p}, \bar{\mu}^{r,m'}(s)\rangle \\
&\quad= \langle \chi^{1+p}, \bar{\mu}^r(m'+s)\rangle \\
&\quad= \langle \chi^{1+p}(\cdot - \bar{S}^r(m'+s)), \bar{\mu}^r(0)\rangle \\
&\qquad + \frac{1}{r}\sum_{i=1}^{r\bar{E}^r(m'+s)} \chi^{1+p}(v_i^r - \bar{S}^r_{U_i^r r^{-1}, m'+s}) \\
&\quad\le \langle \chi^{1+p}, \bar{\mu}^r(0)\rangle + \frac{1}{r}\sum_{i=1}^{r\bar{E}^r(m'+s)} \chi^{1+p}(v_i^r) \\
&\quad\le M_0 + \sum_{m=0}^{m'}\left(\frac{1}{r}\sum_{i=r\bar{E}^{r,m}(0)+1}^{r\bar{E}^{r,m}(1)} \chi^{1+p}(v_i^r)\right) \\
&\quad\le M_0 + 4\alpha\langle \chi^{1+p}, \nu\rangle \\
&\quad\le M_T,
\end{aligned} \tag{3.69}$$

where the second equality above uses (3.2), with $g \equiv \chi^{1+p}$, $m = 0$, $t = 0$ and $h = m' + s$ there, the first inequality follows by simply ignoring any processing which has taken place by time $m' + s$, the second inequality follows from the definition of $A^r$ for the first term and by including additional summands in the second term, and the third inequality is by (3.60) and the fact that $m' \le 1$.

Last, we consider the case $m' \in \{2, \dots, \lfloor rT \rfloor\}$ and $s \in [0, L]$. Let $N = \lfloor m' + s \rfloor - 2$ and note that $N \in \{0, \dots, \lfloor rT + L \rfloor - 2\}$. Note also that by replacing $m$ by $\lfloor m' + s \rfloor - m - 1$, the following estimate holds:

$$\begin{aligned}
&\sum_{m=0}^{N} \left( \frac{M_1}{\lfloor m' + s \rfloor - m - 1} \right)^{1+p} \\
&\qquad = \sum_{m=1}^{N+1} \left( \frac{M_1}{m} \right)^{1+p} \\
&\qquad \le M_1^{1+p} \left( 1 + \int_1^\infty x^{-(1+p)} \, dx \right) \\
&\qquad = M_1^{1+p} (1 + p^{-1}).
\end{aligned} \tag{3.70}$$

For $m' \in \{2, \dots, \lfloor rT \rfloor\}$ and $s \in [0, L]$ we have on $A^r$,

$$\begin{aligned}
&\langle \chi^{1+p}, \bar{\mu}^{r,m'}(s) \rangle \\
&\quad \le M_0 + \frac{1}{r} \sum_{i=1}^{r\bar{E}^r(m'+s)} \chi^{1+p}(v_i^r - \bar{S}^r_{U_i^r r^{-1}, m'+s}) \\
&\quad \le M_0 + \sum_{m=0}^{N} \left( \frac{1}{r} \sum_{i=r\bar{E}^{r,m}(0)+1}^{r\bar{E}^{r,m}(1)} \chi^{1+p}(v_i^r - \bar{S}^r_{U_i^r r^{-1}, m'+s}) \right) \\
&\qquad + \frac{1}{r} \sum_{i=r\bar{E}^r(N+1)+1}^{r\bar{E}^r(m'+s)} \chi^{1+p}(v_i^r) \\
&\quad \le M_0 + \sum_{m=0}^{N} \left( \left( \frac{M_1}{\lfloor m' + s \rfloor - m - 1} \right)^{(1+p)} \frac{1}{r} \sum_{i=r\bar{E}^{r,m}(0)+1}^{r\bar{E}^{r,m}(1)} \chi^{2+2p}(v_i^r) \right) \\
&\qquad + \sum_{m=N+1}^{N+2} \left( \frac{1}{r} \sum_{i=r\bar{E}^{r,m}(0)+1}^{r\bar{E}^{r,m}(1)} \chi^{1+p}(v_i^r) \right) \\
&\quad \le M_0 + M_1^{1+p}(1 + p^{-1}) 2\alpha \langle \chi^{2+2p}, \nu \rangle + 4\alpha \langle \chi^{1+p}, \nu \rangle \\
&\quad \le M_T,
\end{aligned} \tag{3.71}$$

where the first inequality is obtained as in (3.69), the second follows by ignoring any processing of jobs in the third term; the third inequality is by (3.68) (with $t_1 = m + 1$ and $t_2 = m' + s$) for the second term and by including additional summands in the third term; the fourth inequality is by (3.60), (3.61) and (3.70) and the last is by (3.58). □

3.6. *Combining estimates.* To aid the reader, we now provide a lemma that combines several of the preliminary results obtained in Sections 3.2–3.5, such that they hold simultaneously, with respect to specific choices of various constants. It is analogous to Lemma 5.2 in [6], in that subsequent references to results obtained thus far will predominantly reference this lemma. We will need one more definition.

DEFINITION 3.7. Let $\mathcal{C} = \{g \in \mathbf{C}_b^1(\mathbb{R}_+) : g(0) = g'(0) = 0\}$, and let $\tilde{\mathcal{C}} = \{g \in \mathcal{C} : g \text{ has compact support}\}$. Let $\mathcal{V} = \{g_k^{\mathcal{V}}\}_{k\in\mathbb{N}} \subset \tilde{\mathcal{C}}$ be a countable subset such that for any $g \in \tilde{\mathcal{C}}$, there is a sequence $\{g_k\}_{k=1}^{\infty} \subset \mathcal{V}$ that together with $\{g_k'\}_{k=1}^{\infty}$ is uniformly bounded and satisfies

$$g_k \longrightarrow g \quad \text{and} \quad g_k' \longrightarrow g', \qquad k \to \infty, \tag{3.72}$$

pointwise on $\mathbb{R}_+$.

An example of how such a set $\mathcal{V} \subset \tilde{\mathcal{C}}$ can be chosen is given in the proof of property (iii) for Theorem 5.1 in [6]. Recall that $\mathbf{M}_\nu$ is defined in Definition 2.1 and the metric $\mathbf{d}[\cdot,\cdot]$ is defined in (1.4).

LEMMA 3.8. *Assume* (Q.1) *and let* $T, L > 1$ *and* $0 < \eta < 1$ *be given. Let* $\{\varepsilon_n\}_{n=1}^{\infty} \subset (0,1)$ *be a sequence such that* $\varepsilon_n \downarrow 0$, *as* $n \to \infty$. *Then for each* $n \in \mathbb{N}$, *there exist strictly positive constants* $l_n$, $p$, $M_T$, $\kappa_n$, $\gamma_n$, $N_n$, $r_n$ *and events* $\{B_n^r\}_{r\in\mathcal{R}}$, *such that:*

(i) $\mathbf{P}^r(B_n^r) \geq 1 - \eta$, *for all* $r > r_n$,
(ii) $r_{n+1} > r_n$ *and* $B_{n+1}^r \subset B_n^r$, *and*
(iii) *for each* $r \in \mathcal{R}$ *and each* $m \leq \lfloor rT \rfloor$, *on* $B_n^r$ *the following hold:*

$$l_n < \frac{\varepsilon_n}{32\alpha}, \tag{3.73}$$

$$\sup_{t\in[0,L]} \bar{E}^{r,m}(t + l_n) - \bar{E}^{r,m}(t) \leq \frac{\varepsilon_n}{16}, \tag{3.74}$$

$$\sup_{\substack{k\leq n \\ j=1,\ldots,\lfloor L/l_n\rfloor}} \left| \frac{1}{r} \sum_{i=r\bar{E}^{r,m}(0)+1}^{r\bar{E}^{r,m}(jl_n)} g_k^{\mathcal{V}}(v_i^r) - \alpha j l_n \langle g_k^{\mathcal{V}}, \nu\rangle \right| \leq \varepsilon_n, \tag{3.75}$$

*for* $g_k^{\mathcal{V}} \in \mathcal{V}$,

$$\sup_{t\in[0,L]} \langle 1, \bar{\mu}^{r,m}(t)\rangle \vee \langle \chi, \bar{\mu}^{r,m}(t)\rangle \vee \langle \chi^{1+p}, \bar{\mu}^{r,m}(t)\rangle \leq M_T, \tag{3.76}$$

$$\inf_{\xi\in\mathbf{M}_\nu} \mathbf{d}[\bar{\mu}^r(0), \xi] \vee |\langle \chi, \bar{\mu}^r(0)\rangle - \langle \chi, \xi\rangle| < \varepsilon_n, \tag{3.77}$$

$$\kappa_n < \frac{l_n}{2M_T}, \tag{3.78}$$

$$\sup_{x\in\mathbb{R}_+} \langle \mathbb{1}_{[x,x+\kappa_n]}, \bar\mu^r(0)\rangle \le \frac{\varepsilon_n}{4}, \tag{3.79}$$

$$\gamma_n < \frac{\kappa_n\varepsilon_n}{4} \wedge \frac{l_n}{4}, \tag{3.80}$$

$$\sup_{t\in[0,L]} |\langle \chi, \bar\mu^{r,m}(t)\rangle - \langle \chi, \bar\mu^{r,m}(0)\rangle| \le \frac{\gamma_n}{4}, \tag{3.81}$$

$$N_n = \lceil 64M_T^2/(\kappa_n\varepsilon_n^2)\rceil + 1, \tag{3.82}$$

$$\begin{aligned}&\sup_{t\in[0,L]} \frac{1}{r}\sum_{i=r\bar E^{r,m}(t)+1}^{r\bar E^{r,m}(t+l_n)} \mathbb{1}_{[k\kappa_n,(k+1)\kappa_n)}(v_i^r)\\ &\qquad\le \frac{\varepsilon_n}{8}\langle \mathbb{1}_{[(k-1/2)\kappa_n,(k+3/2)\kappa_n)}, \nu\rangle\end{aligned} \tag{3.83}$$

*for all* $k\in\{0,1,\dots,N_n\}$.

Proof. Let $\tilde\eta=\eta/6$. For each $n\in\mathbb{N}$, choose $0<l_n<\varepsilon_n/(32\alpha)$ such that $l_{n+1}<l_n$ and let $\tilde B^{r,m}_{n,1}$ be the event on which both (3.74) and (3.75) hold. Let $\tilde B^r_{n,1}=\bigcap_{m\le\lfloor rT\rfloor}\tilde B^{r,m}_{n,1}$ and let $B^r_{n,1}=\bigcap_{\dot n\le n}\tilde B^r_{\dot n,1}$. The first goal is to show that for each $n\in\mathbb{N}$,

$$\liminf_{r\to\infty}\mathbf{P}^r(B^r_{n,1})\ge 1-\tilde\eta.$$

Fix $n\in\mathbb{N}$ and consider any $\dot n\le n$. Since (3.74) (with $n=\dot n$) holds on $\tilde B^r_{\dot n,1}$ for each $m\le\lfloor rT\rfloor$, and since (3.75) (with $n=\dot n$) holds on $\tilde B^r_{\dot n,1}$ for each $m\le\lfloor rT\rfloor$, we must apply Lemma 3.4 $1+\dot n\lfloor L/l_{\dot n}\rfloor$ times in order to estimate the asymptotic probability of $\tilde B^r_{\dot n,1}$. More precisely, we must apply it once to guarantee (3.74) and also once for each case $k\le\dot n$ and $j\le\lfloor L/l_{\dot n}\rfloor$ in (3.75). Observe that by (2.25), the assumptions of Lemma 3.4 are satisfied with $l=l_{\dot n}$, $\varepsilon'=\varepsilon_{\dot n}/16$, $\eta'=\tilde\eta(n(1+n\lfloor L/l_n\rfloor))^{-1}$, $p'>0$ arbitrary and $g\equiv1$. Also observe that for each $k=1,\dots,\dot n$ and $j=1,\dots,\lfloor L/l_{\dot n}\rfloor$, the assumptions of Lemma 3.4 are satisfied with $l=jl_{\dot n}$, $\varepsilon'=\varepsilon_{\dot n}$, $\eta'=\tilde\eta(n(1+n\lfloor L/l_n\rfloor))^{-1}$, $p'>0$ arbitrary and $g\equiv g_k^{\mathcal{V}}$. Thus, $1+\dot n\lfloor L/l_{\dot n}\rfloor$ applications of Lemma 3.4 yield

$$\liminf_{r\to\infty}\mathbf{P}^r(\tilde B^r_{\dot n,1})\ge 1-(1+\dot n\lfloor L/l_{\dot n}\rfloor)\frac{\tilde\eta}{n(1+n\lfloor L/l_n\rfloor)}\ge 1-\frac{\tilde\eta}{n}. \tag{3.84}$$

This implies that for each $n\in\mathbb{N}$,

$$\liminf_{r\to\infty}\mathbf{P}^r(B^r_{n,1})\ge 1-\sum_{\dot n=1}^n\frac{\tilde\eta}{n}=1-\tilde\eta. \tag{3.85}$$

Note that by definition, $B^r_{n+1,1}\subset B^r_{n,1}$ for $n\in\mathbb{N}$.

Next, let $p$ and $M_T$ be the constants given by Lemma 3.6 and let $B_2^{r,m}$ be the event appearing in (3.76). Let $B_2^r = \bigcap_{m \le \lfloor rT \rfloor} B_2^{r,m}$. Then by Lemma 3.6,

$$\liminf_{r\to\infty} \mathbf{P}^r(B_2^r) \ge 1 - \tilde{\eta}. \tag{3.86}$$

For each $n \in \mathbb{N}$, let $B_{n,3}^r$ be the event appearing in (3.77). Recall that "$\sim$" denotes equivalence in distribution. By (2.35) and the Skorohod representation theorem, there exist random pairs $(\Theta^r, X^r), (\tilde{\Theta}, \tilde{X})$, taking values in $\mathcal{M}_{\mathrm{F}} \times \mathbb{R}_+$ and defined on some common probability space $(\Omega, \mathscr{F}, \mathbf{P})$, such that $(\tilde{\Theta}, \tilde{X}) \sim (\Theta, \langle \chi, \Theta \rangle)$ and $(\Theta^r, X^r) \sim (\bar{\mu}^r(0), \langle \chi, \bar{\mu}^r(0) \rangle)$ for each $r \in \mathcal{R}$ and such that as $r \to \infty$,

$$(\Theta^r, X^r) \longrightarrow (\tilde{\Theta}, \tilde{X}) \qquad \text{a.s.} \tag{3.87}$$

So by (2.30) and (3.87),

$$\begin{aligned} &\liminf_{r\to\infty} \mathbf{P}^r(B_{n,3}^r) \\ &\quad = \liminf_{r\to\infty} \mathbf{P}\left( \inf_{\xi \in \mathbf{M}_\nu} \mathbf{d}[\Theta^r, \xi] \vee |X^r - \langle \chi, \xi \rangle| < \varepsilon_n \right) \\ &\quad \ge 1 - \tilde{\eta}. \end{aligned} \tag{3.88}$$

Since $\{\varepsilon_n\}$ is a decreasing sequence, $B_{n+1,3}^r \subset B_{n,3}^r$ for $n \in \mathbb{N}$.

Next, for each $n \in \mathbb{N}$, we use (2.34) (with $\varepsilon = \varepsilon_n/4$) and (2.31) to choose $\kappa_n > 0$ and $M_n > 0$ so that (3.78) holds, $\kappa_{n+1} < \kappa_n$ for $n \in \mathbb{N}$ and

$$\mathbf{P}\left( \sup_{x \in \mathbb{R}_+} \langle \mathbb{1}_{[x, x+2\kappa_n]}, \Theta \rangle \vee \langle \mathbb{1}_{[M_n, \infty)}, \Theta \rangle < \frac{\varepsilon_n}{4} \right) \ge 1 - \frac{\tilde{\eta}}{4n^2}. \tag{3.89}$$

For $k \in \{0, \ldots, \lceil M_n/\kappa_n \rceil - 1\}$, define $I_k = [k\kappa_n, (k+2)\kappa_n]$ and let $I_{\lceil M_n/\kappa_n \rceil} = [M_n, \infty)$. Note that for any $x \in \mathbb{R}_+$, $[x, x+\kappa_n] \subset I_k$ for some $k \in \{0, \ldots, \lceil M_n/\kappa_n \rceil\}$. Consider the random measures $\Theta^r, \tilde{\Theta}$ defined above such that (3.87) holds. Since for each $n \in \mathbb{N}$, $\{I_k\}$ is a finite set of closed intervals, a trivial generalization of the Portmanteau theorem ([1], Theorem 2.1) to finite measures yields for each $n \in \mathbb{N}$,

$$\limsup_{r\to\infty} \max_{k \le \lceil M_n/\kappa_n \rceil} \langle \mathbb{1}_{I_k}, \Theta^r \rangle \le \max_{k \le \lceil M_n/\kappa_n \rceil} \langle \mathbb{1}_{I_k}, \tilde{\Theta} \rangle \qquad \text{a.s.} \tag{3.90}$$

Since $\tilde{\Theta} \sim \Theta$, (3.89) and (3.90) imply that for each $n \in \mathbb{N}$,

$$\begin{aligned} &\liminf_{r\to\infty} \mathbf{P}\left( \sup_{x \in \mathbb{R}_+} \langle \mathbb{1}_{[x, x+\kappa_n]}, \Theta^r \rangle \le \frac{\varepsilon_n}{4} \right) \\ &\quad \ge \liminf_{r\to\infty} \mathbf{P}\left( \max_{k \le \lceil M_n/\kappa_n \rceil} \langle \mathbb{1}_{I_k}, \Theta^r \rangle \le \frac{\varepsilon_n}{4} \right) \\ &\quad \ge 1 - \frac{\tilde{\eta}}{2n^2}. \end{aligned} \tag{3.91}$$

Let $\tilde{B}^r_{n,4}$ be the event appearing in (3.79). We have by (3.91),

$$\begin{aligned}\liminf_{r\to\infty} \mathbf{P}^r(\tilde{B}^r_{n,4}) \\ &= \liminf_{r\to\infty} \mathbf{P}\left(\sup_{x\in\mathbb{R}_+} \langle \mathbb{1}_{[x,x+\kappa_n]}, \Theta^r\rangle \le \frac{\varepsilon_n}{4}\right) \\ &\ge 1 - \frac{\tilde{\eta}}{2n^2}.\end{aligned} \tag{3.92}$$

Letting $B^r_{n,4} = \bigcap_{i\le n} \tilde{B}^r_{i,4}$, we see that for each $n\in\mathbb{N}$,

$$\liminf_{r\to\infty} \mathbf{P}^r(B^r_{n,4}) \ge 1 - \sum_{i=1}^{n} \frac{\tilde{\eta}}{2i^2} \ge 1 - \tilde{\eta}. \tag{3.93}$$

By definition, $B^r_{n+1,4} \subset B^r_{n,4}$ for $n\in\mathbb{N}$.

For each $n\in\mathbb{N}$, choose $\gamma_n > 0$ such that (3.80) holds and such that $\gamma_{n+1} < \gamma_n$ for $n\in\mathbb{N}$. Let $B^{r,m}_{n,5}$ be the event appearing in (3.81) and let $B^r_{n,5} = \bigcap_{m\le\lfloor rT\rfloor} B^{r,m}_{n,5}$. Then by Corollary 3.2, using $\eta' = \tilde{\eta}$ and $\gamma = \gamma_n$, we have for each $n\in\mathbb{N}$,

$$\liminf_{r\to\infty} \mathbf{P}^r(B^r_{n,5}) \ge 1 - \tilde{\eta}. \tag{3.94}$$

Since $\{\gamma_n\}$ is a decreasing sequence, $B^r_{n+1,5} \subset B^r_{n,5}$ for $n\in\mathbb{N}$.

Finally, define $N_n$ as in (3.82) for each $n\in\mathbb{N}$, and for $k\in\{0,\dots,N_n\}$, choose $g^n_k \in \mathbf{C}_b(\mathbb{R}_+)$ such that $g^n_k$ is nonnegative and such that for all $x\in\mathbb{R}_+$, $\mathbb{1}_{[k\kappa_n,(k+1)\kappa_n)}(x) \le g^n_k(x) \le \mathbb{1}_{[(k-1/2)\kappa_n,(k+3/2)\kappa_n)}(x)$. Note that for $k=0$ and $x\in\mathbb{R}_+$, $\mathbb{1}_{[(k-1/2)\kappa_n,(k+3/2)\kappa_n)}(x) = \mathbb{1}_{[0,3/2\kappa_n)}(x)$. Let $\tilde{B}^{r,m,k}_{n,6}$ be the event in (3.83) and let $\tilde{B}^r_{n,6} = \bigcap_{k\le N_n}\bigcap_{m\le\lfloor rT\rfloor} \tilde{B}^{r,m,k}_{n,6}$. Define

$$G^r_{n,6} = \bigcap_{k=0}^{N_n} \left\{ \sup_{m\le\lfloor rT\rfloor} \sup_{t\in[0,L]} \frac{1}{r} \sum_{i=r\bar{E}^{r,m}(t)+1}^{r\bar{E}^{r,m}(t+l_n)} g^n_k(v^r_i) \le \frac{\varepsilon_n}{8}\langle g^n_k, \nu\rangle \right\}. \tag{3.95}$$

Fix $n\in\mathbb{N}$ and consider $\dot{n}\le n$. Observe that by (2.25), we have that for each $k\in\{0,\dots,N_{\dot{n}}\}$, the assumptions of Lemma 3.4 are satisfied with $l = l_{\dot{n}}$, $\varepsilon' = (\varepsilon_{\dot{n}}/16)\langle g^{\dot{n}}_k, \nu\rangle$, $\eta' = \tilde{\eta}(n(N_n+1))^{-1}$, $p' > 0$ arbitrary, and $g \equiv g^{\dot{n}}_k$. Thus $N_{\dot{n}}+1$ applications of Lemma 3.4 yield

$$\liminf_{r\to\infty} \mathbf{P}^r(G^r_{\dot{n},6}) \ge 1 - \frac{\tilde{\eta}(N_{\dot{n}}+1)}{n(N_n+1)} \ge 1 - \frac{\tilde{\eta}}{n}, \tag{3.96}$$

where we have used the fact that $l_{\dot{n}} < \varepsilon_{\dot{n}}/(32\alpha)$, and the second inequality follows from (3.82) since $\{\kappa_n\}$ and $\{\varepsilon_n\}$ are decreasing sequences. Since $G^r_{\dot{n},6} \subset \tilde{B}^r_{\dot{n},6}$ for each $\dot{n}\le n$, we have

$$\liminf_{r\to\infty} \mathbf{P}^r(\tilde{B}^r_{\dot{n},6}) \ge 1 - \frac{\tilde{\eta}}{n}. \tag{3.97}$$

Letting $B^r_{n,6} = \bigcap_{\dot n \le n} \tilde B^r_{\dot n,6}$, we see that $B^r_{n+1,6} \subset B^r_{n,6}$ for $n \in \mathbb{N}$, and

$$\liminf_{r\to\infty} \mathbf{P}^r(B^r_{n,6}) \ge 1 - \sum_{\dot n=1}^{n} \frac{\tilde\eta}{n} = 1 - \tilde\eta. \tag{3.98}$$

For each $r \in \mathcal{R}$, define $B^r_n = B^r_{n,1} \cap B^r_2 \cap B^r_{n,3} \cap B^r_{n,4} \cap B^r_{n,5} \cap B^r_{n,6}$ and note that by construction, $B^r_{n+1} \subset B^r_n$ for $n \in \mathbb{N}$. Then by choice of $\tilde\eta$, we see that for each $n \in \mathbb{N}$, there exists $r_n \in \mathcal{R}$ such that $r_{n+1} > r_n$ and $r > r_n$ implies

$$\mathbf{P}^r(B^r_n) \ge 1 - \eta. \tag{3.99}$$

□

3.7. *Fluid scale analysis.* In this section we build on the results of Sections 3.1–3.6 to prove several important properties of the fluid scaled state descriptors $\{\bar\mu^{r,m}(\cdot)\}$. These properties will be used at the end of the section to prove the desired relative compactness property (see Corollary 3.16). Lemmas 3.9 and 3.10 give an upper bound for the amount of mass that $\bar\mu^{r,m}(t)$ can have concentrated near zero. They are analogous to Lemma 5.5 in [6]. In Lemmas 3.11 and 3.12 we give an upper bound for the shifted fluid scaled queue length $\langle 1, \bar\mu^{r,m}(\cdot)\rangle$ when the shifted fluid scaled initial workload $\langle \chi, \bar\mu^{r,m}(0)\rangle$ is small, and a lower bound for $\langle 1, \bar\mu^{r,m}(\cdot)\rangle$ when $\langle \chi, \bar\mu^{r,m}(0)\rangle$ is bounded away from zero. These two lemmas are analogous to Lemmas 5.3 and 5.4 in [6]. We then combine the four aforementioned lemmas to prove Theorem 3.14 which provides an oscillation bound for the sample paths of the processes $\{\bar\mu^{r,m}(\cdot)\}$, and yields Corollary 3.16 as a consequence. We begin with the following technical lemma.

LEMMA 3.9. *Assume* (Q.1) *and let* $T, L > 1$ *and* $0 < \eta < 1$ *be given. Let* $\{\varepsilon_n\}_{n=1}^\infty \subset (0,1)$ *be a sequence such that* $\varepsilon_n \downarrow 0$, *as* $n \to \infty$. *Let* $l_n$, $p$, $M_T$, $\kappa_n$, $\gamma_n$, $N_n$, $r_n$ *be the constants, and* $\{B^r_n\}$ *be the events, given by Lemma* 3.8. *Fix* $n \in \mathbb{N}$, $r \in \mathcal{R}$, $m \le \lfloor rT \rfloor$, *and* $t \in [m, m+L]$. *Define the following random times:*

$$\tau_{\varepsilon_n} = \sup\left\{ s \in [0,t] : \langle 1, \bar\mu^r(s)\rangle < \frac{\varepsilon_n}{8} \right\},$$

$$\tau = \max\left\{ 0, \tau_{\varepsilon_n}, t - \frac{8M_T^2}{\varepsilon_n} \right\},$$

*where we define* $\tau_{\varepsilon_n} = -\infty$ *for the supremum of the empty set. Then on* $B^r_n$,

$$\langle \mathbb{1}_{(0,\kappa_n]}(\cdot - \bar S^r_{\tau,t}), \bar\mu^r(\tau)\rangle \le \frac{\varepsilon_n}{4}. \tag{3.100}$$

PROOF. We only consider realizations in $B^r_n$ in the following. We treat each possibility for $\tau$ separately, and suppose first that $\tau = 0$. Then by (3.79) we have

$$\langle \mathbb{1}_{(0,\kappa_n]}(\cdot - \bar S^r_{\tau,t}), \bar\mu^r(\tau)\rangle \le \sup_{x\in\mathbb{R}_+} \langle \mathbb{1}_{[x,x+\kappa_n]}, \bar\mu^r(0)\rangle \le \frac{\varepsilon_n}{4}. \tag{3.101}$$

Next suppose $\tau = \tau_{\varepsilon_n} > 0$. Then there exists a $\tau' \in [(\tau_{\varepsilon_n} - l_n) \vee 0, \tau_{\varepsilon_n}]$ such that $\langle 1, \bar{\mu}^r(\tau')\rangle \le \varepsilon_n/8$. This implies by (3.4) (with $m = 0$) and (3.74) that

$$\begin{aligned}\langle \mathbb{1}_{(0,\kappa_n]}(\cdot - \bar{S}^r_{\tau,t}), \bar{\mu}^r(\tau)\rangle &\le \langle 1, \bar{\mu}^r(\tau)\rangle \\ &\le \langle 1, \bar{\mu}^r(\tau')\rangle + \bar{E}^r(\tau) - \bar{E}^r(\tau') \\ &\le \frac{\varepsilon_n}{8} + \bar{E}^r(\tau' + l_n) - \bar{E}^r(\tau') \\ &\le \frac{\varepsilon_n}{8} + \frac{\varepsilon_n}{16} < \frac{\varepsilon_n}{4}.\end{aligned} \tag{3.102}$$

Lastly, suppose $\tau = t - (8M_T^2/\varepsilon_n)$. Let $\tau'' = t - (4M_T^2/\varepsilon_n)$. Then since $\tau'' > \tau_{\varepsilon_n}$, we have for $s \in [\tau'', t]$ that $\langle 1, \bar{\mu}^r(s)\rangle \ge \varepsilon_n/8$. This implies by (3.76) that

$$\bar{S}^r_{\tau,t} \ge \bar{S}^r_{\tau'',t} = \int_{\tau''}^t \langle 1, \bar{\mu}^r(s)\rangle^{-1}\, ds \ge \frac{t - \tau''}{M_T} \ge \frac{4M_T}{\varepsilon_n}. \tag{3.103}$$

So by (3.103), Markov's inequality and (3.76),

$$\begin{aligned}\langle \mathbb{1}_{(0,\kappa_n]}(\cdot - \bar{S}^r_{\tau,t}), \bar{\mu}^r(\tau)\rangle &\le \langle \mathbb{1}_{[0,\infty)}(\cdot - \bar{S}^r_{\tau,t}), \bar{\mu}^r(\tau)\rangle \\ &\le \langle \mathbb{1}_{[4M_T/\varepsilon_n,\infty)}, \bar{\mu}^r(\tau)\rangle \\ &\le \frac{\varepsilon_n}{4M_T}\langle \chi, \bar{\mu}^r(\tau)\rangle \\ &\le \frac{\varepsilon_n}{4}.\end{aligned} \tag{3.104}$$

□

The next lemma gives, on $B^r_n$, an upper bound for the amount of mass that $\bar{\mu}^{r,m}(t)$ can have concentrated near zero for $t \in [0, L]$.

LEMMA 3.10. *Assume* (Q.1) *and let* $T, L > 1$ *and* $0 < \eta < 1$ *be given. Let* $\{\varepsilon_n\}_{n=1}^\infty \subset (0,1)$ *be a sequence such that* $\varepsilon_n \downarrow 0$, *as* $n \to \infty$. *Let* $l_n$, $p$, $M_T$, $\kappa_n$, $\gamma_n$, $N_n$, $r_n$ *be the constants, and* $\{B^r_n\}$ *be the events, given by Lemma* 3.8. *Then for each* $n \in \mathbb{N}$, $r \in \mathcal{R}$ *and each* $m \le \lfloor rT \rfloor$, *on* $B^r_n$,

$$\sup_{t\in[0,L]} \langle \mathbb{1}_{[0,\kappa_n]}, \bar{\mu}^{r,m}(t)\rangle \le \frac{\varepsilon_n}{2}. \tag{3.105}$$

PROOF. In this proof we only consider realizations in $B^r_n$. Fix $n \in \mathbb{N}$, $r \in \mathcal{R}$, $m \le \lfloor rT \rfloor$ and $t \in [m, m+L]$, and let $\tau$ be defined as in Lemma 3.9. We must show that

$$\langle \mathbb{1}_{[0,\kappa_n]}, \bar{\mu}^r(t)\rangle \le \frac{\varepsilon_n}{2}. \tag{3.106}$$

If $\tau = t$, then the result follows from Lemma 3.9, since $\bar{\mu}^{r,m}(t)$ does not charge $\{0\}$, and $\bar{S}^r_{\tau,\tau} = 0$. Thus, it suffices to consider the case $\tau < t$. Consider an arrival $i$ such that $U_i^r r^{-1} \in (\tau, t]$ and $v_i^r \geq (64M_T^2/\varepsilon_n^2) + 2\kappa_n$. Since $t - \tau \leq 8M_T^2/\varepsilon_n$, and since $\langle 1, \bar{\mu}^r(s)\rangle \geq \varepsilon_n/8$ for $s \in (\tau, t]$, we have

$$\bar{S}^r_{U_i^r r^{-1},t} \leq \bar{S}^r_{\tau,t} = \int_\tau^t \langle 1, \bar{\mu}^r(s)\rangle^{-1}\, ds \leq (t-\tau)\frac{8}{\varepsilon_n} \leq \frac{64M_T^2}{\varepsilon_n^2}. \tag{3.107}$$

Thus, $v_i^r - \bar{S}^r_{U_i^r r^{-1},t} \geq (64M_T^2/\varepsilon_n^2) + 2\kappa_n - (64M_T^2/\varepsilon_n^2) > \kappa_n$. Recalling that $N_n = \lceil 64M_T^2/(\kappa_n \varepsilon_n^2)\rceil + 1$, this implies that for $k \geq N_n + 1$,

$$\mathbb{1}_{[k\kappa_n,(k+1)\kappa_n)}(v_i^r)\mathbb{1}_{(0,\kappa_n]}(v_i^r - \bar{S}^r_{U_i^r r^{-1},t}) = 0. \tag{3.108}$$

Next, consider two jobs $i < j$ for which $v_i^r, v_j^r \in [k\kappa_n, (k+1)\kappa_n)$ for some integer $k \in \{0, \ldots, N_n\}$. If $U_i^r r^{-1}, U_j^r r^{-1} \in (\tau, t]$ and $U_j^r r^{-1} - U_i^r r^{-1} \geq l_n$, then, using the fact that $\langle 1, \bar{\mu}^r(s)\rangle > 0$ for $s \in (\tau, t]$ again, we have at time $t$ that

$$\begin{aligned}
&(v_j^r - \bar{S}^r_{U_j^r r^{-1},t}) - (v_i^r - \bar{S}^r_{U_i^r r^{-1},t})\\
&\quad = \bar{S}^r_{U_i^r r^{-1}, U_j^r r^{-1}} + v_j^r - v_i^r\\
&\quad \geq \frac{U_j^r r^{-1} - U_i^r r^{-1}}{\sup_{s\in[0,t]}\langle 1, \bar{\mu}^r(s)\rangle} - \kappa_n\\
&\quad \geq \frac{l_n}{M_T} - \kappa_n\\
&\quad > 2\kappa_n - \kappa_n = \kappa_n,
\end{aligned}$$

where the last two inequalities are by (3.76) and (3.78), respectively. This implies that at most one of

$$\mathbb{1}_{(0,\kappa_n]}(v_i^r - \bar{S}^r_{U_i^r r^{-1},t}) \quad \text{and} \quad \mathbb{1}_{(0,\kappa_n]}(v_j^r - \bar{S}^r_{U_j^r r^{-1},t})$$

is nonzero. So for each $k \in \{0, \ldots, N_n\}$, all jobs $i$ satisfying $U_i^r r^{-1} \in (\tau, t]$, $v_i^r \in [k\kappa_n, (k+1)\kappa_n)$ and $v_i^r - \bar{S}^r_{U_i^r r^{-1},t} \in (0, \kappa_n]$ must also satisfy $U_i^r r^{-1} \in (s, s + l_n]$, for some $s \in [\tau, (t - l_n) \vee \tau]$. (Note: $s$ is random in general.) This yields the following estimate at time $t$, for each $k \in \{0, \ldots, N_n\}$:

$$\begin{aligned}
&\frac{1}{r}\sum_{i=r\bar{E}^r(\tau)+1}^{r\bar{E}^r(t)} \mathbb{1}_{[k\kappa_n,(k+1)\kappa_n)}(v_i^r)\mathbb{1}_{(0,\kappa_n]}(v_i^r - \bar{S}^r_{U_i^r r^{-1},t})\\
&\quad \leq \sup_{s\in[\tau,(t-l_n)\vee\tau]} \frac{1}{r}\sum_{i=r\bar{E}^r(s)+1}^{r\bar{E}^r(s+l_n)} \mathbb{1}_{[k\kappa_n,(k+1)\kappa_n)}(v_i^r)\mathbb{1}_{(0,\kappa_n]}(v_i^r - \bar{S}^r_{U_i^r r^{-1},t}) \qquad (3.109)\\
&\quad \leq \sup_{m\leq\lfloor rT\rfloor}\sup_{s\in[0,L]} \frac{1}{r}\sum_{i=r\bar{E}^{r,m}(s)+1}^{r\bar{E}^{r,m}(s+l_n)} \mathbb{1}_{[k\kappa_n,(k+1)\kappa_n)}(v_i^r).
\end{aligned}$$

Using (3.2) (with $m=0$), we have that on $B_n^r$,

$$
\begin{aligned}
&\langle \mathbb{1}_{[0,\kappa_n]}, \bar{\mu}^r(t)\rangle \\
&\quad = \langle \mathbb{1}_{(0,\kappa_n]}(\cdot - \bar{S}^r_{\tau,t}), \bar{\mu}^r(\tau)\rangle + \frac{1}{r}\sum_{i=r\bar{E}^r(\tau)+1}^{r\bar{E}^r(t)} \mathbb{1}_{(0,\kappa_n]}(v_i^r - \bar{S}^r_{U_i^r r^{-1},t}) \\
&\quad \le \frac{\varepsilon_n}{4} + \sum_{k=0}^{\infty}\frac{1}{r}\sum_{i=r\bar{E}^r(\tau)+1}^{r\bar{E}^r(t)} \mathbb{1}_{[k\kappa_n,(k+1)\kappa_n)}(v_i^r)\mathbb{1}_{(0,\kappa_n]}(v_i^r - \bar{S}^r_{U_i^r r^{-1},t}) \\
&\quad \le \frac{\varepsilon_n}{4} + \sum_{k=0}^{N_n}\sup_{m\le\lfloor rT\rfloor}\sup_{s\in[0,L]}\frac{1}{r}\sum_{i=r\bar{E}^{r,m}(s)+1}^{r\bar{E}^{r,m}(s+l_n)} \mathbb{1}_{[k\kappa_n,(k+1)\kappa_n)}(v_i^r) \\
&\quad \le \frac{\varepsilon_n}{4} + \sum_{k=0}^{N_n}\frac{\varepsilon_n}{8}\langle \mathbb{1}_{[(k-1/2)\kappa_n,(k+3/2)\kappa_n)}, \nu\rangle \\
&\quad \le \frac{\varepsilon_n}{4} + \frac{\varepsilon_n}{8}\left\langle \sum_{k=0}^{\infty}\mathbb{1}_{[(k-1/2)\kappa_n,(k+3/2)\kappa_n)}, \nu\right\rangle \\
&\quad = \frac{\varepsilon_n}{4} + \frac{\varepsilon_n}{4},
\end{aligned}
$$

where the first inequality is by Lemma 3.9, the second inequality is by (3.108) and (3.109), the third inequality is by (3.83) and the last line follows since $\nu$ is a probability measure. □

The next lemma gives an upper bound on $[0,L]$ for the process $\langle 1, \bar{\mu}^{r,m}(\cdot)\rangle$, on the event that $\langle \chi, \bar{\mu}^{r,m}(0)\rangle$ is below a threshold.

LEMMA 3.11. *Assume* (Q.1) *and let* $T, L > 1$ *and* $0 < \eta < 1$ *be given. Let* $\{\varepsilon_n\}_{n=1}^{\infty} \subset (0,1)$ *be a sequence such that* $\varepsilon_n \downarrow 0$, *as* $n \to \infty$. *Let* $l_n$, $p$, $M_T$, $\kappa_n$, $\gamma_n$, $N_n$, $r_n$ *be the constants, and* $\{B_n^r\}$ *be the events, given by Lemma* 3.8. *For each* $n \in \mathbb{N}$, $r \in \mathcal{R}$ *and* $m \le \lfloor rT\rfloor$, *define the event* $D_{\gamma_n}^{r,m} = \{\langle \chi, \bar{\mu}^{r,m}(0)\rangle \le \gamma_n/2\}$. *Then for each* $n \in \mathbb{N}$, $r \in \mathcal{R}$ *and* $m \le \lfloor rT\rfloor$, *on* $B_n^r \cap D_{\gamma_n}^{r,m}$,

$$
\sup_{t\in[0,L]}\langle 1, \bar{\mu}^{r,m}(t)\rangle \le 2\varepsilon_n. \tag{3.110}
$$

PROOF. Fix $n \in \mathbb{N}$, $r \in \mathcal{R}$, and $m \le \lfloor rT\rfloor$. We must show that on $B_n^r \cap D_{\gamma_n}^{r,m}$,

$$
\sup_{t\in[m,m+L]}\langle 1, \bar{\mu}^{r}(t)\rangle \le 2\varepsilon_n. \tag{3.111}
$$

Let $A^r = B_n^r \cap D_{\gamma_n}^{r,m}$, and let $t_0 = m$, $t_1 = L$, $c = 2\varepsilon_n$ and $l = l_n$. Then condition (i) of Lemma 3.3 holds by (3.74) and condition (ii) holds by (3.80), (3.81) and the definition of $D^{r,m}_{\gamma_n}$. To see that condition (iii) holds, observe that on $B^r_n \cap D^{r,m}_{\gamma_n}$,

$$\begin{aligned} \langle 1, \bar{\mu}^r(m)\rangle &= \langle 1, \bar{\mu}^{r,m}(0)\rangle \\ &= \langle \mathbb{1}_{[0,\kappa_n]}, \bar{\mu}^{r,m}(0)\rangle + \langle \mathbb{1}_{(\kappa_n,\infty)}, \bar{\mu}^{r,m}(0)\rangle \\ &\leq \frac{\varepsilon_n}{2} + \frac{1}{\kappa_n}\langle \chi, \bar{\mu}^{r,m}(0)\rangle \\ &\leq \frac{\varepsilon_n}{2} + \frac{1}{\kappa_n}\frac{\gamma_n}{2} \\ &\leq \varepsilon_n, \end{aligned} \tag{3.112}$$

where the first inequality is by Lemma 3.10 and Markov's inequality, the second uses the definition of $D^{r,m}_{\gamma_n}$ and the last is by (3.80). So the result follows by Lemma 3.3. □

The next lemma provides for each $m \leq \lfloor rT \rfloor$, a lower bound for the process $\langle 1, \bar{\mu}^{r,m}(\cdot)\rangle$ on the event where $\langle \chi, \bar{\mu}^{r,m}(0)\rangle$ is above the threshold $\gamma_n/2$. A consequence of this is an upper bound, on this event, for the rate at which $\bar{S}^{r,m}_{t,t+h}$ can increase as a function of $h$.

Lemma 3.12. *Assume* (Q.1) *and let* $T, L > 1$ *and* $0 < \eta < 1$ *be given. Let* $\{\varepsilon_n\}_{n=1}^{\infty} \subset (0,1)$ *be a sequence such that* $\varepsilon_n \downarrow 0$, *as* $n \to \infty$. *Let* $l_n$, $p$, $M_T$, $\kappa_n$, $\gamma_n$, $N_n$, $r_n$ *be the constants, and* $\{B^r_n\}$ *be the events, given by Lemma* 3.8. *For each* $n \in \mathbb{N}$, $r \in \mathcal{R}$ *and* $m \leq \lfloor rT \rfloor$, *let* $\breve{D}^{r,m}_{\gamma_n}$ *be the complement of* $D^{r,m}_{\gamma_n}$, *that is,* $\breve{D}^{r,m}_{\gamma_n} = \{\langle \chi, \bar{\mu}^{r,m}(0)\rangle > \gamma_n/2\}$. *Then for each* $n \in \mathbb{N}$, *there exists* $\Gamma_n > 0$ *such that for all* $r \in \mathcal{R}$ *and* $m \leq \lfloor rT \rfloor$, *on* $B^r_n \cap \breve{D}^{r,m}_{\gamma_n}$,

$$\inf_{t\in[0,L]} \langle 1, \bar{\mu}^{r,m}(t)\rangle \geq \frac{1}{\Gamma_n} \tag{3.113}$$

*and*

$$\sup_{t\in[0,L-h]} \bar{S}^{r,m}_{t,t+h} \leq h\Gamma_n \tag{3.114}$$

*for any* $0 < h < L$.

Proof. Fix $n \in \mathbb{N}$, $r \in \mathcal{R}$ and $m \leq \lfloor rT \rfloor$. We have for any $K > 0$, on $B^r_n \cap \breve{D}^{r,m}_{\gamma_n}$,

$$\begin{aligned} \frac{\gamma_n}{4} &\leq \inf_{t\in[0,L]} \langle \chi, \bar{\mu}^{r,m}(t)\rangle \\ &= \inf_{t\in[0,L]} (\langle \chi \mathbb{1}_{[0,K]}, \bar{\mu}^{r,m}(t)\rangle + \langle \chi \mathbb{1}_{(K,\infty)}, \bar{\mu}^{r,m}(t)\rangle) \end{aligned}$$

$$\leq \inf_{t\in[0,L]} (K\langle \mathbb{1}_{[0,K]}, \bar{\mu}^{r,m}(t)\rangle + K^{-p}\langle \chi^{1+p}, \bar{\mu}^{r,m}(t)\rangle)$$
$$\leq \inf_{t\in[0,L]} K\langle 1, \bar{\mu}^{r,m}(t)\rangle + K^{-p}M_T,$$

where the first inequality is by (3.81) and the definition of $\breve{D}^{r,m}_{\gamma_n}$, and the third is by (3.76). So for any $K>0$, on $B^r_n \cap \breve{D}^{r,m}_{\gamma_n}$,

$$\inf_{t\in[0,L]} \langle 1, \bar{\mu}^{r,m}(t)\rangle \geq \frac{\gamma_n/4 - K^{-p}M_T}{K}. \tag{3.115}$$

Setting $K = K_n$ large enough so that $K_n^{-p}M_T < \gamma_n/5$, and letting $\Gamma_n = K_n(\frac{\gamma_n}{4} - \frac{\gamma_n}{5})^{-1}$ proves (3.113). To prove (3.114), we have by (3.113) on $B^r_n \cap \breve{D}^{r,m}_{\gamma_n}$,
for $0<h<L$,

$$\sup_{t\in[0,L-h]} \bar{S}^{r,m}_{t,t+h} \leq h \sup_{t\in[0,L]} \varphi(\langle 1, \bar{\mu}^{r,m}(t)\rangle)$$
$$\leq h\Big(\inf_{t\in[0,L]} \langle 1, \bar{\mu}^{r,m}(t)\rangle\Big)^{-1} \leq h\Gamma_n. \qquad \square$$

We are now ready to apply the previous four lemmas to prove an oscillation bound on $B^r_n$ for the sample paths of $\bar{\mu}^{r,m}(\cdot)$, $m \leq \lfloor rT\rfloor$. We use Lemma 3.11 to prove it on $B^r_n \cap D^{r,m}_{\gamma_n}$ and use Lemmas 3.10 and 3.12 to prove it on $B^r_n \cap \breve{D}^{r,m}_{\gamma_n}$. Recall that $\mathbf{D}_L(\mathcal{M}_{\mathrm{F}})$ is the Skorohod space of r.c.l.l. functions defined on $[0,L]$ and taking values in $\mathcal{M}_{\mathrm{F}}$. In order to state the theorem, we define a modulus of continuity on $\mathbf{D}_L(\mathcal{M}_{\mathrm{F}})$, with respect to the metric $\mathbf{d}[\cdot,\cdot]$ as follows.

DEFINITION 3.13. For any $L>1$, $\zeta(\cdot) \in \mathbf{D}_L(\mathcal{M}_{\mathrm{F}})$ and $\delta>0$, define

$$\mathbf{w}_L(\zeta(\cdot),\delta) = \sup_{t\in[0,L-\delta]} \sup_{h\in[0,\delta]} \mathbf{d}[\zeta(t+h), \zeta(t)]. \tag{3.116}$$

THEOREM 3.14 (Oscillation bound). *Assume* (Q.1) *and let* $T, L > 1$ *and* $0<\eta<1$ *be given. Let* $\{\varepsilon_n\}_{n=1}^\infty \subset (0,1)$ *be a sequence such that* $\varepsilon_n \downarrow 0$, *as* $n\to\infty$. *Let* $l_n$, $p$, $M_T$, $\kappa_n$, $\gamma_n$, $N_n$, $r_n$ *be the constants, and* $\{B^r_n\}$ *be the events, given by Lemma* 3.8. *Then for each* $\varepsilon>0$ *there exists* $\delta>0$ *and* $n_\varepsilon \in \mathbb{N}$ *such that for all* $n > n_\varepsilon$ *and all* $r\in\mathcal{R}$, *on* $B^r_n$,

$$\sup_{m\leq\lfloor rT\rfloor} \mathbf{w}_L(\bar{\mu}^{r,m}(\cdot),\delta) \leq \varepsilon. \tag{3.117}$$

PROOF. Fix $\varepsilon>0$. Choose $k_\varepsilon \in \mathbb{N}$ large enough so that $\sum_{k=k_\varepsilon+1}^\infty 2^{-k} \leq \varepsilon/4$, and let $K_\varepsilon = (\max_{k=1}^{k_\varepsilon} \|(g^{\mathcal{G}}_k)'\|_\infty \vee 2)$. Recall that $\mathcal{G} = \{g^{\mathcal{G}}_k\}_{k=1}^\infty \cup \{h^{\mathcal{G}}_k\}_{k=1}^\infty$

is the set of functions used to define the metric $\mathbf{d}[\cdot,\cdot]$ [see (1.4)]. Let $\mathcal{G}^\varepsilon = \{g_k^{\mathcal{G}}\}_{k=1}^{k_\varepsilon} \cup \{h_k^{\mathcal{G}}\}_{k=1}^{\infty}$ so that for any $f \in \mathcal{G}^\varepsilon$, $\|f\|_\infty \le 1$ and $\|f'\|_\infty \le K_\varepsilon$. Choose $n_\varepsilon \in \mathbb{N}$ large enough so that $\varepsilon_{n_\varepsilon} \le \varepsilon/(8k_\varepsilon)$, and let

$$\delta = \min\left\{\frac{L}{2}, l_{n_\varepsilon}, \frac{\varepsilon}{8\Gamma_{n_\varepsilon} M_T k_\varepsilon (K_\varepsilon \vee 1)}, \frac{\kappa_{n_\varepsilon}}{\Gamma_{n_\varepsilon}}, 1\right\}. \tag{3.118}$$

Fix $n > n_\varepsilon$, $r \in \mathcal{R}$ and $m \le \lfloor rT \rfloor$, and consider the event $D_{\gamma_{n_\varepsilon}}^{r,m}$. On $B_n^r \cap D_{\gamma_{n_\varepsilon}}^{r,m}$, we have for all $f \in \mathcal{G}^\varepsilon$,

$$\begin{aligned}
&\sup_{t\in[0,L-\delta]} \sup_{h\in[0,\delta]} |\langle f, \bar{\mu}^{r,m}(t+h)\rangle - \langle f, \bar{\mu}^{r,m}(t)\rangle| \\
&\quad \le 2 \sup_{t\in[0,L]} |\langle f, \bar{\mu}^{r,m}(t)\rangle| \le 2\|f\|_\infty \sup_{t\in[0,L]} \langle 1, \bar{\mu}^{r,m}(t)\rangle \\
&\quad \le 4\varepsilon_{n_\varepsilon} \le \frac{\varepsilon}{2k_\varepsilon}.
\end{aligned} \tag{3.119}$$

The third inequality above follows by Lemma 3.11, since $B_n^r \cap D_{\gamma_{n_\varepsilon}}^{r,m} \subset B_{n_\varepsilon}^r \cap D_{\gamma_{n_\varepsilon}}^{r,m}$ for $n > n_\varepsilon$ [see Lemma 3.8(ii)], and by the fact that $\|f\|_\infty \le 1$ by definition. The last inequality above follows by choice of $n_\varepsilon$. We must show that the above estimate also holds on $B_n^r \cap \breve{D}_{\gamma_{n_\varepsilon}}^{r,m}$ for all $f \in \mathcal{G}^\varepsilon$. First, observe that on $B_n^r \cap \breve{D}_{\gamma_{n_\varepsilon}}^{r,m}$, a first-order Taylor expansion of $f$ gives the following estimate for all $0 < h < L$, $t \in [0, L-h]$ and $y \in (\bar{S}_{t,t+h}^{r,m}, \infty)$:

$$|f(y - \bar{S}_{t,t+h}^{r,m}) - f(y)| = |-\bar{S}_{t,t+h}^{r,m} f'(w_y)| \le h\Gamma_{n_\varepsilon} \|f'\|_\infty \le h\Gamma_{n_\varepsilon} K_\varepsilon, \tag{3.120}$$

for some $w_y \in [y - \bar{S}_{t,t+h}^{r,m}, y]$, where the first inequality follows by Lemma 3.12, since $B_n^r \cap \breve{D}_{\gamma_{n_\varepsilon}}^{r,m} \subset B_{n_\varepsilon}^r \cap \breve{D}_{\gamma_{n_\varepsilon}}^{r,m}$ for $n > n_\varepsilon$, and the second inequality follows by definition of $\mathcal{G}^\varepsilon$. Now subtracting $\langle f, \bar{\mu}^{r,m}(t)\rangle$ from both sides of equation (3.2) and using the fact that $(\mathbb{1}_{(0,\infty)} f)(\cdot - \bar{S}_{t,t+h}^{r,m}) = \mathbb{1}_{(\bar{S}_{t,t+h}^{r,m},\infty)}(\cdot) f(\cdot - \bar{S}_{t,t+h}^{r,m})$ for $t \in [0, L-h]$ yields that on $B_n^r \cap \breve{D}_{\gamma_{n_\varepsilon}}^{r,m}$,

$$\begin{aligned}
&|\langle f, \bar{\mu}^{r,m}(t+h)\rangle - \langle f, \bar{\mu}^{r,m}(t)\rangle| \\
&\quad = \Bigg| \langle \mathbb{1}_{(\bar{S}_{t,t+h}^{r,m},\infty)}(\cdot)(f(\cdot - \bar{S}_{t,t+h}^{r,m}) - f(\cdot)), \bar{\mu}^{r,m}(t)\rangle \\
&\qquad - \langle \mathbb{1}_{[0,\bar{S}_{t,t+h}^{r,m}]} f, \bar{\mu}^{r,m}(t)\rangle \\
&\qquad + \frac{1}{r} \sum_{i=r\bar{E}^{r,m}(t)+1}^{r\bar{E}^{r,m}(t+h)} (\mathbb{1}_{(0,\infty)} f)(v_i^r - \bar{S}_{U_i^r r^{-1}-m,t+h}^{r,m}) \Bigg| \\
&\quad \le \langle |\mathbb{1}_{(\bar{S}_{t,t+h}^{r,m},\infty)}(\cdot)(f(\cdot - \bar{S}_{t,t+h}^{r,m}) - f(\cdot))|, \bar{\mu}^{r,m}(t)\rangle \\
&\qquad + \|f\|_\infty \langle \mathbb{1}_{[0,h\Gamma_{n_\varepsilon}]}, \bar{\mu}^{r,m}(t)\rangle + \|f\|_\infty (\bar{E}^{r,m}(t+h) - \bar{E}^{r,m}(t))
\end{aligned} \tag{3.121}$$

$$\le h\Gamma_{n_\varepsilon} K_\varepsilon \langle 1, \bar{\mu}^{r,m}(t)\rangle + \langle \mathbb{1}_{[0,h\Gamma_{n_\varepsilon}]}, \bar{\mu}^{r,m}(t)\rangle + (\bar{E}^{r,m}(t+h) - \bar{E}^{r,m}(t)),$$

where the first inequality is by Lemma 3.12 (again using the fact that $B_n^r \cap \breve{D}_{\gamma_{n_\varepsilon}}^{r,m} \subset B_{n_\varepsilon}^r \cap \breve{D}_{\gamma_{n_\varepsilon}}^{r,m}$ for $n > n_\varepsilon$) and the second inequality is by (3.120) and the fact that $\|f\|_\infty \le 1$. Taking the supremum over $h \in [0,\delta]$ and $t \in [0, L-\delta]$, we see that on $B_n^r \cap \breve{D}_{\gamma_{n_\varepsilon}}^{r,m}$,

$$\begin{aligned}
&\sup_{t\in[0,L-\delta]} \sup_{h\in[0,\delta]} |\langle f, \bar{\mu}^{r,m}(t+h)\rangle - \langle f, \bar{\mu}^{r,m}(t)\rangle| \\
&\quad\le \sup_{t\in[0,L-\delta]} \Big(\delta\Gamma_{n_\varepsilon} K_\varepsilon \langle 1, \bar{\mu}^{r,m}(t)\rangle \\
&\qquad\qquad + \langle \mathbb{1}_{[0,\delta\Gamma_{n_\varepsilon}]}, \bar{\mu}^{r,m}(t)\rangle + (\bar{E}^{r,m}(t+\delta) - \bar{E}^{r,m}(t))\Big) \\
&\quad\le \sup_{t\in[0,L-\delta]} \left(\frac{\varepsilon}{8M_T k_\varepsilon}\langle 1, \bar{\mu}^{r,m}(t)\rangle \right.\\
&\qquad\qquad \left. + \langle \mathbb{1}_{[0,\kappa_{n_\varepsilon}]}, \bar{\mu}^{r,m}(t)\rangle + (\bar{E}^{r,m}(t+l_{n_\varepsilon}) - \bar{E}^{r,m}(t))\right) \\
&\quad\le \frac{\varepsilon}{8M_T k_\varepsilon} M_T + \frac{\varepsilon_{n_\varepsilon}}{2} + \frac{\varepsilon_{n_\varepsilon}}{16} \\
&\quad\le \frac{\varepsilon}{8k_\varepsilon} + \frac{\varepsilon}{16k_\varepsilon} + \frac{\varepsilon}{128k_\varepsilon} \\
&\quad< \frac{\varepsilon}{2k_\varepsilon},
\end{aligned} \tag{3.122}$$

where the second inequality is by (3.118) and the third is by (3.76), Lemma 3.10, and (3.74) (once again using the fact that $B_n^r \cap \breve{D}_{\gamma_{n_\varepsilon}}^{r,m} \subset B_{n_\varepsilon}^r \cap \breve{D}_{\gamma_{n_\varepsilon}}^{r,m}$ for $n > n_\varepsilon$). So the desired estimate (3.119) holds on both $B_n^r \cap D_{\gamma_{n_\varepsilon}}^{r,m}$ and $B_n^r \cap \breve{D}_{\gamma_{n_\varepsilon}}^{r,m}$ and therefore on $B_n^r$ for all $n > n_\varepsilon$, $r \in \mathcal{R}$ and $m \le \lfloor rT \rfloor$. Thus, combining (3.119) and (3.122) with the definitions (3.116), (1.4) and the definition of $\mathcal{G}^\varepsilon$, we have on $B_n^r$, for all $n > n_\varepsilon$ and $r \in \mathcal{R}$,

$$\begin{aligned}
&\sup_{m\le\lfloor rT\rfloor} \mathbf{w}_L(\bar{\mu}^{r,m}(\cdot),\delta) \\
&\quad= \sup_{m\le\lfloor rT\rfloor} \sup_{t\in[0,L-\delta]} \sup_{h\in[0,\delta]} \mathbf{d}[\bar{\mu}^{r,m}(t+h), \bar{\mu}^{r,m}(t)] \\
&\quad\le \sup_{m\le\lfloor rT\rfloor} \sup_{t\in[0,L-\delta]} \sup_{h\in[0,\delta]} \left(\sum_{k=1}^{k_\varepsilon} 2^{-k}(|\langle g_k^{\mathcal{G}}, \bar{\mu}^{r,m}(t+h)\rangle - \langle g_k^{\mathcal{G}}, \bar{\mu}^{r,m}(t)\rangle| \wedge 1)\right. \\
&\qquad\qquad + \sum_{k=k_\varepsilon+1}^{\infty} 2^{-k} + \sup_{k\in\mathbb{N}} |\langle h_k^{\mathcal{G}}, \bar{\mu}^{r,m}(t+h)\rangle
\end{aligned}$$

$$- \langle h_k^{\mathcal{G}}, \bar{\mu}^{r,m}(t)\rangle|\Bigg)$$

$$\leq \frac{k_\varepsilon}{2}\frac{\varepsilon}{2k_\varepsilon} + \frac{\varepsilon}{4} + \frac{\varepsilon}{2k_\varepsilon}$$

$$\leq \varepsilon. \qquad \square$$

The preceding theorem has the following important consequence, which will be used to identify certain measure valued paths which approximate the sample paths of the fluid scaled state descriptors $\{\bar{\mu}^{r,m}(\cdot)\}$ on the events $B_n^r$. This will be the first important step in our state space collapse argument. We will need the following definition.

DEFINITION 3.15. Assume (Q.1) and let $T, L > 1$ and $0 < \eta < 1$ be given. Let $\{\varepsilon_n\}_{n=1}^{\infty} \subset (0,1)$ be a sequence such that $\varepsilon_n \downarrow 0$, as $n \to \infty$. Let $l_n$, $p$, $M_T$, $\kappa_n$, $\gamma_n$, $N_n$, $r_n$ be the constants, and $\{B_n^r\}$ be the events given by Lemma 3.8. Then for each $r \in \mathcal{R}$, we define

$$n(r) = \begin{cases} \sup\{n \in \mathbb{N} : r_n < r\}, & r > r_1, \\ 1, & r \leq r_1, \end{cases} \tag{3.123}$$

and we define

$$\mathscr{B}_L^r = \{\zeta(\cdot) \in \mathbf{D}_L(\mathcal{M}_{\mathrm{F}}) : \zeta(\cdot) \equiv \bar{\mu}^{r,m}(\cdot)(\omega) \text{ on } [0, L], \\ \text{for some } \omega \in B_{n(r)}^r \text{ and } m \leq \lfloor rT \rfloor\}. \tag{3.124}$$

Note that since $r_n \to \infty$ as $n \to \infty$ (see Lemma 3.8), $n(r)$ is finite for each $r \in \mathcal{R}$. Note also that for $r > r_1$, $\mathscr{B}_L^r$ is nonempty since by Lemma 3.8, $\mathbf{P}^r(B_{n(r)}^r) \geq 1 - \eta$.

COROLLARY 3.16 (Relative compactness). *Assume* (Q.1) *and let* $T, L > 1$ *and* $0 < \eta < 1$ *be given. Let* $\{\varepsilon_n\}_{n=0}^{\infty} \subset (0,1)$ *be a sequence such that* $\varepsilon_n \downarrow 0$, *as* $n \to \infty$. *Let* $\tilde{\mathcal{R}} \subset \mathcal{R}$ *be a subsequence and suppose* $\{\zeta^{\tilde{r}}(\cdot)\}_{\tilde{r} \in \tilde{\mathcal{R}}} \subset \mathbf{D}_L(\mathcal{M}_{\mathrm{F}})$ *is a sequence of paths such that for each* $\tilde{r} > r_1$, $\zeta^{\tilde{r}}(\cdot) \in \mathscr{B}_L^{\tilde{r}}$. *Then* $\{\zeta^{\tilde{r}}(\cdot)\}_{\tilde{r} \in \tilde{\mathcal{R}}}$ *is relatively compact in* $\mathbf{D}_L(\mathcal{M}_{\mathrm{F}})$. *Moreover, any limit point* $\zeta(\cdot)$ *of the sequence is continuous.*

PROOF. Define

$$\tilde{C}_T = \{\zeta \in \mathcal{M}_{\mathrm{F}} : \langle 1, \zeta\rangle \vee \langle \chi, \zeta\rangle \leq M_T\}.$$

Since $\langle \chi, \zeta\rangle \leq M_T$ implies $\langle \mathbb{1}_{[K,\infty)}, \zeta\rangle \leq M_T/K$ for any $K > 0$, we have

$$\sup_{\zeta \in \tilde{C}_T} \langle \mathbb{1}_{[K,\infty)}, \zeta\rangle \to 0 \qquad \text{as } K \to \infty,$$

which implies that $\tilde{C}_T \subset \mathcal{M}_{\mathrm{F}}$ is relatively compact ([10], Theorem A7.5). Let $C_T$ be the closure of $\tilde{C}_T$ and observe that by Definition 3.15 and (3.76), $\zeta^{\tilde{r}}(t) \in C_T$ for all $\tilde{r} > r_1$ and $t \in [0, L]$. Since $n(\tilde{r}) \to \infty$ as $\tilde{r} \to \infty$, Theorem 3.14 implies that

$$\lim_{\delta \to 0} \limsup_{\tilde{r} \to \infty} \mathbf{w}_L(\zeta^{\tilde{r}}(\cdot), \delta) = 0. \tag{3.125}$$

Thus, the relative compactness of $\{\zeta^{\tilde{r}}(\cdot)\}_{\tilde{r} \in \tilde{\mathcal{R}}}$ follows from [4], Chapter 3, Theorem 6.3, by noting that the modulus of continuity used there is bounded above by $\mathbf{w}_L(\cdot, \cdot)$, and that the result stated there still holds if one replaces $\mathbf{D}_\infty(\mathcal{M}_{\mathrm{F}})$ by $\mathbf{D}_L(\mathcal{M}_{\mathrm{F}})$ and $T > 0$ there by $L$. By (3.125) and the definition of $\mathbf{w}_L(\cdot, \cdot)$, we see that any limit point $\zeta(\cdot)$ of the sequence must be continuous. □

**4. Diffusion limit.** Having established the desired relative compactness property for sample paths of the processes $\{\bar{\mu}^{r,m}(\cdot)\}$, we now turn to the second main task in our strategy for proving Theorem 2.3. In Section 4.1, we identify limit points of the sequences specified in Corollary 3.16 as *fluid model solutions on* $[0, L]$ (see Definition 4.2 and Lemma 4.3). We then show that for large $r$, sample paths of the processes $\{\bar{\mu}^{r,m}(\cdot)\}$, can be uniformly approximated by fluid model solutions which are arbitrarily close to being in steady state (see Lemma 4.4 and Proposition 4.5). This leads us to the proof of state space collapse and finally to the proof of Theorem 2.3, both of which are contained in Section 4.2.

4.1. *Local fluid approximations.*

DEFINITION 4.1. Assume (Q.1) and let $T, L > 1$ and $0 < \eta < 1$ be given. Let $\{\varepsilon_n\}_{n=1}^{\infty} \subset (0, 1)$ be a sequence such that $\varepsilon_n \downarrow 0$, as $n \to \infty$. Define $\mathscr{B}_L$ to be the set of all $\zeta(\cdot) \in \mathbf{D}_L(\mathcal{M}_{\mathrm{F}})$ such that there exists a subsequence $\tilde{\mathcal{R}} \subset \mathcal{R}$ and a sequence $\{\zeta^{\tilde{r}}(\cdot)\}_{\tilde{r} \in \tilde{\mathcal{R}}} \subset \mathbf{D}_L(\mathcal{M}_{\mathrm{F}})$ satisfying $\zeta^{\tilde{r}}(\cdot) \in \mathscr{B}_L^{\tilde{r}}$ for each $\tilde{r} \in \tilde{\mathcal{R}}$, and

$$\zeta^{\tilde{r}}(\cdot) \xrightarrow{J_1} \zeta(\cdot) \qquad \text{as } \tilde{r} \to \infty. \tag{4.1}$$

DEFINITION 4.2. Given $1 < L < \infty$, a fluid model solution on $[0, L]$ for the critical data $(\alpha, \nu)$ is a function $\zeta(\cdot) : [0, L] \longrightarrow \mathcal{M}_{\mathrm{F}}$ such that:

(i) $\zeta(\cdot)$ is continuous,
(ii) $\langle \mathbb{1}_{\{0\}}, \zeta(t) \rangle = 0$ for all $t \in [0, L]$,
(iii) if $\zeta(0) \neq \mathbf{0}$, then $\zeta(t) \neq \mathbf{0}$ for all $t \in [0, L]$ and $\zeta(\cdot)$ satisfies

$$\langle g, \zeta(t) \rangle = \langle g, \zeta(0) \rangle - \int_0^t \frac{\langle g', \zeta(s) \rangle}{\langle 1, \zeta(s) \rangle} \, ds + \alpha t \langle g, \nu \rangle, \tag{4.2}$$

for all $t \in [0, L]$ and all $g \in \mathcal{C} = \{g \in \mathbf{C}_b^1(\mathbb{R}_+) : g(0) = 0, \, g'(0) = 0\}$ and

(iv) if $\zeta(0)=\mathbf{0}$, then $\zeta(t)=\mathbf{0}$ for all $t\in[0,L]$.

We note that Definition 4.2 differs slightly from the definition of a fluid model solution for critical data $(\alpha,\nu)$ given in [6], Section 3.1. Besides the fact that we only consider fluid model solutions defined over finite time intervals here, the time $t^*$ used in [6] (the first time at which a fluid model solution reaches the zero measure) is not present in our definition. In fact, it is shown in [6] that if $\zeta(0)\neq\mathbf{0}$, then $t^*=\infty$. Indeed, if we replace $[0,L]$ by $[0,\infty)$ in Definition 4.2 above, we obtain an equivalent definition to that given in [6]. Similarly, by restricting a fluid model solution as defined in [6] to the finite interval $[0,L]$, one obtains a fluid model solution on $[0,L]$ as defined here. The next lemma asserts that the elements of $\mathscr{B}_L$ are in fact fluid model solutions on $[0,L]$.

LEMMA 4.3. *Assume* (Q.1) *and let* $T>1$ *and* $0<\eta<1$ *be given. Let* $\{\varepsilon_n\}_{n=1}^{\infty}\subset(0,1)$ *be a sequence such that* $\varepsilon_n\downarrow 0$, *as* $n\to\infty$. *Then there exists* $q\in(0,p)$ *and* $1<M_q<\infty$ *such that for any* $L>1$ *and any* $\zeta(\cdot)\in\mathscr{B}_L$,

(i) $\sup_{t\in[0,L]}\langle 1,\zeta(t)\rangle\vee\langle\chi,\zeta(t)\rangle\vee\langle\chi^{1+q},\zeta(t)\rangle\le M_q$,
(ii) $\langle\chi,\zeta(\cdot)\rangle$ *is constant on* $[0,L]$ *and*
(iii) $\zeta(\cdot)$ *is a fluid model solution on* $[0,L]$ *for the critical data* $(\alpha,\nu)$.

PROOF. Fix $L>1$ and let $l_n$, $p$, $M_T$, $\kappa_n$, $\gamma_n$, $N_n$, $r_n$ be the constants and $\{B_n^r\}$ be the events given by Lemma 3.8. Fix $\zeta(\cdot)\in\mathscr{B}_L$ and let $\tilde{\mathcal{R}}\subset\mathcal{R}$ be a subsequence and $\{\zeta^{\tilde r}(\cdot)\}_{\tilde r\in\tilde{\mathcal{R}}}\subset\mathbf{D}_L(\mathcal{M}_{\mathrm{F}})$ be a sequence such that $\zeta^{\tilde r}(\cdot)\in\mathscr{B}_L^{\tilde r}$ for each $\tilde r\in\tilde{\mathcal{R}}$ and $\zeta^{\tilde r}(\cdot)\xrightarrow{J_1}\zeta(\cdot)$ as $\tilde r\to\infty$. Recall that for each $\tilde r\in\tilde{\mathcal{R}}$, $\zeta^{\tilde r}(\cdot)$ is a realization of $\bar\mu^{\tilde r,m}(\cdot)$ on $B_{n(\tilde r)}^{\tilde r}$ for some $m\le\lfloor\tilde r T\rfloor$. This fact will be used throughout the proof. By Corollary 3.16, $\zeta(\cdot)$ is continuous, so as $\tilde r\to\infty$, we have

$$\|\mathbf{d}[\zeta^{\tilde r}(\cdot),\zeta(\cdot)]\|_L\to 0. \tag{4.3}$$

By (3.76),

$$\limsup_{\tilde r\to\infty}\sup_{t\in[0,L]}\langle 1,\zeta^{\tilde r}(t)\rangle\vee\langle\chi,\zeta^{\tilde r}(t)\rangle\vee\langle\chi^{1+p},\zeta^{\tilde r}(t)\rangle\le M_T. \tag{4.4}$$

Combined with (4.3), this implies that

$$\sup_{t\in[0,L]}\langle 1,\zeta(t)\rangle\le M_T. \tag{4.5}$$

Similarly, combining (4.3) and (4.4) with two applications of Lemma 3.5 (with $q'=p$ and $M=M_T$) we see that there exists $q\in(0,p)$ and $M_q\in[M_T,\infty)$ such that

$$\sup_{t\in[0,L]}\langle\chi,\zeta(t)\rangle\vee\langle\chi^{1+q},\zeta(t)\rangle\le M_q, \tag{4.6}$$

and such that for each $t \in [0, L]$,

$$|\langle \chi, \zeta^{\tilde{r}}(t)\rangle - \langle \chi, \zeta(t)\rangle| \to 0 \qquad \text{as } \tilde{r} \to \infty. \tag{4.7}$$

Together, (4.5) and (4.6) imply (i) above. Notice that the constant $M_q$ given for $p$ and $M_T$ by Lemma 3.5 does not depend on $L$, and that $M_q \ge M_T > 1$ [see (3.51) and Lemma 3.6].

Next, we have for any $t \in [0, L]$,

$$\begin{aligned} |\langle \chi, \zeta(t)\rangle - \langle \chi, \zeta(0)\rangle| &\le \liminf_{\tilde{r}\to\infty} (|\langle \chi, \zeta(t)\rangle - \langle \chi, \zeta^{\tilde{r}}(t)\rangle| \\ &\qquad + |\langle \chi, \zeta^{\tilde{r}}(t)\rangle - \langle \chi, \zeta^{\tilde{r}}(0)\rangle| + |\langle \chi, \zeta^{\tilde{r}}(0)\rangle - \langle \chi, \zeta(0)\rangle|) \\ &\le \liminf_{\tilde{r}\to\infty} \frac{\gamma_{n(\tilde{r})}}{4} \\ &= 0, \end{aligned} \tag{4.8}$$

where the second inequality is by (4.7) for the first and third terms and (3.81) for the second term. The last line follows since $n(\tilde{r}) \to \infty$, as $\tilde{r} \to \infty$, and $\gamma_n \to 0$, as $n \to \infty$ (see Lemma 3.8). Thus $\langle \chi, \zeta(\cdot)\rangle$ is constant on $[0, L]$, which proves (ii) above.

We now show that $\zeta(\cdot)$ is a fluid model solution on $[0, L]$. We have already verified property (i) of Definition 4.2. For any fixed $n \in \mathbb{N}$,

$$\begin{aligned} \sup_{t\in[0,L]} \langle \mathbb{1}_{\{0\}}, \zeta(t)\rangle &\le \sup_{t\in[0,L]} \langle \mathbb{1}_{[0,\kappa_n)}, \zeta(t)\rangle \\ &\le \sup_{t\in[0,L]} \limsup_{\tilde{r}\to\infty} \langle \mathbb{1}_{[0,\kappa_n)}, \zeta^{\tilde{r}}(t)\rangle \\ &\le \limsup_{\tilde{r}\to\infty} \sup_{t\in[0,L]} \langle \mathbb{1}_{[0,\kappa_n]}, \zeta^{\tilde{r}}(t)\rangle \\ &\le \frac{\varepsilon_n}{2}, \end{aligned} \tag{4.9}$$

where the second inequality follows by (4.3) and the Portmanteau theorem ([1], Theorem 2.1). The last inequality follows by Lemma 3.10, since $n(\tilde{r}) \to \infty$, as $\tilde{r} \to \infty$, which implies that $B^{\tilde{r}}_{n(\tilde{r})} \subset B^{\tilde{r}}_n$ for all sufficiently large $\tilde{r}$. Since $n \in \mathbb{N}$ is arbitrary and $\varepsilon_n \downarrow 0$ as $n \to \infty$, this proves property (ii) of Definition 4.2. It remains to verify properties (iii) and (iv) there. We must show that either $\zeta(\cdot) \equiv \mathbf{0}$ on $[0, L]$, or that for all $t \in [0, L]$, $\zeta(t) \ne \mathbf{0}$ and (4.2) holds. Since we have shown that $\langle \chi, \zeta(\cdot)\rangle$ is constant on $[0, L]$ and that property (ii) of Definition 4.2 holds, it suffices to show that (4.2) holds when $\zeta(0) \ne \mathbf{0}$. For this, suppose that $\zeta(0) \ne \mathbf{0}$ and note that property (ii) of Definition 4.2 implies that $\langle\chi,\zeta(t)\rangle=\langle\chi,\zeta(0)\rangle>0$ for all $t\in[0,L]$. So since $\zeta(\cdot)$ is continuous,

$$\inf_{t\in[0,L]}\langle 1,\zeta(t)\rangle>0. \tag{4.10}$$

To show (4.2), we follow the analogous proof in ([6], Section 5.3) closely. We first restrict our attention to the case $g\in\mathcal{V}$ (see Definition 3.7) and derive a prelimit version of (4.2) which is satisfied by $\zeta^{\tilde r}(\cdot)$ for all sufficiently large $\tilde r$. We will then pass to the limit in this relation to obtain (4.2) for $\zeta(\cdot)$ and $g\in\mathcal{V}$. Finally, we make a simple extension from $\mathcal{V}$ to $\mathcal{C}$.

Let $g\in\mathcal{V}$ and note that since $\mathcal{V}\subset\tilde{\mathcal{C}}$, $g'$ is uniformly continuous on $\mathbb{R}$ [recall that we extend $g$ to be identically zero on $(-\infty,0)$]. So, there is a continuous nondecreasing function $\psi_g:[0,\infty)\longrightarrow[0,\infty)$ with $\psi_g(0)=0$, such that for any $h\in\mathbb{R}$,

$$\sup_{x\in\mathbb{R}}\|g'(x+h)-g'(x)\|\le\psi_g(|h|). \tag{4.11}$$

Recall that for each $\tilde r\in\tilde{\mathcal{R}}$, $\zeta^{\tilde r}(\cdot)$ is a realization of $\bar\mu^{\tilde r,m}(\cdot)$ for some $m\le\lfloor\tilde rT\rfloor$ and some $\tilde\omega\in B^{\tilde r}_{n(\tilde r)}$. Since the sequence $\{\zeta^{\tilde r}(\cdot)\}_{\tilde r\in\tilde{\mathcal{R}}}$ remains fixed for the remainder of this proof, it is understood that all random objects indexed by $\tilde r\in\tilde{\mathcal{R}}$ and $m\le\lfloor\tilde rT\rfloor$ are evaluated at the particular $\tilde\omega\in B^{\tilde r}_{n(\tilde r)}$ and $m$ for which $\zeta^{\tilde r}(\cdot)$ is a realization. Fix $t\in[0,L]$. For any $N\in\mathbb{N}$ and $j\in\{0,1,\dots,N-1\}$, define $t_j=\frac{jt}{N}$ and $t^j=t_{j+1}$. Then we have for each $\tilde r\in\tilde{\mathcal{R}}$,

$$\begin{aligned}
\langle g,\zeta^{\tilde r}(t)\rangle-\langle g,\zeta^{\tilde r}(0)\rangle&=\sum_{j=0}^{N-1}(\langle g,\zeta^{\tilde r}(t^j)\rangle-\langle g,\zeta^{\tilde r}(t_j)\rangle)\\
&=\sum_{j=0}^{N-1}(\langle g,\zeta^{\tilde r}(t^j)\rangle-\langle g(\cdot-\bar S^{\tilde r,m}_{t_j,t^j}),\zeta^{\tilde r}(t_j)\rangle)\\
&\quad+\sum_{j=0}^{N-1}(\langle g(\cdot-\bar S^{\tilde r,m}_{t_j,t^j}),\zeta^{\tilde r}(t_j)\rangle-\langle g,\zeta^{\tilde r}(t_j)\rangle)\\
&=\sum_{j=0}^{N-1}\frac{1}{\tilde r}\sum_{i=\tilde r\bar E^{\tilde r,m}(t_j)+1}^{\tilde r\bar E^{\tilde r,m}(t^j)}g(v_i^{\tilde r}-\bar S^{\tilde r,m}_{U_i^{\tilde r}\tilde r^{-1}-m,t^j})\\
&\quad+\sum_{j=0}^{N-1}\langle g(\cdot-\bar S^{\tilde r,m}_{t_j,t^j})-g(\cdot),\zeta^{\tilde r}(t_j)\rangle,
\end{aligned} \tag{4.12}$$

where the first term in the last equality is by (3.2) (using $t=t_j$, $h=t^j-t_j$ there) and the fact that $(\mathbb{1}_{(0,\infty)}g)\equiv g$, since $g(0)=0$ for $g\in\tilde{\mathcal{C}}$. By (4.3),

(4.4) and (4.10), we have

$$\limsup_{\tilde{r}\to\infty} \|\langle 1, \zeta^{\tilde{r}}(\cdot)\rangle\|_L \le M_T, \tag{4.13}$$

$$\limsup_{\tilde{r}\to\infty} \|\langle 1, \zeta^{\tilde{r}}(\cdot)\rangle^{-1}\|_L < \infty. \tag{4.14}$$

By (4.13) and (4.14), we can assume for the remainder of the proof that $\tilde{r}$ is large enough so that for some $M_\zeta > 0$,

$$\sup_{s\in[0,L]} \langle 1, \zeta^{\tilde{r}}(s)\rangle \vee \langle 1, \zeta^{\tilde{r}}(s)\rangle^{-1} \le M_\zeta. \tag{4.15}$$

We handle the two terms in (4.12) separately. To begin with, since $g \in \mathcal{V} \subset \tilde{\mathcal{C}}$ has been extended to be an element of $\mathbf{C}_b^1(\mathbb{R})$, we have the following first-order Taylor expansion for each $j \in \{0, 1, \ldots, N-1\}$ and each $x \in \mathbb{R}_+$:

$$g(x - \bar{S}^{\tilde{r},m}_{t_j,t^j}) - g(x) = g'(w_j^x) h_j, \tag{4.16}$$

where $h_j = -\bar{S}^{\tilde{r},m}_{t_j,t^j}$ and $w_j^x \in \mathbb{R}$ is in the interval $[x - \bar{S}^{\tilde{r},m}_{t_j,t^j}, x]$. Note that by (2.20) and (4.15),

$$\max_{j=0,\ldots,N-1} |h_j| = \max_{j=0,\ldots,N-1} |\bar{S}^{\tilde{r},m}_{t_j,t^j}| \le \frac{t}{N} \|\langle 1, \zeta^{\tilde{r}}(\cdot)\rangle^{-1}\|_L \le \frac{tM_\zeta}{N}. \tag{4.17}$$

For each $j \in \{0, \ldots, N-1\}$, let $z_j = \sup_{s\in[t_j,t^j)} \langle 1, \zeta^{\tilde{r}}(s)\rangle^{-1}$ and define $\tilde{h}_j = -z_j \frac{t}{N}$. Then

$$\begin{aligned} \sum_{j=0}^{N-1} |h_j - \tilde{h}_j| &= \sum_{j=0}^{N-1} \left| z_j \frac{t}{N} - \bar{S}^{\tilde{r},m}_{t_j,t^j} \right| \\ &= \sum_{j=0}^{N-1} \left( z_j \frac{t}{N} - \bar{S}^{\tilde{r},m}_{t_j,t^j} \right) \\ &= \sum_{j=0}^{N-1} \left( z_j \frac{t}{N} \right) - \bar{S}^{\tilde{r},m}_{0,t}. \end{aligned} \tag{4.18}$$

For each $N \in \mathbb{N}$ and $s \in [0, t)$, let $f_N(s) = \sum_{j=0}^{N-1} z_j \mathbb{1}_{[t_j,t^j)}(s)$ and define $f_N(t) = 0$. We can make the following estimate for the second term in (4.12):

$$\begin{aligned} &\left| \sum_{j=0}^{N-1} \langle g(\cdot - \bar{S}^{\tilde{r},m}_{t_j,t^j}) - g(\cdot), \zeta^{\tilde{r}}(t_j)\rangle - \sum_{j=0}^{N-1} \langle g'(\cdot)\tilde{h}_j, \zeta^{\tilde{r}}(t_j)\rangle \right| \\ &\quad \le \sum_{j=0}^{N-1} \sup_{x\in\mathbb{R}_+} |g(x - \bar{S}^{\tilde{r},m}_{t_j,t^j}) - g(x) - g'(x)\tilde{h}_j| \langle 1, \zeta^{\tilde{r}}(t_j)\rangle, \end{aligned}$$

$$= \sum_{j=0}^{N-1} \sup_{x\in\mathbb{R}_+} |g'(w_j^x)h_j - g'(x)\tilde{h}_j|\langle 1, \zeta^{\tilde{r}}(t_j)\rangle$$

$$\leq \|\langle 1, \zeta^{\tilde{r}}(\cdot)\rangle\|_L \sum_{j=0}^{N-1} \sup_{x\in\mathbb{R}_+} (|g'(w_j^x) - g'(x)||h_j| + |g'(x)||h_j - \tilde{h}_j|), \tag{4.19}$$

$$\leq M_\zeta\left(N\psi_g\left(\frac{tM_\zeta}{N}\right)\frac{tM_\zeta}{N} + \|g'\|_\infty\left(\sum_{j=0}^{N-1}\left(z_j\frac{t}{N}\right) - \bar{S}_{0,t}^{\tilde{r},m}\right)\right)$$

$$= M_\zeta\left(\psi_g\left(\frac{tM_\zeta}{N}\right)tM_\zeta + \|g'\|_\infty\left(\int_0^t f_N(s)\,ds - \int_0^t \langle 1, \zeta^{\tilde{r}}(s)\rangle^{-1}\,ds\right)\right).$$

In the third line above we have used the Taylor expansion (4.16). The last inequality then follows from (4.15), (4.11), (4.17) and (4.18). The substitution of the second integral in the last line follows by (2.20) and (4.15). Let $N\to\infty$ in the above inequality. By the continuity of $\psi_g$ and the fact that $\psi_g(0)=0$, we see that the first term in the outer parentheses tends to zero. Note that $f_N(s)\longrightarrow\langle 1,\zeta^{\tilde{r}}(s)\rangle^{-1}$, as $N\to\infty$, for any $s\in[0,t)$ at which $\langle 1,\zeta^{\tilde{r}}(\cdot)\rangle^{-1}$ is continuous. Since it is continuous for almost every $s$ [the path $\zeta^{\tilde{r}}(\cdot)$ is r.c.l.l.], the second term in the outer parantheses tends to zero by (4.15) and bounded convergence. Furthermore, we note that

$$\sum_{j=0}^{N-1}\langle g'(\cdot)\tilde{h}_j, \zeta^{\tilde{r}}(t_j)\rangle = -\sum_{j=0}^{N-1}\langle g', \zeta^{\tilde{r}}(t_j)\rangle z_j\frac{t}{N},$$

and that as $N\to\infty$,

$$-\sum_{j=0}^{N-1}\langle g', \zeta^{\tilde{r}}(t_j)\rangle z_j\frac{t}{N} \longrightarrow -\int_0^t \frac{\langle g', \zeta^{\tilde{r}}(s)\rangle}{\langle 1, \zeta^{\tilde{r}}(s)\rangle}\,ds,$$

also by (4.15) and bounded convergence, since the function $\langle g',\zeta^{\tilde{r}}(\cdot)\rangle\langle 1,\zeta^{\tilde{r}}(\cdot)\rangle^{-1}$ is also continuous for almost every $s$. Together with the estimate (4.19), this implies that as $N\to\infty$,

$$\sum_{j=0}^{N-1}\langle g(\cdot - \bar{S}_{t_j,t^j}^{\tilde{r},m}) - g(\cdot), \zeta^{\tilde{r}}(t_j)\rangle \longrightarrow -\int_0^t \frac{\langle g', \zeta^{\tilde{r}}(s)\rangle}{\langle 1, \zeta^{\tilde{r}}(s)\rangle}\,ds. \tag{4.20}$$

We handle the first term of (4.12) in a similar (although simpler) fashion. Once again we can use a first-order Taylor expansion for each summand

appearing in this term:

$$g(v_i^{\tilde r} - \bar S^{\tilde r,m}_{U_i^{\tilde r}\tilde r^{-1}-m,t^j}) = g(v_i^{\tilde r}) + g'(w_j^i)h_j^i, \tag{4.21}$$

where $h_j^i = -\bar S^{\tilde r,m}_{U_i^{\tilde r}\tilde r^{-1}-m,t^j}$, and $w_j^i \in [v_i^{\tilde r} - \bar S^{\tilde r,m}_{U_i^{\tilde r}\tilde r^{-1}-m,t^j}, v_i^{\tilde r}]$. Since $|t^j - (U_i^{\tilde r}\tilde r^{-1} - m)| \le t/N$ for each pair $j,i$ in the first term of (4.12), we have by (2.20) and (4.15) as before that

$$\max_{j,i} |h_j^i| \le \frac{t}{N}\|\langle 1, \zeta^{\tilde r}(\cdot)\rangle^{-1}\|_L \le \frac{tM_\zeta}{N}. \tag{4.22}$$

Using the Taylor expansion (4.21) along with (4.22), we have

$$\begin{aligned}
&\left|\left(\sum_{j=0}^{N-1}\frac{1}{\tilde r}\sum_{i=\tilde r\bar E^{\tilde r,m}(t_j)+1}^{\tilde r\bar E^{\tilde r,m}(t^j)} g(v_i^{\tilde r} - \bar S^{\tilde r,m}_{U_i^{\tilde r}\tilde r^{-1}-m,t^j})\right) - \frac{1}{\tilde r}\sum_{i=\tilde r\bar E^{\tilde r,m}(0)+1}^{\tilde r\bar E^{\tilde r,m}(t)} g(v_i^{\tilde r})\right| \\
&\qquad = \left|\sum_{j=0}^{N-1}\frac{1}{\tilde r}\sum_{i=\tilde r\bar E^{\tilde r,m}(t_j)+1}^{\tilde r\bar E^{\tilde r,m}(t^j)} g'(w_j^i)h_j^i\right| \qquad (4.23)\\
&\qquad \le \bar E^{\tilde r,m}(t)\|g'\|_\infty \frac{tM_\zeta}{N},
\end{aligned}$$

which tends to zero as $N \to \infty$. By combining (4.20) and (4.23) above, we can let $N \to \infty$ in (4.12) to obtain the relation

$$\langle g, \zeta^{\tilde r}(t)\rangle = \langle g, \zeta^{\tilde r}(0)\rangle - \int_0^t \frac{\langle g', \zeta^{\tilde r}(s)\rangle}{\langle 1, \zeta^{\tilde r}(s)\rangle}\,ds + \frac{1}{\tilde r}\sum_{i=\tilde r\bar E^{\tilde r,m}(0)+1}^{\tilde r\bar E^{\tilde r,m}(t)} g(v_i^{\tilde r}), \tag{4.24}$$

for all sufficiently large $\tilde r$. We would like to let $\tilde r \to \infty$ in this relation to obtain

$$\langle g, \zeta(t)\rangle = \langle g, \zeta(0)\rangle - \int_0^t \frac{\langle g', \zeta(s)\rangle}{\langle 1, \zeta(s)\rangle}\,ds + \alpha t\langle g, \nu\rangle. \tag{4.25}$$

By (4.3), we have that the left-hand side, as well as the first term on the right-hand side of (4.24), converges to the corresponding term in (4.25). Similarly, (4.3), (4.13), (4.14) and (4.10) imply that the integrands in the second term on the right-hand side of (4.24) are uniformly bounded, and converge pointwise on $[0,t]$ to the integrand in the second term on the right-hand side of (4.25). Thus the integrals converge by bounded convergence. To see that the third term on the right-hand side converges, note that since $g \in \mathcal V$, $g \equiv g_k^{\mathcal V}$ for some $k \in \mathbb N$ (see Definition 3.7). So for any $\tilde r \in \tilde{\mathcal R}$ large enough so that $n(\tilde r) \ge k$,

$$\left|\frac{1}{\tilde r}\sum_{i=\tilde r\bar E^{\tilde r,m}(0)+1}^{\tilde r\bar E^{\tilde r,m}(t)} g_k^{\mathcal V}(v_i^{\tilde r}) - \alpha t\langle g_k^{\mathcal V}, \nu\rangle\right|$$

$$\le \left| \frac{1}{\tilde r} \sum_{i=\tilde r \bar E^{\tilde r,m}(0)+1}^{\tilde r \bar E^{\tilde r,m}(\lfloor t/l_{n(\tilde r)}\rfloor l_{n(\tilde r)})} g_k^{\mathcal V}(v_i^{\tilde r}) - \alpha \lfloor t/l_{n(\tilde r)}\rfloor l_{n(\tilde r)} \langle g_k^{\mathcal V}, \nu\rangle \right|$$

$$+ \left| \frac{1}{\tilde r} \sum_{i=\tilde r \bar E^{\tilde r,m}(\lfloor t/l_{n(\tilde r)}\rfloor l_{n(\tilde r)})+1}^{\tilde r \bar E^{\tilde r,m}(t)} g_k^{\mathcal V}(v_i^{\tilde r}) - \alpha (t - \lfloor t/l_{n(\tilde r)}\rfloor l_{n(\tilde r)}) \langle g_k^{\mathcal V}, \nu\rangle \right|$$

$$\le \varepsilon_{n(\tilde r)} + \|g_k^{\mathcal V}\|_\infty \Big( (\bar E^{\tilde r,m}(t) - \bar E^{\tilde r,m}(\lfloor t/l_{n(\tilde r)}\rfloor l_{n(\tilde r)})) + \alpha l_{n(\tilde r)} \Big) \tag{4.26}$$

$$\le \varepsilon_{n(\tilde r)} + \|g_k^{\mathcal V}\|_\infty \left( \sup_{t\in[0,L]} (\bar E^{\tilde r,m}(t + l_{n(\tilde r)}) - \bar E^{\tilde r,m}(t)) + \alpha l_{n(\tilde r)} \right)$$

$$\le \varepsilon_{n(\tilde r)} + \|g_k^{\mathcal V}\|_\infty \varepsilon_{n(\tilde r)},$$

where the second inequality is by (3.75), and the last inequality follows [after relaxing the bound to $\varepsilon_{n(\tilde r)}$] by (3.74) and (3.73). Since (4.26) holds for all $\tilde r \in \tilde{\mathcal R}$ sufficiently large, and since $\varepsilon_{n(\tilde r)} \downarrow 0$ as $\tilde r \to \infty$, we obtain (4.25) from (4.24) by letting $\tilde r \to \infty$.

We have shown that if $\zeta(0) \neq \mathbf{0}$, then (4.2) holds for all $t \in [0, L]$ and $g \in \mathcal V$. We now extend to the case $g \in \mathcal C$, using the case $g \in \tilde{\mathcal C}$ as an intermediate step. Suppose $g \in \tilde{\mathcal C}$ and choose functions $g_k \in \mathcal V$, $k \in \mathbb N$, such that $\{g_k\}_{k=1}^\infty$ and $\{g_k'\}_{k=1}^\infty$ are uniformly bounded, and such that as $k \to \infty$,

$$g_k \longrightarrow g \quad \text{and} \quad g_k' \longrightarrow g', \tag{4.27}$$

pointwise on $\mathbb R_+$ (see Definition 3.7). We have established that for each $k \in \mathbb N$ and $t \in [0, L]$,

$$\langle g_k, \zeta(t)\rangle = \langle g_k, \zeta(0)\rangle - \int_0^t \frac{\langle g_k', \zeta(s)\rangle}{\langle 1, \zeta(s)\rangle} ds + \alpha t \langle g_k, \nu\rangle. \tag{4.28}$$

It then follows by (4.27) and bounded convergence that, as $k \to \infty$, the left-hand side as well as the first and third terms on the right-hand side of (4.28) converge respectively to the corresponding terms of (4.2). Similarly, since $\zeta(\cdot)$ is continuous and $\zeta(t) \neq \mathbf{0}$ for all $t \in [0, L]$, the integral term also converges by (4.27) and two applications of bounded convergence. Thus, we obtain (4.2) from (4.28) by letting $k \to \infty$.

Finally, fix $g \in \mathcal C$. For $k \in \mathbb N$, choose a function $\psi_k \in \mathbf{C}_b^1(\mathbb R_+)$ such that $\psi_k(x) \in [0, 1]$ and $|\psi_k'(x)| \le 2$ for all $x \in \mathbb R_+$, $\psi_k \equiv 1$ on $[0, k]$ and $\psi_k \equiv 0$ on $[k+1, \infty)$. Let $g_k = \psi_k g$, and note that $\{g_k\}_{k=1}^\infty$ and $\{g_k'\}_{k=1}^\infty$ are uniformly bounded and that $g_k \longrightarrow g$, and $g_k' \longrightarrow g'$ are pointwise on $\mathbb R_+$ as $k \to \infty$. Since $g_k$ has compact support, we see that $g_k \in \tilde{\mathcal C}$ for each $k \in \mathbb N$. Therefore, (4.2) holds for all $t \in [0, L]$ and for each of the functions $g_k$. By the same argument as that appearing after (4.28), this implies that (4.2) also holds for $g$. This completes the proof. □

The next lemma asserts that for large $r$, the fluid scaled processes $\bar{\mu}^{r,m}(\cdot)$ can be uniformly approximated on $[0, L]$ by paths in $\mathscr{B}_L$. This fact, together with a result due to Puha and Williams [16] about the uniform convergence to steady state for elements of $\mathscr{B}_L$ (see Proposition 4.5) constitutes the second main ingredient for proving state space collapse.

LEMMA 4.4. *Assume* (Q.1) *and let* $T, L > 1$ *and* $0 < \eta < 1$ *be given. Let* $\{\varepsilon_n\}_{n=1}^{\infty} \subset (0, 1)$ *be a sequence such that* $\varepsilon_n \downarrow 0$, *as* $n \to \infty$. *Let* $l_n$, $p$, $M_T$, $\kappa_n$, $\gamma_n$, $N_n$, $r_n$ *be the constants, and* $\{B_n^r\}$ *be the events given by Lemma* 3.8. *Then for each* $n \in \mathbb{N}$, *there exists* $r'_n > r_n$ *such that for every* $r > r'_n$ *and* $\zeta^r(\cdot) \in \mathscr{B}_L^r$, *there exists* $\zeta(\cdot) \in \mathscr{B}_L$ *such that*

$$\|\mathbf{d}[\zeta^r(\cdot), \zeta(\cdot)]\|_L \le \varepsilon_n, \tag{4.29}$$

$$\|\langle \chi, \zeta^r(\cdot)\rangle - \langle \chi, \zeta(\cdot)\rangle\|_L \le \varepsilon_n. \tag{4.30}$$

PROOF. We follow the proof of ([2], Lemma 4.1). Suppose on the contrary that there exists $n \in \mathbb{N}$, a subsequence $\tilde{\mathcal{R}} \subset \mathcal{R}$, and a sequence $\{\zeta^{\tilde{r}}(\cdot)\}_{\tilde{r} \in \tilde{\mathcal{R}}} \subset \mathbf{D}_L(\mathcal{M}_{\mathrm{F}})$, such that $\zeta^{\tilde{r}}(\cdot) \in \mathscr{B}_L^{\tilde{r}}$ for each $\tilde{r}$, and such that either

$$\begin{aligned} &\inf_{\tilde{r} \in \tilde{\mathcal{R}}} \inf_{\zeta(\cdot) \in \mathscr{B}_L} \|\mathbf{d}[\zeta^{\tilde{r}}(\cdot), \zeta(\cdot)]\|_L > \varepsilon_n \quad \text{or} \\ &\inf_{\tilde{r} \in \tilde{\mathcal{R}}} \inf_{\zeta(\cdot) \in \mathscr{B}_L} \|\langle \chi, \zeta^{\tilde{r}}(\cdot)\rangle - \langle \chi, \zeta(\cdot)\rangle\|_L > \varepsilon_n. \end{aligned} \tag{4.31}$$

By Corollary 3.16 and Definition 4.1, there exists a $\zeta(\cdot) \in \mathscr{B}_L$ and a further subsequence $\{\tilde{r}_j\}_{j \in \mathbb{N}} \subset \tilde{\mathcal{R}}$ such that as $j \to \infty$,

$$\zeta^{\tilde{r}_j}(\cdot) \xrightarrow{J_1} \zeta(\cdot). \tag{4.32}$$

By Lemma 4.3, $\zeta(\cdot)$ is continuous. So as $j \to \infty$,

$$\|\mathbf{d}[\zeta^{\tilde{r}_j}(\cdot), \zeta(\cdot)]\|_L \to 0, \tag{4.33}$$

which implies by Lemma 3.5 and (3.76) that for each $t \in [0, L]$,

$$|\langle \chi, \zeta^{\tilde{r}_j}(t)\rangle - \langle \chi, \zeta(t)\rangle| \to 0 \qquad \text{as } j \to \infty. \tag{4.34}$$

Since $\langle \chi, \zeta(\cdot)\rangle$ is constant by (ii) of Lemma 4.3, we have by (4.34) and (3.81) that in fact,

$$\|\langle \chi, \zeta^{\tilde{r}_j}(\cdot)\rangle - \langle \chi, \zeta(\cdot)\rangle\|_L \to 0 \qquad \text{as } j \to \infty. \tag{4.35}$$

Together, (4.33) and (4.35) contradict (4.31). □

We now apply a result obtained in [16], on the rate at which fluid model solutions converge to their steady state, to choose $L > 1$ large enough so that elements of $\mathscr{B}_L$ are nearly in steady state for $t \in [L-1, L]$. We first recall the definition of the mapping $\Xi : \mathcal{M}_{\mathrm{F}}^c \longrightarrow \mathbf{D}_\infty(\mathcal{M}_{\mathrm{F}})$ introduced in [6], Lemma 4.9. Recall that the set $\mathcal{M}_{\mathrm{F}}^c \subset \mathcal{M}_{\mathrm{F}}$ consists of all nonatomic $\xi \in \mathcal{M}_{\mathrm{F}}$, that is, all $\xi \in \mathcal{M}_{\mathrm{F}}$ satisfying $\langle \mathbb{1}_{\{x\}}, \xi\rangle = 0$, for all $x \in \mathbb{R}_+$. Given $\xi \in \mathcal{M}_{\mathrm{F}}^c$, we define $\Xi(\xi)(\cdot) = \bar{\zeta}_\xi(\cdot)$, where $\bar{\zeta}_\xi(\cdot) \in \mathbf{D}_\infty(\mathcal{M}_{\mathrm{F}})$ is the unique fluid model solution for critical data $(\alpha, \nu)$ (defined on $[0, \infty)$ as in [6], Section 3.1), such that $\bar{\zeta}_\xi(0) = \xi$. The following proposition is a direct consequence of Theorem 1.3 in [16]. Recall that the lifting map $\Delta_\nu$ was defined in Definition 2.2.

PROPOSITION 4.5. *Assume* (Q.1) *and let* $T > 1$ *and* $0 < \eta < 1$ *be given. Let* $\{\varepsilon_n\}_{n=1}^\infty \subset (0, 1)$ *be a sequence such that* $\varepsilon_n \downarrow 0$, *as* $n \to \infty$. *Then for any* $\varepsilon > 0$, *there exists* $L^* > 1$ *such that for all* $\zeta(\cdot) \in \mathscr{B}_{L^*}$, $t \in [L^* - 1, L^*]$ *implies*

$$\mathbf{d}[\zeta(t), \Delta_\nu\langle \chi, \zeta(t)\rangle] \le \varepsilon. \tag{4.36}$$

PROOF. Let $l_n$, $p$, $M_T$, $\kappa_n$, $\gamma_n$, $N_n$, $r_n$ be the constants, and $\{B_n^r\}$ be the events given by Lemma 3.8. Let $q$ and $M_q$ be the constants given by Lemma 4.3, and define

$$\mathbf{S} = \{\xi \in \mathcal{M}_{\mathrm{F}}^c : \langle 1, \xi\rangle \vee \langle \chi, \xi\rangle \vee \langle \chi^{1+q}, \xi\rangle \le M_q\}. \tag{4.37}$$

Let $\mathscr{B} \subset \mathbf{D}_\infty(\mathcal{M}_{\mathrm{F}})$ be the set of all fluid model solutions $\bar{\zeta}(\cdot)$ for critical data $(\alpha, \nu)$ (defined on $[0, \infty)$ as in [6], Section 3.1) such that $\bar{\zeta}(0) \in \mathbf{S}$. Then as a direct consequence of Theorem 1.3(i) in [16] (with $M = M_q$ and $\varepsilon = q$ there), we have that given $\delta > 0$, there exists $L^* > 1$ such that

$$\sup_{\bar{\zeta}(\cdot) \in \mathscr{B}} \rho(\bar{\zeta}(t), \Delta_\nu\langle \chi, \bar{\zeta}(0)\rangle) \le \delta \qquad \text{for all } t \ge L^* - 1, \tag{4.38}$$

where $\rho(\cdot, \cdot)$ is the Prohorov metric on $\mathcal{M}_{\mathrm{F}}$ defined in [16], Section 1. Note that for any $\bar{\zeta}(\cdot) \in \mathscr{B}$, $\Delta_\nu\langle \chi, \bar{\zeta}(0)\rangle \in \mathbf{M}_\nu^q$, where

$$\mathbf{M}_\nu^q = \{\xi \in \mathbf{M}_\nu : \langle \chi, \xi\rangle \le M_q\}. \tag{4.39}$$

Thus, to replace $\delta$ and $\rho(\cdot, \cdot)$ with the given $\varepsilon$ and the metric $\mathbf{d}[\cdot, \cdot]$ in (4.38) above, it suffices to show that there exists $\delta > 0$ such that for all $\xi \in \mathbf{M}_\nu^q$, the $\rho$-ball $B_\rho(\xi, \delta)$ of radius $\delta$ centered at $\xi$, is contained in the $\mathbf{d}$-ball $B_{\mathbf{d}}(\xi, \varepsilon)$ of radius $\varepsilon$ centered at $\xi$. To this end, let $\xi^q = M_q\langle \chi, \nu_e\rangle^{-1}\nu_e$ be the element of $\mathbf{M}_\nu^q$ with the greatest mass. Note that the metrics $\rho(\cdot, \cdot)$ and $\mathbf{d}[\cdot, \cdot]$ both induce the same topology on $\mathcal{M}_{\mathrm{F}}$. So, we can choose $\delta > 0$ such that $B_\rho(\xi^q, \delta) \subset B_{\mathbf{d}}(\xi^q, \varepsilon)$. Suppose that $\xi \in \mathbf{M}_\nu^q$ and $\rho(\zeta, \xi) < \delta$ for some $\zeta \in \mathcal{M}_{\mathrm{F}}$. Since $\mathbf{M}_\nu^q = \{c\nu_e : c \in [0, M_q\langle \chi, \nu_e\rangle^{-1}]\}$, there is a $\xi' \in \mathbf{M}_\nu^q$ such that $\xi + \xi' = \xi^q$, where for two elements $\xi_1, \xi_2 \in \mathcal{M}_{\mathrm{F}}$, we define $\xi_1 + \xi_2 \in \mathcal{M}_{\mathrm{F}}$ by $(\xi_1 + \xi_2)(A) = \xi_1(A) + \xi_2(A)$, for any Borel set $A \subset \mathbb{R}_+$. Using the definition of the metric $\rho(\cdot, \cdot)$, it is not difficult to verify that

$$\rho(\zeta + \xi', \xi + \xi') \le \rho(\zeta, \xi) < \delta,$$

which implies that $\zeta + \xi' \in B_\rho(\xi^q, \delta) \subset B_{\mathbf{d}}(\xi^q, \varepsilon)$. By the definition (1.4) of $\mathbf{d}[\cdot,\cdot]$, this yields

$$\mathbf{d}[\zeta, \xi] = \mathbf{d}[\zeta + \xi', \xi + \xi'] < \varepsilon.$$

Thus, we can choose $L^* > 1$ so that

$$\sup_{\bar\zeta(\cdot)\in\mathscr{B}} \mathbf{d}[\bar\zeta(t), \Delta_\nu\langle\chi, \bar\zeta(0)\rangle] \le \varepsilon \qquad \text{for all } t \ge L^* - 1. \tag{4.40}$$

Let $\mathscr{B}_{L^*}$ be given by Definition 4.1 of this section. In [6], it was shown that for a fluid model solution $\bar\zeta(\cdot) \in \mathbf{D}_\infty(\mathcal{M}_{\mathrm{F}})$, $\bar\zeta(t)$ has no atoms for all $t \ge 0$ ([6], Lemma 4.3 and equation (4.33) ff.). This result carries over to fluid model solutions on $[0, L^*]$ in a straightforward manner. So by Lemma 4.3 of this section, we have for any $\zeta(\cdot) \in \mathscr{B}_{L^*}$, that $\zeta(t) \in \mathbf{S} \subset \mathcal{M}_{\mathrm{F}}^c$ for all $t \in [0, L^*]$. For any $\zeta(\cdot) \in \mathscr{B}_{L^*}$, let $\bar\zeta(\cdot) \in \mathbf{D}_\infty(\mathcal{M}_{\mathrm{F}})$ be defined by

$$\bar\zeta(t) = \begin{cases} \zeta(t), & t \in [0, L^*], \\ \Xi(\zeta(L^*))(t - L^*), & t \in (L^*, \infty). \end{cases} \tag{4.41}$$

Note that if $\zeta(0) = \mathbf{0}$, then by Lemma 4.3, Definition 4.2 and [6], Theorem 3.1, $\bar\zeta(t) = \mathbf{0}$ for all $t \ge 0$. If $\zeta(0) \ne \mathbf{0}$, then $\zeta(\cdot)$ satisfies (4.2) for all $t \in [0, L^*]$ and $\Xi(\zeta(L^*))(\cdot)$ satisfies (4.2) for all $t \ge 0$. It is not difficult to see that in this case $\bar\zeta(\cdot)$ must also satisfy (4.2) for all $t \ge 0$. So, $\bar\zeta(\cdot)$ is a fluid model solution for critical data $(\alpha, \nu)$ on $[0, \infty)$ such that $\bar\zeta(t) = \zeta(t)$ for all $t \in [0, L^*]$. Since $\zeta(0) \in \mathbf{S}$ for all $\zeta(\cdot) \in \mathscr{B}_{L^*}$, we have $\bar\zeta(\cdot) \in \mathscr{B}$ for all $\zeta(\cdot) \in \mathscr{B}_{L^*}$. Thus, by Lemma 4.3(ii) and (4.40), for any $t \in [L^* - 1, L^*]$,

$$\sup_{\zeta(\cdot)\in\mathscr{B}_{L^*}} \mathbf{d}[\zeta(t), \Delta_\nu\langle\chi, \zeta(t)\rangle] \le \varepsilon. \tag{4.42}$$

□

4.2. *State space collapse.* Before proceeding to the proof of state space collapse, we will need the following technical lemma, which provides a uniform continuity property for the mapping $\Xi : \mathcal{M}_{\mathrm{F}}^c \longrightarrow \mathbf{D}_\infty(\mathcal{M}_{\mathrm{F}})$ on a set of measures $\xi \in \mathcal{M}_{\mathrm{F}}^c$ that are close to the truncated invariant manifold $\mathbf{M}_\nu^q$. Recall that for any $\xi \in \mathcal{M}_{\mathrm{F}}^c$, $\bar\zeta_\xi(\cdot) = \Xi(\xi)(\cdot)$ is the unique fluid model solution defined on $[0, \infty)$ such that $\bar\zeta_\xi(0) = \xi$.

LEMMA 4.6. *Assume* (Q.1) *and let* $T, L > 1$ *and* $0 < \eta < 1$ *be given. Let* $\{\varepsilon_n\}_{n=1}^\infty \subset (0, 1)$ *be a sequence such that* $\varepsilon_n \downarrow 0$, *as* $n \to \infty$. *Let* $l_n$, $p$, $M_T$, $\kappa_n$, $\gamma_n$, $N_n$, $r_n$ *be the constants, and* $\{B_n^r\}$ *be the events given by Lemma* 3.8, *and let* $q$ *and* $M_q > 1$ *be the constants given by Lemma* 4.3. *Let* $M_1 \ge M_q$ *and define*

$$\mathbf{M}_\nu^1 = \{\xi \in \mathbf{M}_\nu : \langle\chi, \xi\rangle \le M_1\}. \tag{4.43}$$

*Then for any* $\varepsilon > 0$, *there exists* $\delta > 0$ *such that for all* $\zeta(\cdot) \in \mathscr{B}_L$ *and* $\xi \in \mathbf{M}_\nu^1$ *satisfying*

$$\mathbf{d}[\zeta(0), \xi] \vee |\langle\chi, \zeta(0)\rangle - \langle\chi, \xi\rangle| \le \delta, \tag{4.44}$$

*we have*

$$\sup_{t\in[0,L]} \mathbf{d}[\zeta(t),\bar{\zeta}_\xi(t)] \le \varepsilon. \tag{4.45}$$

PROOF. Fix $\varepsilon > 0$ and choose $n_0$ such that $\varepsilon_{n_0} \le \varepsilon/4$. Let

$$\delta_0 = \frac{\varepsilon(\kappa_{n_0} \wedge 1)(\langle \chi,\nu_e\rangle \wedge 1)}{8}. \tag{4.46}$$

Note that by (4.9), any $\zeta(\cdot) \in \mathscr{B}_L$ satisfies

$$\sup_{t\in[0,L]} \langle \mathbb{1}_{[0,\kappa_{n_0})}, \zeta(t)\rangle \le \frac{\varepsilon_{n_0}}{2} \le \frac{\varepsilon}{8}. \tag{4.47}$$

Note also that for any $\zeta(\cdot) \in \mathscr{B}_L$, $\zeta(t) = \Xi(\zeta(0))(t)$ for all $t \in [0,L]$. By Theorem 3.8 in [15], $\Xi$ is continuous on $\mathcal{M}_{\mathrm{F}}^c$. Let $\xi_1 = M_1\langle \chi,\nu_e\rangle^{-1}\nu_e$. Since $\xi_1 \in \mathcal{M}_{\mathrm{F}}^c$, we can choose $0 < \delta_1 < 1$ such that for any $\xi \in \mathcal{M}_{\mathrm{F}}^c$ satisfying $\mathbf{d}[\xi,\xi_1] \le \delta_1$, we have

$$\sup_{t\in[0,2M_1L/\delta_0]} \mathbf{d}[\bar{\zeta}_\xi(t),\bar{\zeta}_{\xi_1}(t)] \le \varepsilon. \tag{4.48}$$

Define

$$\delta = \frac{\delta_0\delta_1}{2M_1}, \tag{4.49}$$

and consider any $\zeta(\cdot) \in \mathscr{B}_L$ such that for some $\xi \in \mathbf{M}_\nu^1$, (4.44) holds. The argument splits into two cases. First, suppose that $\langle \chi,\zeta(0)\rangle < \delta_0$. Then, by (4.44), (4.49) and the fact that $M_1 > 1$, $\langle \chi,\xi\rangle \le 2\delta_0$. This implies by (4.46) that $\langle 1,\xi\rangle = \langle \chi,\xi\rangle\langle \chi,\nu_e\rangle^{-1} \le \varepsilon/4$. Thus,

$$\begin{aligned}
&\sup_{t\in[0,L]} \mathbf{d}[\zeta(t),\bar{\zeta}_\xi(t)] \\
&\qquad \le \sup_{t\in[0,L]} \mathbf{d}[\zeta(t),\mathbf{0}] + \mathbf{d}[\mathbf{0},\xi] \\
&\qquad \le \sup_{t\in[0,L]} 2\langle 1,\zeta(t)\rangle + 2\langle 1,\xi\rangle \\
&\qquad \le 2\sup_{t\in[0,L]} (\langle \mathbb{1}_{[0,\kappa_{n_0})},\zeta(t)\rangle + \langle \mathbb{1}_{(\kappa_{n_0},\infty)},\zeta(t)\rangle) + \frac{\varepsilon}{2} \\
&\qquad \le \frac{\varepsilon}{4} + \frac{2}{\kappa_{n_0}}\langle \chi,\zeta(t)\rangle + \frac{\varepsilon}{2} \\
&\qquad \le \varepsilon,
\end{aligned} \tag{4.50}$$

where the first inequality uses the fact that $\bar{\zeta}_\xi(t) = \bar{\zeta}_\xi(0) = \xi$ for all $t \ge 0$, since $\xi$ is on the invariant manifold ([16], Theorem 1.1). The second inequality uses definition (1.4), the fourth uses (4.47) and Markov's inequality, and the last inequality uses (4.46) and the fact that since $\zeta(\cdot)$ is a fluid model solution, $\langle\chi,\zeta(t)\rangle = \langle\chi,\zeta(0)\rangle < \delta_0$ for all $t \in [0,L]$ ([6], Theorem 3.1). This establishes (4.45) in the first case.

For the second case, suppose that $\langle\chi,\zeta(0)\rangle \geq \delta_0$. Since $\delta_1 < 1 < M_1$, (4.44) and (4.49) imply that $\langle\chi,\xi\rangle \geq \delta_0/2$. So $\xi \neq \mathbf{0}$. Let $c_1 = M_1\langle\chi,\xi\rangle^{-1}$ and define the transformed function $\tilde{\zeta}(\cdot) : [0, c_1 L] \longrightarrow \mathcal{M}_{\mathrm{F}}$ by

$$\tilde{\zeta}(t) = c_1\zeta(t/c_1), \qquad t \in [0, c_1 L]. \tag{4.51}$$

Using a change of variables, it is not difficult to verify that $\tilde{\zeta}(\cdot)$ satisfies (4.2) for all $t \in [0, c_1L]$ and is therefore a fluid model solution on $[0, c_1L]$. Note that $\tilde{\zeta}(0) = c_1\zeta(0)$, and that $c_1\xi = \xi_1$. Also note that since $\xi \in \mathbf{M}^1_\nu$, we have $c_1 \geq 1$. Thus, using definition (1.4), we have

$$\begin{aligned}\mathbf{d}[\tilde{\zeta}(0),\xi_1] &= \mathbf{d}[c_1\zeta(0), c_1\xi]\\ &\leq c_1\,\mathbf{d}[\zeta(0),\xi]\\ &\leq \frac{2M_1}{\delta_0}\delta\\ &= \delta_1.\end{aligned} \tag{4.52}$$

This implies that

$$\begin{aligned}\sup_{t\in[0,L]} \mathbf{d}[\zeta(t),\bar{\zeta}_\xi(t)] &\leq \sup_{t\in[0,L]} \mathbf{d}[c_1\zeta(t), c_1\xi]\\ &= \sup_{t\in[0,c_1L]} \mathbf{d}[\tilde{\zeta}(t),\bar{\zeta}_{\xi_1}(t)]\\ &\leq \varepsilon.\end{aligned} \tag{4.53}$$

The first line in the above inequality uses (1.4), the fact that $c_1 \geq 1$ and the fact that $\bar{\zeta}_\xi(t) = \xi$, since $\xi \in \mathbf{M}^1_\nu$. The second line follows by definition of $\tilde{\zeta}(\cdot)$, and the fact that $c_1\xi = \xi_1 \in \mathbf{M}^1_\nu$. The last line then follows by (4.52) and (4.48) [with $\xi = \tilde{\zeta}(0)$ there] by noting that $c_1L \leq 2M_1L/\delta_0$. $\square$

We are now ready to prove state space collapse for the sequence of diffusion scaled state descriptors $\{\hat{\mu}^r(\cdot)\}$. This leads directly to the proof of Theorem 2.3, which appears at the end of the section.

THEOREM 4.7 (State space collapse). *Assume* (Q.1) *and let* $T > 1$ *be given. Then as* $r \to \infty$,

$$\|\mathbf{d}[\hat{\mu}^r(\cdot), \Delta_\nu\hat{W}^r(\cdot)]\|_T \Longrightarrow 0. \tag{4.54}$$

PROOF. We must show that for any $0 < \varepsilon^*, \eta < 1$, there exists $r^* \in \mathcal{R}$ such that $r > r^*$ implies

$$\mathbf{P}^r(\|\mathbf{d}[\hat{\mu}^r(\cdot), \Delta_\nu\hat{W}^r(\cdot)]\|_T \leq \varepsilon^*) \geq 1-\eta. \tag{4.55}$$

Fix $0 < \varepsilon^*, \eta < 1$. We first note that by definition of the metric $\mathbf{d}[\cdot,\cdot]$ and the lifting map $\Delta_\nu$,

$$\mathbf{d}[\Delta_\nu w_1, \Delta_\nu w_2] \le c_\nu |w_1 - w_2|, \qquad w_1, w_2 \in \mathbb{R}_+, \tag{4.56}$$

for some constant $c_\nu \ge 1$. Choose a sequence $\{\varepsilon_n\}_{n=1}^\infty$ such that $\varepsilon_n \downarrow 0$ as $n \to \infty$. By Proposition 4.5, using $\varepsilon = \varepsilon^*/3$ there, we can choose $L^* > 1$ such that for all $\zeta(\cdot) \in \mathscr{B}_{L^*}$, $t \in [L^* - 1, L^*]$ implies

$$\mathbf{d}[\zeta(t), \Delta_\nu \langle \chi, \zeta(t)\rangle] \le \frac{\varepsilon^*}{3}. \tag{4.57}$$

Let $l_n$, $p$, $M_T$, $\kappa_n$, $\gamma_n$, $N_n$, $r_n$ be the constants, and $\{B_n^r\}$ be the events given by Lemma 3.8. Let $q$ and $M_q > 1$ be the constants given by Lemma 4.3, and let $\delta > 0$ be given by Lemma 4.6 for $L = L^*$, $M_1 = 2M_q$ and $\varepsilon = \varepsilon^*/12$. Fix $n^* \in \mathbb{N}$ large enough so that

$$\varepsilon_{n^*} \le \frac{\delta}{2} \wedge \frac{\varepsilon^*}{12 c_\nu}. \tag{4.58}$$

Choose $r^*$ large enough so that $n(r^*) \ge n^*$ and $r^* \ge r'_{n^*}$, where $r'_{n^*} > r_{n^*}$ is given by Lemma 4.4. Fix $r > r^*$ for the remainder of the proof. The argument splits into two cases. We first consider the process $\bar{\mu}^{r,0}(\cdot)$ on $B_{n(r)}^r$, and then consider the processes $\bar{\mu}^{r,m}(\cdot)$ for $m \in \{1, \ldots, \lfloor rT \rfloor\}$ on $B_{n(r)}^r$.

For the first case, consider any realization $\zeta^r(\cdot) \in \mathscr{B}_{L^*}^r$ of the process $\bar{\mu}^{r,0}(\cdot)$ on the event $B_{n(r)}^r$. By (3.77), there exists a $\xi \in \mathbf{M}_\nu$ such that

$$\mathbf{d}[\zeta^r(0), \xi] \vee |\langle \chi, \zeta^r(0)\rangle - \langle \chi, \xi\rangle| < \varepsilon_{n(r)} \le \varepsilon_{n^*}. \tag{4.59}$$

By (3.76) and the fact that $\varepsilon_{n^*} < 1 < M_T \le M_q < M_1$, we see that $\xi \in \mathbf{M}_\nu^1 = \{\xi \in \mathbf{M}_\nu : \langle \chi, \xi\rangle \le 2M_q\}$. Let $\bar{\zeta}_\xi(\cdot) = \Xi(\xi)(\cdot)$ and note that since $\xi$ is on the invariant manifold, $\xi = \bar{\zeta}_\xi(0) = \bar{\zeta}_\xi(t)$ for all $t \in [0, L^*]$ ([16], Theorem 1.1). So by definition of $\mathbf{M}_\nu$ and $\Delta_\nu$,

$$\bar{\zeta}_\xi(t) = \Delta_\nu \langle \chi, \bar{\zeta}_\xi(t)\rangle \qquad \text{for all } t \in [0, L^*]. \tag{4.60}$$

Let $\zeta(\cdot) \in \mathscr{B}_{L^*}$ be given by Lemma 4.4 for $\zeta^r(\cdot)$ and $n = n^*$. Note that since $\bar{\zeta}_\xi(\cdot)$ and $\zeta(\cdot)$ are fluid model solutions on $[0, L^*]$, we have for all $t \in [0, L^*]$ ([6], Theorem 3.1),

$$\langle \chi, \bar{\zeta}_\xi(t)\rangle = \langle \chi, \bar{\zeta}_\xi(0)\rangle \quad \text{and} \quad \langle \chi, \zeta(t)\rangle = \langle \chi, \zeta(0)\rangle. \tag{4.61}$$

By combining the result of Lemma 4.4 cited above with (4.59) and (4.58), one obtains

$$\begin{aligned} \mathbf{d}[\zeta(0), \xi] &\le \mathbf{d}[\zeta(0), \zeta^r(0)] + \mathbf{d}[\zeta^r(0), \xi] \\ &\le 2\varepsilon_{n^*} \\ &\le \delta, \end{aligned} \tag{4.62}$$

and also

$$
\begin{aligned}
&|\langle \chi, \bar{\zeta}_\xi(0)\rangle - \langle \chi, \zeta(0)\rangle| \\
&\qquad \le |\langle \chi, \bar{\zeta}_\xi(0)\rangle - \langle \chi, \zeta^r(0)\rangle| + |\langle \chi, \zeta^r(0)\rangle - \langle \chi, \zeta(0)\rangle| \\
&\qquad \le 2\varepsilon_{n^*} \\
&\qquad \le \delta.
\end{aligned} \tag{4.63}
$$

This implies, by choice of $\delta$ and Lemma 4.6, that for any $t \in [0, L^*]$,

$$
\mathbf{d}[\zeta(t), \bar{\zeta}_\xi(t)] \le \frac{\varepsilon^*}{12}. \tag{4.64}
$$

Since $r > r^*$, we have for any $t \in [0, L^*]$,

$$
\begin{aligned}
&\mathbf{d}[\zeta^r(t), \Delta_\nu\langle \chi, \zeta^r(t)\rangle] \\
&\quad \le \mathbf{d}[\zeta^r(t), \zeta(t)] + \mathbf{d}[\zeta(t), \bar{\zeta}_\xi(t)] + \mathbf{d}[\bar{\zeta}_\xi(t), \Delta_\nu\langle \chi, \bar{\zeta}_\xi(t)\rangle] \\
&\qquad + \mathbf{d}[\Delta_\nu\langle \chi, \bar{\zeta}_\xi(t)\rangle, \Delta_\nu\langle \chi, \zeta(t)\rangle] + \mathbf{d}[\Delta_\nu\langle \chi, \zeta(t)\rangle, \Delta_\nu\langle \chi, \zeta^r(t)\rangle] \\
&\quad \le \varepsilon_{n^*} + \frac{\varepsilon^*}{12} + 0 + c_\nu 2\varepsilon_{n^*} + c_\nu \varepsilon_{n^*} \\
&\quad \le \frac{\varepsilon^*}{2},
\end{aligned} \tag{4.65}
$$

where the first estimate in the second inequality above follows by Lemma 4.4, the second estimate follows by (4.64), the third by (4.60), the fourth by (4.56), (4.61) and (4.63), and the fifth estimate uses (4.56) and Lemma 4.4. The last inequality follows by (4.58) and the fact that $c_\nu \ge 1$.

We now proceed to the second case. Let $\zeta^r(\cdot) \in \mathscr{B}^r_{L^*}$ be a realization of the process $\bar{\mu}^{r,m}(\cdot)$ on $B^r_{n(r)}$, where $m \in \{1, \dots, \lfloor rT \rfloor\}$. Let $\zeta(\cdot) \in \mathscr{B}_{L^*}$ be given by Lemma 4.4 for $\zeta^r(\cdot)$ and $n = n^*$. Since $r > r^*$, we have for any $t \in [L^* - 1, L^*]$,

$$
\begin{aligned}
&\mathbf{d}[\zeta^r(t), \Delta_\nu\langle \chi, \zeta^r(t)\rangle] \\
&\quad \le \mathbf{d}[\zeta^r(t), \zeta(t)] + \mathbf{d}[\zeta(t), \Delta_\nu\langle \chi, \zeta(t)\rangle] \\
&\qquad + \mathbf{d}[\Delta_\nu\langle \chi, \zeta(t)\rangle, \Delta_\nu\langle \chi, \zeta^r(t)\rangle] \\
&\quad \le \varepsilon_{n^*} + \frac{\varepsilon^*}{3} + c_\nu \varepsilon_{n^*} \\
&\quad \le \frac{\varepsilon^*}{2},
\end{aligned} \tag{4.66}
$$

where the three estimates in the second inequality are by Lemma 4.4, (4.57), (4.56) and Lemma 4.4, and the last inequality is by (4.58) and the fact that $c_\nu \ge 1$.

Finally, combining the estimates (4.65) and (4.66), we have on $B^r_{n(r)}$, for $r > r^*$,

$$
\begin{aligned}
&\|\mathbf{d}[\hat{\mu}^r(\cdot), \Delta_\nu \hat{W}^r(\cdot)]\|_T \\
&\quad = \sup_{t\in[0,rT]} \mathbf{d}[\bar{\mu}^r(t), \Delta_\nu\langle\chi, \bar{\mu}^r(t)\rangle] \\
&\quad \le \sup_{t\in[0,L^*]} \mathbf{d}[\bar{\mu}^{r,0}(t), \Delta_\nu\langle\chi, \bar{\mu}^{r,0}(t)\rangle] \\
&\qquad + \sup_{1\le m\le\lfloor rT\rfloor,\ t\in[L^*-1,L^*]} \mathbf{d}[\bar{\mu}^{r,m}(t), \Delta_\nu\langle\chi, \bar{\mu}^{r,m}(t)\rangle] \\
&\quad \le \frac{\varepsilon^*}{2} + \frac{\varepsilon^*}{2}.
\end{aligned} \tag{4.67}
$$

Since $\mathbf{P}^r(B^r_{n(r)}) \ge 1 - \eta$ for $r > r^*$, (4.55) is proved. □

PROOF OF THEOREM 2.3. By Proposition 3.1, we have that $\hat{W}^r(\cdot) \Longrightarrow W^*(\cdot)$ as $r \to \infty$, where $W^*(\cdot)$ is a reflected Brownian motion in $\mathbb{R}_+$ with drift $-\lambda$, variance $\alpha a^2 + \beta b^2$ and initial condition $W^*(0)$ equal in distribution to $\langle\chi, \Theta\rangle$. Since $\Delta_\nu : \mathbb{R}_+ \longrightarrow \mathcal{M}_{\mathrm{F}}$ is continuous, the continuous mapping theorem implies that $\Delta_\nu \hat{W}^r(\cdot) \Longrightarrow \mu^*(\cdot) = \Delta_\nu W^*(\cdot)$ as $r \to \infty$. Thus, Theorem 4.7, combined with the "converging together lemma" ([1], Theorem 4.1), implies that $\hat{\mu}^r(\cdot) \Longrightarrow \mu^*(\cdot)$ as $r \to \infty$. □

**Acknowledgments.** The author thanks Ruth Williams, Amber Puha, Pat Fitzsimmons and Łukasz Kruk for helpful discussions and comments during the preparation of this work. Many thanks are extended to the referee as well, for helpful suggestions and considerable effort.

EURANDOM
5612 AZ Eindhoven
The Netherlands
e-mail: gromoll@eurandom.tue.nl